\documentclass[12pt,a4paper]{article}

\usepackage{amsfonts}
\usepackage{amsmath}
\usepackage{amsbsy}
\usepackage{theorem}
\usepackage{graphicx}

\newtheorem{theorem}{Theorem}
\newtheorem{proposition}{Proposition}
\newtheorem{lemma}{Lemma}
\newtheorem{corollary}{Corollary}
\newtheorem{remark}{Remark}

\begin{document}
\sf

\title{\bfseries Testing LRD in the spectral domain for  functional time series in manifolds
}
\date{}

 \maketitle

\author{\hspace*{-1cm} M.D. Ruiz-Medina$^{1}$ and Rosa M. Crujeiras$^{2}$\\
$^{1}$ University of Granada\\
$^{2}$ University of Santiago de Compostela}

\begin{abstract}
A statistical hypothesis test  for long range dependence (LRD) is formulated  in the spectral domain for  functional time series  in manifolds.
 The elements of the spectral density operator family are assumed to be  invariant with respect to the group of isometries of the manifold.
The proposed  test statistic   is based on the weighted periodogram operator. A Central Limit Theorem is   derived  to obtain the asymptotic Gaussian distribution of the proposed test statistic operator under the null hypothesis. The rate of convergence to zero, in the  Hilbert--Schmidt operator norm, of the bias of the integrated  empirical second and fourth order cumulant  spectral  density  operators  is obtained under the alternative hypothesis. The consistency of the test follows from  the consistency of the  integrated  weighted periodogram  operator  under LRD.   Practical implementation of our     testing  approach is based on the random projection methodology.
 A simulation study illustrates, in the context of spherical functional time series,  the  asymptotic normality of the test statistic  under the null hypothesis, and its consistency under the alternative. The empirical size and power properties are also computed for different functional sample sizes, and under different scenarios.
\end{abstract}

\medskip

\noindent \textbf{MSC2020 subject classifications:}

\noindent {\em Primary.} 60G10, 60G12, 60G18 , 60G20, 60G22

\noindent {\em Secondary.} 60G60

\medskip


\section{Introduction}

\label{intro}
Spherical functional time series analysis helps in understanding the dynamics and spatiotemporal patterns of data that are embedded into the sphere, providing valuable insights for prediction, monitoring, and decision-making.
 Time series analysis of   global temperature data distribution among other climate variables, usually arising in Climate Science and Meteorology,
 can be performed in a more  efficient way by adopting a  functional time series framework
 (see   \cite{Shah2024}). That is the  case of ocean currents, and other marine functional time series
      to be analyzed in Oceanography studies (see, e.g.,  \cite{Tirellietal24}; \cite{Yarger}). Other areas  demanding  this type of techniques   are   Geophysics, Astronomy and Astrophysics. In the last  few decades, the cosmic microwave background radiation variation analysis over time has gained special attention (see   \cite{Duke24}; \cite{Marinucci04}; \cite{Marinucci}).  In a more general manifold setting, functional time series analysis is often  applied  in
    Medical Imaging, Computer Vision and Graphics  (see  \cite{Wu2024};  \cite{Zhou2013};  \cite{Zhu2009}, among others).
   This paper focuses on the spectral analysis of   functional time series in manifolds, with special attention to LRD analysis.

The spectral analysis of functional time series has  mainly been developed under Short Range Dependence (SRD).  In this context, based on the weighted periodogram operator, a nonparametric framework  is adopted in  \cite{Panaretos13}. Particularly, the asymptotic normality of the functional discrete Fourier transform (fDFT), and the weighted periodogram operator of the curve data are proved under suitable summability conditions on the $L^{2}$ norm of the cumulant spectral density operators. The consistency of the weighted periodogram operator, in  the Hilbert--Schmidt operator norm,  is  derived under SRD.  In \cite{Panaretos13b}, a  harmonic principal component analysis  of functional time series    in the temporal functional spectral domain is also obtained, based on a  Karhunen–Lo\'eve--like decomposition,   the so–called Cram\'er–Karhunen–Lo\'eve representation.   In the context of functional regression, some  applications are presented  in \cite{Pham2018}, \cite{Rubin20a} and  \cite{Tavakoli2014}. Hypothesis testing for detecting modelling differences in functional time series dynamics  is addressed in \cite{Tavakoli2016} in the functional spectral domain.

In   LRD  analysis of functional time series  several problems   still remain open. One of the key approaches in the current literature
is presented in \cite{LiRobinsonShang19}, where
  the
eigendecomposition of the
 long-run covariance operator is considered,  under an  asymptotic semiparametric  functional principal component framework.
   The consistent estimation of the dimension and the orthonormal
functions  spanning the dominant subspace, where the projected curve
process displays the largest dependence range is derived. Fractionally integrated functional autoregressive moving averages processes  constitute an  interesting example of this modelling framework.

 A first attempt to characterize LRD in functional time series in the  spectral domain can be found  in \cite{RuizMedina2022}, adopting
  the theoretical framework of operator--valued random fields, including fractional Brownian motion with   operator--valued Hurst coefficient
(see, e.g.,  \cite{duker18};  \cite{MarinucciRV}; \cite{Rackauskas}  and  \cite{Rackauskasv2}).  The  eigenvalues of the LRD operator are parameterized. These eigenvalues  induce different levels of singularity at zero frequency, corresponding to different levels of temporal persistent of the  process  projected into different eigenspaces of the Laplace Beltrami operator. Under this LRD scenario, the integrated periodogram operator
 is proved to be asymptotically  unbiased in the Hilbert--Schmidt operator norm. Minimum contrast estimation of the LRD operator is achieved in the  spectral domain in a weak--consistent way under a Gaussian scenario. Interesting examples of this setting are analyzed in
 \cite{Ovalle23}, where    the spectral analysis of multifractionally integrated  functional time series in manifolds is considered. In particular, multifractionally integrated spherical functional ARMA models (i.e., multifractionally integrated SPHARMA models) are analyzed through simulations. In this modelling framework, SRD and LRD can coexist at different spherical scales.

 Up to our knowledge,  no further  developments have been achieved in the  spectral analysis of LRD functional time series.
 Alternative  contributions for stationary  LRD functional sequences are based on the diagonalization of the heavy tail autocovariance kernel of the time--varying functional error term  (see, e.g., \cite{Beran16a}).  However, under this modelling framework, functional spectral  analysis can not be achieved in the time domain due to the assumed independence between the random components of the error term  (see also \cite{Beran14}). Similar assertions hold for the two sample problem analyzed in
 \cite{Beran16b}. In \cite{LiRobinsonShang21}, a special family of LRD linear functional time series is analyzed  with scalar coefficients displaying slow decay. Under stationarity, this LRD scenario  constitutes a particular case  of our framework when the elements of the spectral density operator family have  degenerated pure point spectra corresponding to  one infinite--dimensional eigenspace. The Cram\'er–Karhunen–Lo\'eve representation above referred constitutes a powerful tool in the functional spectral analysis of weak-dependent functional time series (see \cite{Rubin20b}). In this paper, this representation is extended to our LRD stationary functional time series  context, assuming the  invariance of the  cross-covariance kernels under the group of isometries of the manifold, given by a connected and compact two points homogeneous space.

  In this paper we perform a weighted periodogram operator based analysis, requiring the asymptotic analysis of the bias of  the integrated empirical fourth--order   cumulant spectral density operators, to prove consistency  of the integrated weighted periodogram operator   under LRD.  Its application to spectral  statistical hypothesis testing of LRD in   $L^{2}\left(\mathbb{M}_{d},d\nu , \mathbb{R}\right)$--valued correlated sequences constitutes one of the main goals of this work. Here,  $L^{2}\left(\mathbb{M}_{d},d\nu , \mathbb{R}\right)$ denotes the space   of real--valued square integrable functions  on a Riemannian manifold $\mathbb{M}_{d},$  embedded into $\mathbb{R}^{d+1},$ given by a connected and compact  two--point homogeneous space. The topological dimension of $\mathbb{M}_{d}$ is
  $d,$ and $d\nu$  denotes the normalized Riemannian measure on $\mathbb{M}_{d}.$  In what follows, we will consider $X=\left\{X_{t},\ t\in \mathbb{Z}\right\}$ to be a  functional sequence such that    $\mathcal{P}\left(X_{t}\in L^{2}\left(\mathbb{M}_{d},d\nu , \mathbb{R}\right)\right)=1,$ for every $t\in \mathbb{Z},$ with $\mathcal{P}$ denoting the probability measure defined on the basic probability space $(\Omega ,\mathcal{Q}, \mathcal{P}),$ i.e., for every $t\in \mathbb{Z},$
\begin{equation}X_{t}: (\Omega ,\mathcal{Q}, \mathcal{P})\longrightarrow L^{2}\left(\mathbb{M}_{d},d\nu , \mathbb{R}\right)
\label{eq1}
\end{equation}
\noindent is a measurable mapping.

The   invariance of the elements of the spectral  density operator family of $X$ under the group of isometries of the manifold $\mathbb{M}_{d}$ is assumed along the paper. A   frequency--varying eigenvalue sequence then characterizes  the pure point spectra of the elements of the  spectral density operator family, with respect to the orthonormal basis of  eigenfunctions of the Laplace--Beltrami operator.  This  invariance assumption is exploited in
 the derived Central Limit Theorem  that characterizes the asymptotic distribution of the proposed  test statistic operator  under the   null hypothesis, which states that $X$ displays SRD. In our formulation of the alternative hypothesis on LRD, we adopt a semiparametric framework in terms of a functional parameter given by the LRD operator.  In contrast with the approach presented in  \cite{RuizMedina2022}, here we do not assume a parameterization of the eigenvalues of the LRD operator.  Under this scenario,   the rate of convergence to zero,  in the corresponding $L^{2}$ norm, of the the bias of the  integrated empirical second and fourth order cumulant spectral density operators  is respectively obtained in  Lemmas \ref{lem1} and \ref{lem4cs} under suitable conditions.  Proposition \ref{cor1}  shows the divergence, in   the Hilbert--Schmidt operator norm, of the mean of the  test statistic operator under the alternative hypothesis. Theorem \ref{consisimse}  derives suitable conditions, in particular,  on the  nonparametric functional  spectral factor,  to ensure consistency of the integrated weighted periodogram operator in the Hilbert--Schmidt operator norm under LRD.   Theorem \ref{th2} then provides the almost surely divergence of the  test statistic in  the Hilbert--Schmidt operator norm under the alternative, yielding the consistency of the test.

Theorem \ref{th2} also plays a crucial role in the implementation  in practice of the proposed  testing procedure, based on  rejecting the null hypothesis when the random Fourier coefficients  of the test operator statistic, suitably standardized according to
 Theorem \ref{lemvfh0}, cross an upper or lower tail standard normal critical value.  The  orthonormal basis involved in the computation of these coefficients  is constructed by the tensor product of the eigenfunctions of the Laplace Beltrami operator. The random projection methodology   (see Theorem 4.1 in \cite{Cuestalbertos}) can be implemented here to  alleviate the dimensionality problem. Specifically, when the moments of our test statistic operator  satisfy the Carleman condition under the null hypothesis,  our testing procedure is equivalent to rejecting the null hypothesis when  the absolute value of a random projection of the test statistic is larger than  an   upper tail standard normal critical value.  In the implementation of the  random projection  methodology, the Karhunen--Lo\'eve expansion  in Lemma \ref{CLTSBT} below can   be considered  for generation of the  involved Gaussian random directions. In the simulation study undertaken, robust empirical sizes, and competitive values of the empirical power are displayed by our testing approach (see Section \ref{sesptp}).

The outline of the paper is as follows. In Section \ref{background}, the functional spectral background material is introduced. Our hypothesis   testing procedure  is formulated in a  functional  semiparametric  spectral framework in Section  \ref{HT}. The asymptotic Gaussian distribution of the test statistic operator under the null hypothesis $H_{0}$ is obtained in Theorem \ref{lemvfh0} in Section \ref{underH0}.  Asymptotics of the bias of   the  integrated empirical second and fourth order  cumulant spectral density operators   under LRD are derived in Section \ref{underalternative}. Section \ref{tlrdfts} provides the preliminary results required for  consistency of the test, which is derived in  Theorem \ref{th2} of this section. Practical implementation is also discussed in Section \ref{tlrdfts}. In Section \ref{sgd},   a simulation study is undertaken   to  illustrate the asymptotic Gaussian distribution  of the proposed test statistic operator under the null hypothesis,  in
the context of SRD spherical functional time series. The consistency of the test is also illustrated in Section \ref{consec},  in the framework of  multifractionally integrated spherical functional time series.  This numerical analysis is extended in Section \ref{Sec5.3} to a wider family of LRD operators allowing stronger persistency in time, displayed by the projected process in the dominant subspace.
Section \ref{sesptp} analyzes empirical size and power properties of the  test. Section \ref{s6} summarizes conclusions of the simulation study,
 focusing on the large  functional sample size properties  of our test statistic operator under different bandwidth parameter scenarios, from additional numerical results. The proofs of the   results  of this paper can be found  in the Appendix.

\subsection{Background}
\label{background}

 Along this work we will assume that  $X=\left\{X_{t},\ t\in \mathbb{Z}\right\}$ in (\ref{eq1}) is
a   stationary zero--mean   functional sequence, with nuclear covariance operator  family  $\left\{ \mathcal{R}_{\tau},\ \tau\in \mathbb{Z}\right\}$ satisfying $\mathcal{R}_{\tau} =E[X_{s}\otimes X_{s+\tau}]=E[ X_{s+ \tau}\otimes X_{s}],$ for every $s,\tau \in \mathbb{Z}.$
   The elements of this family are  characterized in the spectral domain by
the spectral density operator family $\{ \mathcal{F}_{\omega },\ \omega \in [-\pi,\pi]\}.$ The assumed invariance of the elements of these families with respect to  the group of isometries of $\mathbb{M}_{d}$   yields to  their  diagonal  series expansion in terms  of   $\{S_{n,j}^{d}\otimes \overline{S_{n,j}^{d}}, \ j=1,\dots, \Gamma (n,d),\ n\in \mathbb{N}_{0}\},$ with  $\{S_{n,j}^{d}, \ j=1,\dots, \Gamma (n,d),\ n\in \mathbb{N}_{0}\}$ being
 the  orthonormal basis of   eigenfunctions of the Laplace--Beltrami operator $\Delta_{d}$ on $L^{2}\left(\mathbb{M}_{d},d\nu , \mathbb{C}\right)$ (see, e.g., \cite{Gine75}; \cite{Helgason59}).   In particular,
\begin{eqnarray}& &
  \mathcal{F}_{\omega}
\underset{\mathcal{S}(L^{2}(\mathbb{M}_{d},d\nu; \mathbb{C}))}{=}  \frac{1}{2\pi} \sum_{\tau \in \mathbb{Z}}\exp\left(-i\omega \tau\right)\mathcal{R}_{\tau}\nonumber\\
 & &\underset{\mathcal{S}(L^{2}(\mathbb{M}_{d},d\nu; \mathbb{C}))}{=}\sum_{n\in \mathbb{N}_{0}}f_{n}(\omega)\sum_{j=1}^{\Gamma (n,d)}S_{n,j}^{d}\otimes \overline{S_{n,j}^{d}},\quad \omega \in [-\pi,\pi],
  \label{sdo2}
\end{eqnarray}
\noindent where, for every $n\in \mathbb{N}_{0},$ $\Gamma (n,d)$ denotes the dimension of the eigenspace $\mathcal{H}_{n}$ associated with the eigenvalue $\lambda_{n}(\Delta_{d})$ of the Laplace Beltrami operator $\Delta_{d}$  (see, e.g., Section 2.1 in \cite{MaMalyarenko}).
 The equality  $ \underset{\mathcal{S}(L^{2}(\mathbb{M}_{d},d\nu; \mathbb{C}))}{=}$ means identity in the norm of the space of Hilbert--Schmidt operators on   \linebreak   $L^{2}(\mathbb{M}_{d},d\nu; \mathbb{C}),$  the space of complex--valued square integrable functions on $\mathbb{M}_{d}.$   Specifically, the equality in (\ref{sdo2})  means that
 $$\int_{\mathbb{M}_{d}\times \mathbb{M}_{d}}\left|f_{\omega}(x,y)-\sum_{n=0}^{\infty}f_{n}(\omega)\sum_{j=1}^{\Gamma (n,d)}S_{n,j}^{d}(x)\overline{S_{n,j}^{d}(y)}\right|^{2}d\nu(x)d\nu(y)=0,$$
 \noindent where $f_{\omega}$ is the kernel of the integral operator $\mathcal{F}_{\omega},$ for every $\omega\in [-\pi,\pi].$

Let $\left\{ X_{t},\ t=0,\dots,T-1\right\}$ be a functional sample of size $T\geq 2$ of $X.$ The    fDFT $\widetilde{X}_{\omega }^{(T)}$ is  defined as
\begin{eqnarray}\widetilde{X}_{\omega }^{(T)}(x)&=&\frac{1}{\sqrt{2\pi T}}\sum_{t=0}^{T-1}X_{t}(x)\exp(-i\omega t),\quad x\in \mathbb{M}_{d},\quad \omega \in [-\pi,\pi].\label{fDFT}
\end{eqnarray}
\noindent   The  kernel  $p_{\omega }^{(T)}$ of the  periodogram operator $\mathcal{P}^{(T)}_{\omega }=\widetilde{X}_{\omega}^{(T)}\otimes  \widetilde{X}^{(T)}_{-\omega}$  satisfies, for each $\omega \in [-\pi,\pi],$
 \begin{eqnarray}
  p_{\omega }^{(T)}(x,y)&=&\frac{1}{2\pi T}\sum_{t=0}^{T-1}\sum_{s=0}^{T-1}X_{t}(x)X_{s}(y)\exp(-i\omega (t-s)),\ \forall  x,y\in \mathbb{M}_{d}. \label{periodogro}
\end{eqnarray}
\noindent We will  denote  by $f_{\omega }^{(T)}(x,y)=\mbox{cum}\left(\widetilde{X}^{(T)}_{\omega}(x), \widetilde{X}^{(T)}_{-\omega}(y)\right)=E\left[p^{(T)}_{\omega }(x,y)\right],$ $x,y\in \mathbb{M}_{d},$ the kernel of the cumulant  operator $\mathcal{F}_{\omega }^{(T)}$ of order $2$ of the fDFT  $\widetilde{X}_{\omega }^{(T)}$   over the diagonal $\omega \in [-\pi,\pi].$
Note that,  for
$\omega \in [-\pi,\pi],$ and $T\geq 2,$  the F\'ejer kernel is given by
\begin{equation}
F_{T}(\omega )=\frac{1}{T}\sum_{t=1}^{T}
\sum_{s=1}^{T}\exp\left(-i(t-s)\omega \right)=\frac{1}{T}\left[\frac{\sin\left(T\omega/2\right)}{\sin(\omega/2)}\right]^{2}.
\label{eqfkd}
\end{equation}

 The weighted periodogram operator, denoted as $\widehat{\mathcal{F}}_{\omega }^{(T)},$  has kernel $\widehat{f}_{\omega }^{(T)}(x,y)$ given,
 for each  $\omega \in [-\pi,\pi],$ by
 \begin{eqnarray}
&&\widehat{f}_{\omega }^{(T)}(x,y)=\left[\frac{2\pi}{T}\right]\sum_{s=1}^{T-1}  W^{(T)}\left(\omega - \frac{2\pi s}{T} \right)
 p_{ \frac{2\pi  s}{T} }^{(T)}(x,y) ,\ x,y\in \mathbb{M}_{d}, \label{enp}
\end{eqnarray}
\noindent where $W^{(T)}$  is a  weight function  satisfying
\begin{equation}\label{Eq5}
W^{(T)}(x) = \sum_{j\in \mathbb{Z}}\frac{1}{B_{T}} W\left(\frac{x + 2\pi j}{B_{T}}\right),
\end{equation}
\noindent with $B_{T}$ being the positive  bandwidth parameter. Function $W$ is  a real--valued function defined on $\mathbb{R}$ such that
  $W$ is positive, even, and bounded in variation; $W(x) =0$, if  $ |x|\geq 1$;
 $\int_{\mathbb{R}} \left|W(x)\right|^{2}dx <\infty;$  $\int_{\mathbb{R}} W(x)dx =1.$

\subsection{Hypothesis   testing}
\label{HT}
The    SRD and LRD scenarios  respectively
tested under  the null  $H_{0}$ and the alternative $H_{1}$ hypotheses are introduced in this section.  The   proposed  test statistic operator based on the weighted periodogram operator   is then  formulated.

  Stationary  SRD functional time series are characterized by the summability of the series of trace norms of the elements of the  family of covariance operators  $\left\{ \mathcal{R}_{\tau},\ \tau\in \mathbb{Z}\right\}$  (see, e.g.,  \cite{Panaretos13}). That is, $X$ displays SRD if and only if
  $\sum_{\tau \in \mathbb{Z}}\|\mathcal{R}_{\tau}\|_{L^{1}(L^{2}(\mathbb{M}_{d},d\nu , \mathbb{R}))}<\infty,$  where $L^{1}(L^{2}(\mathbb{M}_{d},d\nu , \mathbb{R}))$ denotes the space of trace operators on $L^{2}(\mathbb{M}_{d},d\nu , \mathbb{R}).$ In our setting, this condition can be formulated as follows:
\begin{equation}\sum_{\tau \in \mathbb{Z}}\|\mathcal{R}_{\tau}\|_{L^{1}(L^{2}(\mathbb{M}_{d},d\nu , \mathbb{R}))}=\sum_{\tau \in \mathbb{Z}}\sum_{n\in \mathbb{N}_{0}}\Gamma (n,d)\left|\int_{-\pi}^{\pi}\exp\left(i\omega \tau \right)f_{n}(\omega )d\omega \right|<\infty.\label{sm}
\end{equation}

When (\ref{sm}) fails,   $X$ is said to display LRD.
In  what follows, we will adopt the LRD  scenario introduced in   \cite{RuizMedina2022}, given by
 \begin{equation}
\mathcal{F}_{\omega }=\mathcal{M}_{\omega }|\omega|^{-\mathcal{A}},\quad \omega \in [-\pi,\pi],\label{LRD}
\end{equation}
\noindent where  the invariant positive self--adjoint operators  $\mathcal{M}_{\omega }$ and $|\omega|^{-\mathcal{A}}$  are composed yielding the definition of $\mathcal{F}_{\omega }.$
Specifically, $\mathcal{A}$ denotes the     LRD operator   on   $L^{2}(\mathbb{M}_{d},d\nu; \mathbb{C}).$ Operator $|\omega|^{-\mathcal{A}}$ in (\ref{LRD})  is interpreted as in \cite{Charac14}, \cite{Rackauskas} and   \cite{Rackauskasv2},
where $\mathcal{A}$ plays the role of  operator--valued Hurst coefficient in the setting of  fractional Brownian motion introduced in this framework.
Moreover,
$\mathcal{M}_{\omega }$ is the regular  spectral operator  reflecting markovianess when the null space of $\mathcal{A}$ coincides with $L^{2}(\mathbb{M}_{d},d\nu; \mathbb{C}).$ Specifically, $\mathcal{M}_{\omega }$  satisfies\begin{equation}\sum_{\tau \in \mathbb{Z}}\left\|\int_{[-\pi,\pi]}\exp(i\omega \tau )\mathcal{M}_{\omega }  d\omega \right\|_{L^{1}(L^{2}(\mathbb{M}_{d},d\nu , \mathbb{R}))}<\infty,\label{regpart}
\end{equation}
 \noindent where
the operator integrals are understood as improper operator Stieltjes integrals which converge strongly (see, e.g., Section 8.2.1 in  \cite{Ramm05}).

We will apply the spectral theory of self--adjoint operators (see, e.g., \cite{Dautray85}) in terms of the  common spectral kernel
$$\Upsilon (x,y)=\sum_{n\in \mathbb{N}_{0}}\sum_{j=1}^{\Gamma (n,d)}
S_{n,j}^{d}(x)\overline{S_{n,j}^{d}(y)},\quad x,y\in \mathbb{M}_{d},$$
\noindent under the assumed invariance property with respect to the group of isometries of $\mathbb{M}_{d}.$

The point spectrum of  $\mathcal{A}$ is given by $\left\{ \alpha (n),\ n\in \mathbb{N}_{0}\right\},$  with
$l_{\alpha }\leq \alpha (n)\leq L_{\alpha },$  for every $n\in \mathbb{N}_{0},$  and  $l_{\alpha }, L_{\alpha }\in (0,1/2).$
 It is assumed that LRD operator $\mathcal{A}$ has kernel $\mathcal{K}_{\mathcal{A}}$ admitting the following series expansion in the weak--sense:
 \begin{equation}
\mathcal{K}_{\mathcal{A}}(x,y)=\sum_{n\in \mathbb{N}_{0}}\alpha (n)\sum_{j=1}^{\Gamma (n,d)}S_{n,j}^{d}\otimes \overline{S_{n,j}^{d}}(x,y).
 \label{lerdopbb}
 \end{equation}
 \noindent  Specifically, identity (\ref{lerdopbb}) is understood  as
 \begin{equation}
 \mathcal{A}(f)(g)=\int_{\mathbb{M}_{d}\times \mathbb{M}_{d}}f(x)g(y)\sum_{n\in \mathbb{N}_{0}}\alpha (n)\sum_{j=1}^{\Gamma (n,d)}S_{n,j}^{d}(x) \overline{S_{n,j}^{d}}(y)d\nu(x)d\nu(y),
 \label{lerdop}
 \end{equation}
\noindent for all $f,g \in C^{\infty }(\mathbb{M}_{d}),$  where $C^{\infty }(\mathbb{M}_{d})$ denotes
 the space of infinitely differentiable functions with compact support contained  in $\mathbb{M}_{d}.$
Note that, under the conditions assumed,  $\mathcal{A}$ and $\mathcal{A}^{-1}$ are bounded, and $\left\|\mathcal{A}\right\|_{\mathcal{L}(L^{2}(\mathbb{M}_{d}, d\nu,\mathbb{C}))}<1/2,$ with $\left\|\cdot \right\|_{\mathcal{L}(L^{2}(\mathbb{M}_{d}, d\nu,\mathbb{C}))}$ denoting the norm in the space
$\mathcal{L}(L^{2}(\mathbb{M}_{d}, d\nu,\mathbb{C}))$ of bounded linear operators on $L^{2}(\mathbb{M}_{d}, d\nu,\mathbb{C}).$

 In a similar way,    operator  $|\omega|^{-\mathcal{A}}$   is interpreted as
\begin{equation}|\omega|^{-\mathcal{A}}(f)(g)=\int_{\mathbb{M}_{d}\times \mathbb{M}_{d}}f(x)g(y)\sum_{n\in \mathbb{N}_{0}}\frac{1}{|\omega |^{\alpha (n)}}\sum_{j=1}^{\Gamma (n,d)}S_{n,j}^{d}\otimes \overline{S_{n,j}^{d}}(x, y),
  d\nu(x)d\nu(y),\label{lerdop0}
 \end{equation}
\noindent  for every $f,g \in C^{\infty }(\mathbb{M}_{d})$ and $\omega\in [-\pi,\pi]\backslash \{0\}.$

 Operator  $\mathcal{M}_{\omega }$ in (\ref{LRD}) is  a   Hilbert--Schmidt operator on $L^{2}(\mathbb{M}_{d},d\nu; \mathbb{C}),$ whose kernel $\mathcal{K}_{\mathcal{M}_{\omega }}(x,y)$ admits the following series expansion in the norm of the space $\mathcal{S}(L^{2}(\mathbb{M}_{d},d\nu; \mathbb{C})):$
\begin{equation}
\mathcal{K}_{\mathcal{M}_{\omega }}(x,y)=\sum_{n\in \mathbb{N}_{0}}M_{n}(\omega )\sum_{j=1}^{\Gamma (n,d)}S_{n,j}^{d}\otimes \overline{S_{n,j}^{d}}(x, y),\quad x,y\in \mathbb{M}_{d},\ \omega \in [-\pi,\pi],
\label{rp}
\end{equation}
\noindent
where   $\{M_{n}(\omega ),\ n\in \mathbb{N}_{0}\}$ denotes the sequence of positive eigenvalues. For each $n\in \mathbb{N}_{0},$ $M_{n}(\omega ),$ $\omega \in [-\pi,\pi],$   is a  continuous positive slowly varying  function at
$\omega =0$ in the Zygmund's sense (see Definition 6.6 in \cite{Beran17}, and  Assumption IV in   \cite{RuizMedina2022}). Equation (\ref{regpart}) can  be equivalently expressed, in terms of   $\left\{M_{n}(\omega ),\ n\in \mathbb{N}_{0},\ \omega \in [-\pi,\pi]\right\},$  as \begin{equation}\sum_{\tau \in \mathbb{Z}}\sum_{n\in \mathbb{N}_{0}}\Gamma (n,d)\left|\int_{-\pi}^{\pi}\exp\left(i\omega \tau \right)M_{n}(\omega )d\omega \right|<\infty.\label{tracep}
 \end{equation}
\noindent  As before,  equation (\ref{tracep}) implies that   $X$ displays SRD, when $\alpha (n)=0,$ for every $n\in \mathbb{N}_{0}.$ Under (\ref{tracep}), $\left\{ \mathcal{M}_{\omega },\ \omega \in [-\pi,\pi]\right\}$ is also included  in the trace class.

Under the above setting of conditions,
 \begin{eqnarray}&&
\int_{-\pi}^{\pi}\left\|\mathcal{F}_{\omega }\right\|_{\mathcal{S}(L^{2}(\mathbb{M}_{d},d\nu, \mathbb{C}))}^{2}d\omega <\infty, \label{shconvrm}\end{eqnarray}
\noindent i.e., $\left\|\mathcal{F}_{\omega }\right\|_{\mathcal{S}(L^{2}(\mathbb{M}_{d},d\nu, \mathbb{C}))}\in L^{2}([-\pi,\pi]),$ with $L^{2}([-\pi,\pi])$ being the space of square integrable functions on the interval $[-\pi,\pi].$
Condition  (\ref{shconvrm})  plays a crucial role in the derivation of the results of this paper under LRD.

From equations (\ref{LRD})--(\ref{rp}),  the elements of the positive frequency varying eigenvalue sequence $\left\{f_{n}(\omega),\ n\in \mathbb{N}_{0}\right\}$ in (\ref{sdo2})   admit the following expression:
\begin{eqnarray}  f_{n}(\omega )&=&\frac{M_{n}(\omega )}{|\omega |^{\alpha (n)}},\quad \omega \in [-\pi,\pi],\ n\in \mathbb{N}_{0}.
  \label{eqscm1}
 \end{eqnarray}

 Note that,  since  $\sin(\omega )\sim \omega,$
$\omega \to 0,$
\begin{equation}
\left|1-\exp\left(-i\omega \right)\right|^{-\mathcal{A}}=[4\sin^{2}(\omega /2)]^{-\mathcal{A}/2}\sim |\omega
|^{-\mathcal{A}},\quad \omega \to 0.\label{ffsvint}
\end{equation}

 Sequence (\ref{eqscm1}) is involved in the formulation of the alternative hypothesis $H_{1}$ stating that $X$ displays LRD against  $H_{0}$ where SRD is assumed. Specifically,
 \begin{eqnarray}
 && H_{0}: \ f_{n}(\omega )=M_{n}(\omega ), \  \omega \in [-\pi,\pi],\ \forall n\in \mathbb{N}_{0}\label{altfhp0}\\
 && H_{1}: \ f_{n}(\omega )=M_{n}(\omega )\left|\omega\right|^{-\alpha (n)},\ \omega \in [-\pi,\pi],\ \forall   n\in \mathbb{N}_{0}.
 \label{reatlrd}
 \end{eqnarray}

\noindent    In our context, the formulation of the test statistic must capture the singularities at zero frequency for different manifold resolution levels under $H_{1}.$     The proposed  test  statistic operator $\mathcal{S}_{B_{T}}$ is then given by
\begin{equation}\mathcal{S}_{B_{T}}=\sqrt{B_{T}T}\int_{[-\sqrt{B_{T}}/2, \sqrt{B_{T}}/2]}\widehat{\mathcal{F}}_{\omega }^{(T)}
\frac{d\omega}{\sqrt{B_{T}}} ,\label{oslrdtintro}
\end{equation}
\noindent where the kernel of the integral operator $\widehat{\mathcal{F}}_{\omega }^{(T)}$ has been introduced in equation (\ref{enp}), with, as before,   $B_{T}$ being   the bandwidth parameter. The indicator function on the interval $[-\sqrt{B_{T}}/2, \sqrt{B_{T}}/2]$ is denoted by  $\mathbb{I}_{[-\sqrt{B_{T}}/2, \sqrt{B_{T}}/2]}(\omega ),$ for $\omega \in [-\pi,\pi].$  Note that,  as $T\to \infty,$  $$\int_{-\pi}^{\pi}\frac{\mathbb{I}_{[-\sqrt{B_{T}}/2, \sqrt{B_{T}}/2]}(\omega )}{\sqrt{B_{T}}}h(\omega  )d\omega \to \int_{-\pi}^{\pi}\delta (0-\omega )h(\omega  )d\omega = h(0),$$ \noindent for every $h\in L^{2}([-\pi ,\pi])$  (see \cite{Gelfand64} for the usual notion of  convergence in the sense of generalized functions or distributions).   Here,  $\delta (0-\omega )$ denotes  the  Dirac Delta distribution
at zero frequency. Hence, in what follows  we adopt the notation $\delta_{T}(0-\omega )=\frac{\mathbb{I}_{[-\sqrt{B_{T}}/2, \sqrt{B_{T}}/2]}(\omega )}{\sqrt{B_{T}}}$ representing a truncated Dirac Delta distribution.

\section{Preliminary results under SRD}
\label{underH0}
The following lemma will be applied in  the proof of the main result of this section,   Theorem \ref{lemvfh0},   deriving the asymptotic Gaussian distribution of $\mathcal{S}_{B_{T}}$ in (\ref{oslrdtintro})
under $H_{0}.$ Specifically,
Lemma \ref{lem2pt}    provides  the   asymptotic Gaussian distribution of the weighted periodogram operator $\widehat{\mathcal{F}}_{\omega }^{(T)}$ under $H_{0}.$  Its proof can be obtained in the same way as in \cite{Panaretos13}, where this result is formulated for the separable Hilbert space $H=L^{2}([0,1],\mathbb{C}).$

 \begin{lemma}
 \label{lem2pt}
 Assume that $E\|X_{0}\|^{k}<\infty,$ for all $k\geq 2,$ and
 \begin{itemize}
 \item[(i)] $\sum_{t_{1},\dots,t_{k-1}\in \mathbb{Z}}\|\mbox{cum}\left(X_{t_{1}},\dots,X_{t_{k-1}},X_{0}\right)\|_{L^{2}(\mathbb{M}^{k}_{d}, \otimes_{i=1}^{k} d\nu (x_{i}),\mathbb{R})}<\infty,$ $k\geq 2$
 \item[(i$^{\prime}$)]  $\sum_{t_{1},\dots,t_{k-1}\in \mathbb{Z}}(1+|t_{j}|)\|\mbox{cum}\left(X_{t_{1}},\dots,X_{t_{k-1}},X_{0}\right)\|_{L^{2}(\mathbb{M}^{k}_{d}, \otimes_{i=1}^{k} d\nu (x_{i}),\mathbb{R})}<\infty,$ \\ for
 $k\in\{2,4\},$ $j<k$
 \item[(ii)] $\sum_{t\in \mathbb{Z}}(1+|t|)\|\mathcal{R}_{t}\|_{L^{1}(L^{2}(\mathbb{M}_{d},  d\nu ,\mathbb{R}))}<\infty$
 \item[(iii)] $\sum_{t_{1},t_{2},t_{3}\in \mathbb{Z}}\|\mathcal{R}_{t_{1},t_{2},t_{3}}\|_{L^{1}(L^{2}(\mathbb{M}_{d}^{3},  \otimes_{i=1}^{3} d\nu (x_{i}) ,\mathbb{R}))}<\infty.$
 \end{itemize}
 Then, for every frequencies $\omega_{j},$ $j=1,\dots,J,$  with $J<\infty,$
\begin{equation}
\sqrt{B_{T}T}(\widehat{f}_{\omega_{j} }^{(T)}-E[\widehat{f}_{\omega_{j} }^{(T)}])\to_{D} \widehat{f}_{\omega_{j} },\ j=1,\dots,J\label{NCLT}
\end{equation}
\noindent where $\to_{D}$ denotes the convergence in distribution. Here,    $\widehat{f}_{\omega_{j} },$ $j=1,\dots,J,$ are jointly zero--mean complex Gaussian elements in $\mathcal{S}(L^{2}(\mathbb{M}_{d},d\nu, \mathbb{C}))$  $=L^{2}(\mathbb{M}_{d}^{2},d\nu\otimes d\nu , \mathbb{C}),$ with covariance kernel:
\begin{eqnarray}
&&\mbox{cov}(\widehat{f}_{\omega_{i}}(x_{1},y_{1}),\widehat{f}_{\omega_{j}}(x_{2},y_{2}))=2\pi \|W\|_{L^{2}(\mathbb{R})}^{2}\left\{
\eta (\omega_{i}-\omega_{j})f_{\omega_{i}}(x_{1},x_{2})f_{-\omega_{i}}(y_{1},y_{2})\right.
\nonumber\\
&&
\left.+\eta (\omega_{i}+\omega_{j})f_{\omega_{i}}(x_{1},y_{2})f_{-\omega_{i}}(y_{1},x_{2})\right\},\quad (x_{i},y_{i})\in \mathbb{M}_{d}^{2},\ i=1,2,
\label{cobk}
\end{eqnarray}
\noindent with $\eta (\omega )=1,$ for $\omega \in 2\pi \mathbb{Z},$ and  $\eta (\omega )=0,$ otherwise.
Thus, $\widehat{f}_{\omega_{i}}$ and $\widehat{f}_{\omega_{j}}$ are independent for $\omega_{i}+\omega_{j}\neq 0,$ mod $2\pi $
and $\omega_{i}-\omega_{j}\neq 0,$ mod $2\pi .$  For zero frequency modulus $2\pi$ the limit Gaussian random element is in $\mathcal{S}(L^{2}(\mathbb{M}_{d},d\nu, \mathbb{R}))=L^{2}(\mathbb{M}_{d}^{2},d\nu\otimes d\nu, \mathbb{R}).$

\end{lemma}

\medskip

\noindent  \textbf{Proof}. See Theorem 3.7  in  \cite{Panaretos13}.

 \medskip

The next result provides the asymptotic Gaussian distribution of the test statistic operator  $\mathcal{S}_{B_{T}}$  under $H_{0}.$ The convergence to a Gaussian random element in  the norm of the space
 $\mathcal{L}^{2}_{\mathcal{S}(L^{2}(\mathbb{M}_{d},d\nu, \mathbb{C}))}(\Omega ,\mathcal{A},P)$ is also obtained. Here,
 $\mathcal{L}_{\mathcal{S}(L^{2}(\mathbb{M}_{d},d\nu; \mathbb{C}))}(\Omega ,\mathcal{A},\mathcal{P})$ denotes the space of zero--mean second--order  \linebreak  $\mathcal{S}(L^{2}(\mathbb{M}_{d},d\nu; \mathbb{C}))$--valued random variables with the norm $\sqrt{E\|\cdot\|^{2}_{\mathcal{S}(L^{2}(\mathbb{M}_{d},d\nu; \mathbb{C}))}}.$
\begin{theorem}
\label{lemvfh0}
Under $H_{0},$ assume that  the conditions of Lemma \ref{lem2pt} hold.  Then,
     \begin{equation}\mathcal{S}_{B_{T}}-E[\mathcal{S}_{B_{T}}]\to_{D} Y_{0}^{(\infty)},\quad T\to \infty,\label{lemcdrso}\end{equation}
\noindent where $\mathcal{S}_{B_{T}}$ has been introduced in (\ref{oslrdtintro}), and $Y_{0}^{(\infty)}$ is a zero--mean  Gaussian random element
in  $\mathcal{S}(L^{2}(\mathbb{M}_{d},d\nu, \mathbb{R})),$ with autocovariance operator  $\mathcal{R}_{Y_{0}^{(\infty)}}= E\left[ Y_{0}^{(\infty)}\otimes Y_{0}^{(\infty)}\right]$ having kernel introduced in  equation (\ref{cobk}) in Lemma
\ref{lem2pt}, with $\omega_{i}=\omega_{j}=0.$
\end{theorem}

\medskip

\noindent \textbf{Proof}. See Appendix \ref{apendicea1}.

\section{Second and fourth order bias  asymptotics under LRD}

\label{underalternative}
 This section provides new results on the bias asymptotics in the Hilbert-Schmidt operator norm of the
integrated   empirical second and fourth order   cumulant spectral density operators  of $X$ under $H_{1}.$ These results are applied
  in the derivation of Theorem  \ref{consisimse} and Corollary \ref{cor22}, providing the consistency of the integrated weighted periodogram operator under $H_{1}.$
  In what follows,  we assume  $B_{T}\to 0$ and $B_{T}T\to \infty$ as $T\to \infty.$

  The rate of convergence to zero of the norm of the  bias in the space \linebreak $\mathcal{S}(L^{2}(\mathbb{M}_{d},d\nu, \mathbb{C}))$   of the integrated periodogram operator    is obtained  under LRD in the next lemma.
  The following well-known  identity will be applied: \begin{eqnarray}&&\mathcal{F}_{\omega}^{(T)}=E\left[\mathcal{P}_{ \omega }^{(T)} \right]=\left[F_{T}*\mathcal{F}_{\bullet
}\right](\omega )=
\int_{-\pi}^{\pi} F_{T}(\omega - \xi)
\mathcal{F}_{\xi} d\xi,\quad T\geq 2,
\nonumber\\
\label{cfk}
\end{eqnarray}
\noindent for  $\omega \in [-\pi,\pi]\backslash \{0\},$   where  $F_{T}(\omega )$ denotes the F\'ejer kernel  introduced in equation
 (\ref{eqfkd}) of Section \ref{background}.

\begin{lemma}
\label{lem1}
 Under  $H_{1},$ as   $T\to \infty,$

\begin{eqnarray}&&
\int_{-\pi}^{\pi}\mathcal{F}^{(T)}_{\omega }d\omega  =\int_{-\pi}^{\pi}E_{H_{1}}\left[\mathcal{P}_{ \omega }^{(T)} \right]d\omega
\underset{\mathcal{S}(L^{2}(\mathbb{M}_{d},d\nu, \mathbb{C}))}{=}
\int_{-\pi}^{\pi}\mathcal{F}_{\omega }d\omega +\mathcal{O}(T^{-1}),\nonumber\\
 \label{orderofconvmpo}
   \end{eqnarray}

 \noindent where $E_{H_{1}}$ denotes expectation under the  alternative $H_{1},$ and, as before,  $\underset{\mathcal{S}(L^{2}(\mathbb{M}_{d},d\nu, \mathbb{C}))}{=}$ denotes the equality in the norm of the space
$\mathcal{S}(L^{2}(\mathbb{M}_{d},d\nu, \mathbb{C})).$

\end{lemma}

\medskip

\noindent \textbf{Proof}. See Appendix \ref{ap2}.

\medskip

The following corollary is obtained from Lemma \ref{lem1}, and   provides the rate of convergence to zero of the bias of the integrated weighted periodogram operator, in the norm of  the space $\mathcal{S}(L^{2}(\mathbb{M}_{d},d\nu, \mathbb{C}))$ under $H_{1}.$
\begin{corollary}\label{lem3}
Under $H_{1},$  as $T\to \infty,$ \begin{equation}\int_{-\pi}^{\pi}E_{H_{1}}[\widehat{\mathcal{F}}^{(T)}_{\omega }]d\omega \underset{\mathcal{S}(L^{2}(\mathbb{M}_{d},d\nu, \mathbb{C}))}{=}\int_{-\pi}^{\pi}\int_{\mathbb{R}}W(\xi)\mathcal{F}_{\omega -\xi B_{T}}d\xi d\omega+\mathcal{O}(B_{T}^{-1}T^{-1})+\mathcal{O}(T^{-1}).\label{eqapphh}\end{equation}
\end{corollary}

\medskip

\noindent \textbf{Proof}. See Appendix  \ref{aplem3}.

\medskip

The  rate of convergence to zero, in the  norm of the space \linebreak $\mathcal{S}\left(L^{2}\left(\mathbb{M}_{d}^{2},\otimes_{i=1}^{2}\nu(dx_{i}), \mathbb{C}\right)\right)\equiv L^{2}\left(\mathbb{M}_{d}^{4},\otimes_{i=1}^{4}\nu(dx_{i}), \mathbb{C}\right),$ of the bias of the  integrated  empirical fourth--order   cumulant spectral density operators of $X$ under LRD
 is obtained in Lemma \ref{lem4cs} below. The following assumption is required:

 \medskip

 \noindent \emph{Assumption I}.   For every $t_{1},t_{2},t_{3}\in \mathbb{Z},$ $\mbox{cum}\left(X_{t_{1}}, X_{t_{2}}, X_{t_{3}},  X_{0}\right)$ defines an isotropic kernel in $L^{2}(\mathbb{M}_{d}^{4},\otimes_{i=1}^{4}d\nu (dx_{i}), \mathbb{R}),$ and the following convergence holds:

\medskip

\begin{eqnarray}&&
   \sum_{t_{1},t_{2},t_{3}\in
\mathbb{Z}}\left\|\mbox{cum}\left(X_{t_{1}}, X_{t_{2}}, X_{t_{3}},  X_{0}\right)\right\|^{2}_{L^{2}(\mathbb{M}^{4}_{d}, \otimes_{i=1}^{4} d\nu (x_{i}),\mathbb{R})}<\infty,
\label{eqfinitthirdcumulant}
\end{eqnarray}
\noindent where
\begin{eqnarray}&&\left\|\mbox{cum}\left(X_{t_{1}}, X_{t_{2}}, X_{t_{3}},  X_{0}\right)\right\|_{L^{2}(\mathbb{M}_{d}^{4},\otimes_{i=1}^{4}d\nu (dx_{i}), \mathbb{R})}^{2}\nonumber\\
&&=\int_{\mathbb{M}_{d}^{4}}
\left|\mbox{cum}\left(X_{t_{1}}(x), X_{t_{2}}(y), X_{t_{3}}(z),  X_{0}(v)\right)\right|^{2}d\nu (x)d\nu (y)d\nu (z)d \nu(v).\nonumber
\end{eqnarray}

\medskip

\begin{lemma}
\label{lem4cs}
Under $H_{1}$  and  \emph{Assumption I}, uniformly in $\omega_{4}\in [-\pi,\pi],$
\begin{eqnarray}&&
\int_{[-\pi ,\pi ]^{3}}
T\mbox{cum}\left(\widetilde{X}_{\omega_{1}}^{(T)}(\tau_{1}), \widetilde{X}_{\omega_{2}}^{(T)}(\tau_{2}),
\widetilde{X}_{\omega_{3}}^{(T)}(\tau_{3}), \widetilde{X}_{\omega_{4}}^{(T)}(\tau_{4})\right)
d\omega_{1}d\omega_{2}d\omega_{3}\nonumber\\
&&
\underset{\mathcal{S}\left(L^{2}\left(\mathbb{M}_{d}^{2},\otimes_{i=1}^{2}\nu(dx_{i}), \mathbb{C}\right)\right)}{=}2\pi \int_{[-\pi ,\pi ]^{3}}\mathcal{F}_{\omega_{1},\omega_{2},\omega_{3}}(\tau_{1},\tau_{2}, \tau_{3}, \tau_{4})
d\omega_{1}d\omega_{2}d\omega_{3}+\mathcal{O}(T^{-1}),\nonumber\\
\label{eqapp}\end{eqnarray}
\noindent where, for $\omega_{i}\in [-\pi,\pi],$ $i=1,2,3,$ \begin{eqnarray}&&\hspace*{-0.35cm} \mathcal{F}_{\omega_{1}, \omega_{2}, \omega_{3}}\underset{\mathcal{S}\left(L^{2}\left(\mathbb{M}_{d}^{2},\otimes_{i=1}^{2}\nu(dx_{i}), \mathbb{C}\right)\right)}{=}\frac{1}{(2\pi)^{3}}\sum_{t_{1},t_{2},t_{3}= -\infty}^{\infty}
\exp\left(\sum_{j=1}^{3}\omega_{j}t_{j}\right)\nonumber\\
&&\hspace*{7cm}\times
\mbox{cum}(X_{t_{1}},X_{t_{2}},X_{t_{3}},X_{0})\nonumber\\
\label{csdoo4}
\end{eqnarray}

\noindent
 denotes
the cumulant spectral density operator of order $4$ of $X,$  and, as before, $\underset{\mathcal{S}\left(L^{2}\left(\mathbb{M}_{d}^{2},\otimes_{i=1}^{2}\nu(dx_{i}), \mathbb{C}\right)\right)}{=}$
means the identity in the norm of the  space \linebreak $\mathcal{S}\left(L^{2}\left(\mathbb{M}_{d}^{2},\otimes_{i=1}^{2}\nu(dx_{i}), \mathbb{C}\right)\right)\equiv L^{2}(\mathbb{M}^{4}_{d}, \otimes_{i=1}^{4} d\nu (x_{i}),\mathbb{C}).$

\end{lemma}

\medskip

\noindent \textbf{Proof}. See  Appendix  \ref{a3b}.

\medskip

\section{A test for LRD in   functional time series on $\mathbb{M}_{d}$}

\label{tlrdfts}
\label{testlrd}  Consistency of  the test based on    $\mathcal{S}_{B_{T}}$ is derived in this section. Specifically, Theorem \ref{th2}  provides the almost surely divergence in the norm of the space $\mathcal{S}(L^{2}(\mathbb{M}_{d},d\nu, \mathbb{C}))$ of $\mathcal{S}_{B_{T}}$ under  $H_{1}.$ The proof of this result follows from Proposition \ref{cor1},  showing the divergence in the norm of the space $\mathcal{S}(L^{2}(\mathbb{M}_{d},d\nu, \mathbb{C}))$ of the centering operator of $\mathcal{S}_{B_{T}},$ and
 from Theorem \ref{consisimse}  and Corollary \ref{cor22}, establishing  the consistency of the integrated weighted periodogram operator under  $H_{1}.$ The implementation of the  testing procedure in practice is also discussed.

\begin{proposition}
\label{cor1}
Under $H_{1},$ as $T\to \infty,$
\begin{eqnarray}&&\hspace*{-0.7cm}\left\|E_{H_{1}}\left[\frac{\mathcal{S}_{B_{T}}}{\sqrt{TB_{T}}}\right]\right\|_{\mathcal{S}(L^{2}(\mathbb{M}_{d},d\nu, \mathbb{C}))}
=
\left\|\int_{[-\sqrt{B_{T}}/2, \sqrt{B_{T}}/2]}E_{H_{1}}[\widehat{\mathcal{F}}_{\omega }^{(T)}]\frac{d\omega}{\sqrt{B_{T}}} \right\|_{\mathcal{S}(L^{2}(\mathbb{M}_{d},d\nu, \mathbb{C}))}\nonumber\\
&&\hspace*{4cm}
\geq g(T)=\mathcal{O}(B_{T}^{-l_{\alpha }-1/2}).\label{eqorexh1hh}\end{eqnarray}
\end{proposition}

\medskip

\noindent \textbf{Proof}. See  Appendix  \ref{apcb}.

\medskip

\begin{theorem}
\label{consisimse} Under $H_{1},$ \emph{Assumption I},  and \begin{equation}\int_{[-\pi,\pi]}\|\mathcal{M}_{\omega }\|^{2}_{L^{1}(L^{2}(\mathbb{M}_{d},d\nu, \mathbb{C}))}|\omega |^{-2L_{\alpha }}d\omega <\infty,\label{ca}\end{equation}
\noindent   as $T\to \infty,$
\begin{eqnarray}
&&\hspace*{-0.7cm}\int_{-\pi}^{\pi}E_{H_{1}}\left\|\widehat{\mathcal{F}}_{\omega }^{(T)}-E_{H_{1}}[\widehat{\mathcal{F}}_{\omega }^{(T)}]\right\|_{\mathcal{S}(L^{2}(\mathbb{M}_{d},d\nu, \mathbb{C}))}^{2}d\omega \leq h(T)=\mathcal{O}(B_{T}^{-1}T^{-1}),
\label{ca2}
\end{eqnarray}
\noindent where, as before,  $\|\cdot\|_{L^{1}(L^{2}(\mathbb{M}_{d},d\nu, \mathbb{C}))}$ denotes the norm in the space
$L^{1}(L^{2}(\mathbb{M}_{d},d\nu, \mathbb{C}))$ of nuclear operators on $L^{2}(\mathbb{M}_{d},d\nu, \mathbb{C}).$
\end{theorem}

\begin{remark}
Note that condition (\ref{ca}) is satisfied, for instance, when the family $\left\{\mathcal{M}_{\omega },\ \omega \in [-\pi,\pi]\right\}$ lies in a ball of
radius $R>0$ of the space $L^{1}(L^{2}(\mathbb{M}_{d},d\nu, \mathbb{C})).$
\end{remark}

\medskip

\noindent \textbf{Proof}. See  Appendix \ref{appdd}.

\medskip

 Theorem \ref{consisimse}  implies the  weak consistency of the  integrated weighted periodogram operator   under $H_{1}$  in the norm of    the space $\mathcal{S}(L^{2}(\mathbb{M}_{d},d\nu, \mathbb{C})).$
\begin{corollary}
\label{cor22}
Under the conditions of Theorem \ref{consisimse}, as $T\to \infty,$
 \begin{eqnarray}
&&\left\|\int_{-\pi}^{\pi}E_{H_{1}}\left[\widehat{\mathcal{F}}_{\omega }^{(T)}-\int_{-\pi}^{\pi}W(\xi)\mathcal{F}_{\omega -B_{T}\xi }d\xi \right] d\omega
\right\|_{\mathcal{S}(L^{2}(\mathbb{M}_{d},d\nu ,\mathbb{C}))}\nonumber\\ &&\hspace*{6cm} \leq \widetilde{g}(T)=\mathcal{O}(T^{-1/2}B_{T}^{-1/2}).
\nonumber\end{eqnarray}
\end{corollary}

\medskip

\noindent \textbf{Proof}. See  Appendix \ref{apseccon}.

\medskip

Under the conditions assumed in the following result,  consistency of the test follows.
\begin{theorem}
\label{th2}
Under $H_{1},$ assume that   $l_{\alpha }> 1/4,$ and  that the  bandwidth parameter $B_{T}=T^{-\beta }$ for $\beta \in (0,1).$
If conditions of Theorem \ref{consisimse}  hold,  then,  as  $T\to \infty,$
$$
\left\|\mathcal{S}_{B_{T}}\right\|_{\mathcal{S}(L^{2}(\mathbb{M}_{d},d\nu, \mathbb{C}))}\to_{\mbox{a.s}} \infty,$$
\noindent where $\to_{\mbox{a.s.}} \infty$ denotes almost surely divergence.
\end{theorem}

\medskip

\noindent \textbf{Proof}. See  Appendix \ref{apd}.

\subsection{Practical implementation}
The practical implementation of the proposed statistical  testing procedure, in terms of the Gaussian  random projection methodology    (see Theorem 4.1 in \cite{Cuestalbertos}),  is now briefly discussed. The  Karhunen--Lo\'eve expansion in Lemma \ref{CLTSBT} below can be applied in such an implementation. Let us consider the random Fourier coefficients

\begin{eqnarray}&&
Y_{n,j,h,l}(\omega )=\frac{(\sqrt{2\pi}\|W\|_{L^{2}(\mathbb{R})})^{-1}}{\sqrt{f_{n}(\omega  )f_{h}(\omega  )}}\int_{\mathbb{M}_{d}^{2}}\widehat{f}_{\omega }(\tau ,\sigma )\overline{S_{n,j}^{d}(\tau)}S_{h,l}^{d}(\sigma)d\nu (\sigma )d\nu (\tau ),
\nonumber\\
&&\quad j=1,\dots, \Gamma (n,d),\quad  l=1,\dots, \Gamma(h,d), \ n,h\in \mathbb{N}_{0},\ \omega \in [-\pi,\pi ]\backslash \{0\},\nonumber\\
\label{KLSH2}\end{eqnarray}
\noindent   where integration is understood in the mean--square sense, and $\widehat{f}_{\omega }$ is the limit Gaussian random element in $\mathcal{S}(L^{2}(\mathbb{M}_{d},d\nu, \mathbb{C}))$ introduced in Lemma \ref{lem2pt}.
\begin{lemma}
\label{CLTSBT}
Let $\widehat{f}_{\omega }$ be defined, as before,  satisfying equation   (\ref{cobk}) in   Lemma \ref{lem2pt}.
Then, the  following series expansion holds in the mean--square sense:  For every $(\tau ,\sigma )\in \mathbb{M}_{d}^{2},$
\begin{eqnarray}
&&\frac{1}{\sqrt{2\pi}\|W\|_{L^{2}(\mathbb{R})}}\widehat{f}_{\omega }(\tau,\sigma)\underset{\mathcal{L}^{2}_{\mathcal{S}(L^{2}(\mathbb{M}_{d},d\nu; \mathbb{C}))}(\Omega ,\mathcal{A},\mathcal{P})}{=}\sum_{n,h\in \mathbb{N}_{0}}\sum_{j=1}^{\Gamma (n,d)}
\sum_{l=1}^{\Gamma (h,d)}
\sqrt{f_{n}(\omega)f_{h}(\omega)}\nonumber\\
&&\hspace*{3cm}\times
Y_{n,j,h,l}(\omega)S_{n,j}^{d}(\tau)\overline{S_{h,l}^{d}(\sigma)},\ \omega \in [-\pi,\pi]\backslash \{0\},
\label{KLSH}\end{eqnarray}
\noindent where, as before,   $\mathcal{L}^{2}_{\mathcal{S}(L^{2}(\mathbb{M}_{d},d\nu; \mathbb{C}))}(\Omega ,\mathcal{A},\mathcal{P})$ denotes the space of zero--mean second--order $\mathcal{S}(L^{2}(\mathbb{M}_{d},d\nu; \mathbb{C}))$--valued random variables with the norm \linebreak $\sqrt{E\|\cdot\|^{2}_{\mathcal{S}(L^{2}(\mathbb{M}_{d},d\nu; \mathbb{C}))}}.$
   The random Fourier coefficients  $$\left\{Y_{n,j,h,l}(\omega),\ j=1,\dots, \Gamma (n,d),\  l=1,\dots, \Gamma(h,d), \ n,h\in \mathbb{N}_{0}\right\},$$ \noindent for  $\omega \in [-\pi,\pi ]\backslash \{0\},$  have been introduced in equation (\ref{KLSH2}). They are
  independent and identically distributed complex--valued  standard  Gaussian    random variables.
\end{lemma}

\noindent \textbf{Proof}. See  Appendix \ref{apendicea0}.

\vspace*{0.5cm}

Theorem \ref{th2} motivates the methodology to be adopted in practice. Specifically, as illustrated in the simulation study undertaken in the next section,  a   consistent test for LRD is obtained by rejecting $H_{0},$ when, for every $j=1,\dots, \Gamma(n,d),$ $l=1,\dots,\Gamma (h,d),$ $n,h\in \mathbb{N}_{0},$  \begin{equation}\frac{\left|[\mathcal{S}_{B_{T}}-E[\mathcal{S}_{B_{T}}]](\overline{S_{h,l}^{d}})(S_{n,j}^{d})\right|}{\sqrt{\mbox{Var}
\left(\mathcal{S}_{B_{T}}(\overline{S_{h,l}^{d}})(S_{n,j}^{d})\right)}}\label{ptso}
 \end{equation}
 \noindent is larger than an upper tail  standard normal critical value.
 Note that  for $T$ sufficiently large, and for $n,h\in \mathbb{N}_{0},$
 \begin{eqnarray}
 &&\hspace*{-0.7cm}\left[\mathcal{S}_{B_{T}}-E[\mathcal{S}_{B_{T}}]\right](\overline{S_{h,l}^{d}})(S_{n,j}^{d})
 =\int_{\mathbb{M}_{d}^{2}}[\mathcal{S}_{B_{T}}-E[\mathcal{S}_{B_{T}}]](\tau,\sigma)\overline{S_{n,j}^{d}}(\tau)S_{h,l}^{d}(\sigma)d\nu(\sigma )d\nu(\tau )
\nonumber\\
&&\mbox{Var}\left(\mathcal{S}_{B_{T}}(\overline{S_{h,l}^{d}})(S_{n,j}^{d})\right)=\mbox{Var}\left(\int_{\mathbb{M}_{d}^{2}}
\mathcal{S}_{B_{T}}(\tau,\sigma)\overline{S_{n,j}^{d}}(\tau )S_{h,l}^{d}(\sigma )d\nu(\sigma )d\nu(\tau )
\right)\nonumber\\
&&=2\pi\int_{\left[-\sqrt{B_{T}}/2,\sqrt{B_{T}}/2\right]^{2}\times [-\pi,\pi] }W\left(\frac{\omega -\alpha }{B_{T}}\right)\left[W\left(\frac{\xi-\alpha }{B_{T}}\right)+W\left(\frac{\xi+\alpha }{B_{T}}\right)\right.\nonumber\\
&&\left.\times \left\langle S_{n,j}^{d}, \overline{S_{h,l}^{d}}\right\rangle_{L^{2}(\mathbb{M}_{d}, d\nu,\mathbb{C})}\left\langle \overline{S_{n,j}^{d}}, S_{h,l}^{d} \right\rangle_{L^{2}(\mathbb{M}_{d}, d\nu,\mathbb{C})}\right] f_{n}(\alpha )f_{h}(\alpha )\frac{d\alpha d\omega d\xi }{B_{T}^{2}}\nonumber\\
&&+\mathcal{O}(B_{T}^{-2}T^{-2})+\mathcal{O}(T^{-1}),\  j=1,\dots,\Gamma (n,d),\quad  l=1,\dots,\Gamma (h,d). \nonumber\\
\label{projts}
\end{eqnarray}

The associated dimensionality problem   can be substantially alleviated if we restrict our attention to the case where all moments of $\mathcal{S}_{B_{T}}$ are finite and satisfy the  Carleman condition.  In that case, Theorem 4.1 in \cite{Cuestalbertos} leads to the following   test statistic, evaluated conditionally to the observed  functional value $\mathbf{k}$ of a non--degenerated functional Gaussian random variable in the space $\mathcal{S}(L^{2}(\mathbb{M}_{d},d\nu; \mathbb{C})),$ whose  probability measure  on
$\mathcal{S}(L^{2}(\mathbb{M}_{d},d\nu; \mathbb{C}))$ is denoted as $\mu.$ Specifically,  consider
\begin{equation} \mathcal{T}_{B_{T}}^{\mathbf{k}}=\frac{\left|\left\langle \mathcal{S}_{B_{T}}-E[\mathcal{S}_{B_{T}}],\mathbf{k}\right\rangle_{\mathcal{S}(L^{2}(\mathbb{M}_{d},d\nu; \mathbb{C}))}\right|}{\sqrt{\mbox{Var}\left(\left\langle \mathcal{S}_{B_{T}}-E[\mathcal{S}_{B_{T}}],\mathbf{k}\right\rangle_{\mathcal{S}(L^{2}(\mathbb{M}_{d},d\nu; \mathbb{C}))}\right)}}.\label{ptsopa}
 \end{equation}
\noindent Then, $H_{0}^{\mathbf{k}}$ will be  rejected if  the observed value of   $\mathcal{T}_{B_{T}}^{\mathbf{k}}$  is larger than an upper tail standard normal critical value.
Note that, if $H_{0}$ holds then $H_{0}^{\mathbf{k}}$ also holds, and if $H_{0}$ fails then $H_{0}^{\mathbf{k}}$ also fails $\mu $--a.s. Thus, with probability one, we will generate  a realization of   random direction $\mathbf{k}$ in  $\mathcal{S}(L^{2}(\mathbb{M}_{d},d\nu; \mathbb{C}))$ for which $H_{0}^{\mathbf{k}}$ fails (see also \cite{Cuesta}).

In the spirit of the Gaussian random degree--$l$ spherical harmonics introduced   in \cite{Marrossiwig} ($l\geq 0$), we will  consider a truncated
 Karhunen--Lo\'eve expansion in the generation of a non--degenerated Gaussian measure characterizing the random direction $\mathbf{k},$
 where our test statistics is projected (see equation (\ref{ptsopa})). Specifically, we consider a zero--mean  Gaussian random element $\mathbf{k}$  in   $\mathcal{S}(L^{2}(\mathbb{M}_{d},d\nu; \mathbb{C}))$  with trace covariance operator $\mathcal{C}_{\mathbf{k}}$ having  kernel
 \begin{eqnarray} &&C_{\mathbf{k}}(\tau_{1},\sigma_{1},\tau_{2},\sigma_{2})=E\left[\mathbf{k}(\tau_{1},\sigma_{1})
\mathbf{k}(\tau_{2},\sigma_{2})\right]
\nonumber\\
&&=\sum_{n\in \mathbb{N}_{0}}\sum_{j=1}^{\Gamma(n,d)}\sum_{h\in \mathbb{N}_{0}}\sum_{l=1}^{\Gamma(h,d)}\lambda_{n, h}
S_{n,j}^{d}(\tau_{1})\overline{S_{h,l}^{d}(\sigma_{1})}\overline{S_{n,j}^{d}(\tau_{2})}S_{h,l}^{d}(\sigma_{2}),
 \end{eqnarray}
\noindent  for every  $(\tau_{1},\sigma_{1}), (\tau_{2},\sigma_{2})\in \mathbb{M}_{d}^{2}.$ The trace property of  $\mathcal{C}_{\mathbf{k}}$ can be equivalently expressed as $\sum_{n\in \mathbb{N}_{0}}\sum_{h\in \mathbb{N}_{0}}\Gamma (n,d)\Gamma(h,d)\lambda_{n,h}<\infty.$ Therefore, $\mathbf{k}$ admits a series expansion,   Karhunen--Lo\'eve expansion
(see, e.g., Lemma \ref{CLTSBT}), whose truncated version is implemented.  Particularly, one can construct centered isotropic Gaussian random fields
by finite--dimensional projection
   \begin{eqnarray}
   f_{n,h}(\tau ,\sigma )=\frac{1}{\Gamma(n,d)\Gamma(h,d)}\sum_{j=1}^{\Gamma(n,d)}\sum_{l=1}^{\Gamma(h,d)}Y_{n,j,h,l}S_{n,j}^{d}\otimes \overline{S_{h,l}^{d}}(\tau ,\sigma ),\nonumber\\
   \end{eqnarray}
\noindent   for
  $(\tau ,\sigma )\in  \mathbb{M}_{d}^{2},$ and $n,h \in \mathbb{N}_{0},$ involving the Gaussian random directions in  $\mathcal{S}(L^{2}(\mathbb{M}_{d},d\nu; \mathbb{C}))$
  \begin{equation}\mathbf{k}_{n,j,h,l}(\tau ,\sigma )=Y_{n,j,h,l}[S_{n,j}^{d}\otimes \overline{S_{h,l}^{d}}(\tau ,\sigma )],\quad (\tau ,\sigma )\in \mathbb{M}_{d}^{2},\label{eqgensph}
\end{equation}
\noindent  where $Y_{n,j,h,l}$ denotes a zero--mean Gaussian random variable with variance $\lambda_{n,h},$ for  $j=1,\dots \Gamma(n,d),$ $l=1,\dots,\Gamma (h,d),$ $n,h\in \mathbb{N}_{0}.$  Kernel $f_{n,h}(\tau ,\sigma )$ plays a key role in our  approach when
   $\mathcal{H}_{n}$ and/or $\mathcal{H}_{h}$ are  dominant eigenspaces of the Laplace Beltrami operator on $\mathbb{M}_{d}.$

\section{Simulation study}
\label{ni}
 Our simulations will be set on   $\mathbb{M}_{d}=\mathbb{S}_{d}\subset \mathbb{R}^{d+1}.$ An alternative   generation algorithm to  the ones considered  in  \cite{Ovalle23} and \cite{Ovalle24} is implemented reducing computational burden allowing for the consideration of large functional sample sizes.  Theorem \ref{lemvfh0}   is illustrated in the context of  SPHARMA(p,q) processes, and  the illustration of  Theorem \ref{th2}   is carried out in the context of    multifractionally integrated SPHARMA(p,q) processes.
  These numerical results are respectively  reported in Sections \ref{sgd}  and \ref{consec} for $\beta =1/4,$ i.e., $B_{T}=T^{-1/4}$   (see also  Section \ref{s6} for a wider analysis for  different $\beta $ values).
When  $L_{\alpha }>1/2,$ Section \ref{Sec5.3} opens new research lines  beyond condition (\ref{shconvrm}), providing empirical evidence of a faster  a.s.  divergence rate of $\mathcal{S}_{B_{T}}$ in the norm of the space $\mathcal{S}(L^{2}(\mathbb{M}_{d},d\nu ,\mathbb{C})).$
 Finally, Section \ref{sesptp} shows empirical  size and power properties of the testing procedure.

\subsection{Asymptotic Gaussian distribution of $\mathcal{S}_{B_{T}}$ under $H_{0}$}
\label{sgd}
Let us consider that the elements of  the family of spectral density operators of $X$ have  frequency varying eigenvalues, with respect to the system of eigenfunctions of the Laplace--Beltrami operator,   obeying the following equation under $H_{0}$ (see \cite{RuizMedina2022}):
\begin{eqnarray}
f_{n}(\omega ) &=&\frac{\lambda_{n}(\mathcal{R}^{\eta}_{0})}{2\pi}\left|\frac{\Psi_{q,n}(\exp(-i\omega ))}{\Phi_{p,n}(\exp(-i\omega ))}\right|^{2},\ n\in \mathbb{N}_{0},\ \omega \in [-\pi,\pi],\nonumber\\
\label{eqsim1}
\end{eqnarray}
\noindent where $\left\{\lambda_{n}(\mathcal{R}^{\eta}_{0}),\ n\in \mathbb{N}_{0}\right\}$ is the system of eigenvalues of the autocovariance operator $\mathcal{R}^{\eta}_{0}$ of the innovation process $\eta=\{\eta_{t},\ t\in \mathbb{Z}\},$ with respect to the system of eigenfunctions of the Laplace--Beltrami operator.  Process   $\eta $ is assumed to be strong--white noise in $L^{2}(\mathbb{S}_{d},d\nu, \mathbb{R}).$ That is, $\eta $ is assumed to be a
sequence of independent and identically distributed $L^{2}(\mathbb{S}_{d},d\nu, \mathbb{R})$--valued random variables
such that $E[\eta_{t}]=0,$ and $E[\eta_{t}\otimes
\eta_{s}]=\delta_{t,s}\mathcal{R}^{\eta}_{0},$ with
$\mathcal{R}^{\eta}_{0}\in L^{1}(L^{2}(\mathbb{S}_{d},d\nu, \mathbb{R})),$ and $\delta_{t,s}=0,$
for $t\neq s,$ and $\delta_{t,s}=1,$ for $t=s.$
For $n\in \mathbb{N}_{0},$
  $\Phi_{p,n}(z)=1-\sum_{j=1}^{p}\lambda_{n}(\varphi_{j})z^{j}$ and $\Psi_{q,n}(z)=\sum_{j=1}^{q}\lambda_{n}(\psi_{j})z^{j},$
 with $\left\{\lambda_{n}(\varphi_{j}),\ n\in \mathbb{N}_{0}\right\}$ and  $\left\{\lambda_{n}(\psi_{l}),\ n\in \mathbb{N}_{0}\right\}$
     denoting the sequences of eigenvalues, with respect to the system of eigenfunctions of the Laplace--Beltrami operator,  of the  self--adjoint invariant integral operators $\varphi_{j}$   and $\psi_{l},$ for $j=1,\dots,p,$ and $l=1,\dots,q,$ respectively.   These operators  satisfy the following equations:
\begin{eqnarray}
\Phi_{p}(B)=1-\sum_{j=1}^{p}\varphi_{j}B^{j},\quad \Psi_{q}(B)=\sum_{j=1}^{q}\psi_{j}B^{j},\nonumber 
\label{rf}
\end{eqnarray}
\noindent where   $B$ is a difference operator such that
\begin{equation}
E\|B^{j}X_{t}-X_{t-j}\|_{L^{2}(\mathbb{S}_{d},d\nu, \mathbb{R})}^{2}=0,\quad \forall t,j\in \mathbb{Z}.
\label{ihnorndo}
\end{equation}
\noindent Here, $\Phi_{p}$ and  $\Psi_{q}$ are the so--called autoregressive and moving average operators, respectively.
 Also, for each $n\in \mathbb{N}_{0},$  $\Phi_{p,n}(z)=1-\sum_{j=1}^{p}\lambda_{n}(\varphi_{j})z^{j}$ and $\Psi_{q,n}(z)=\sum_{j=1}^{q}\lambda_{n}(\psi_{j})z^{j}$ have not common
roots, and  their roots are outside of the unit circle (see also Corollary 6.17 in \cite{Beran17}). Thus, $X$ satisfies   an SPHARMA(p,q) equation (see also \cite{Caponera21}; \cite{CaponeraMarinucci}).

\begin{figure}[h]
\begin{center}
\includegraphics[width=10cm,height=11cm]{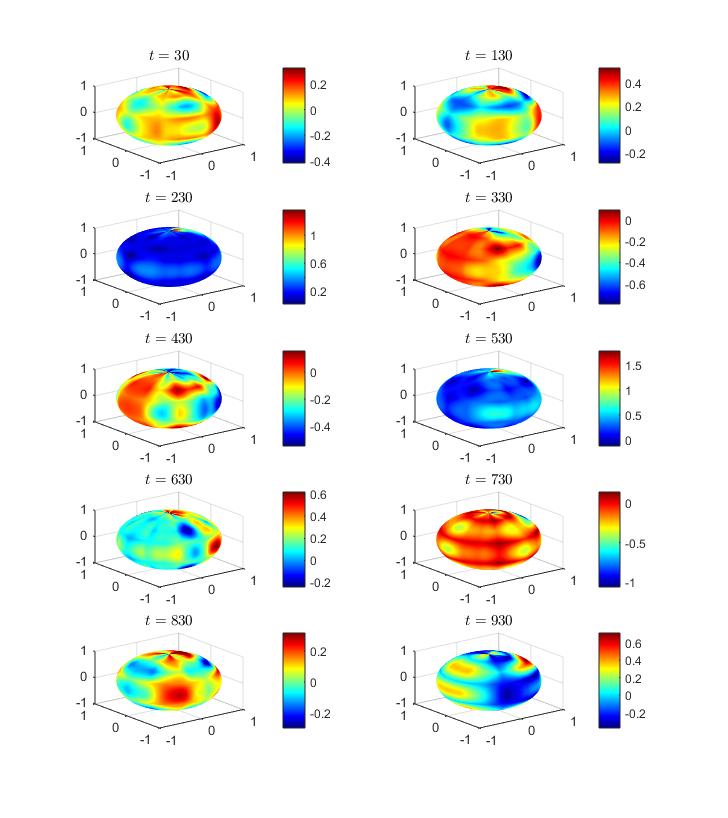}
\end{center}
\caption{\scriptsize{One realization at times}  $t=30, 130, 230, 330, 430, 530, 630, 730, 830,930$  \scriptsize{of    SPHARMA(1,1) process }
$\left(\lambda_{n}(\varphi_{1})=0.7\left(\frac{n+1}{n}\right)^{-3/2}\right.,$ \scriptsize{and}
$\left.\lambda_{n}(\psi_{1})= (0.4)\left(\frac{n+1}{n}\right)^{-5/1.95},\ n=1,2,3,4,5,6,7,8\right),$ \scriptsize{projected  into the direct sum} $\bigoplus_{n=1}^{8}\mathcal{H}_{n}$ \scriptsize{of eigenspaces} $\mathcal{H}_{n},$ $n=1,\dots,8,$ of $\Delta_{2}$}
\label{f1CLT}
\end{figure}

In the simulations we have generated an SPHARMA(1,1) process, i.e.,   $p=q=1,$ with $H=L^{2}(\mathbb{S}_{2},d\nu, \mathbb{R}),$  and $\lambda_{n}(\varphi_{1})=0.7\left(\frac{n+1}{n}\right)^{-3/2}$ and
$\lambda_{n}(\psi_{1})= (0.4)\left(\frac{n+1}{n}\right)^{-5/1.95},$ $n\in \mathbb{N}_{0}.$  Figure \ref{f1CLT} displays  one realization of the   generated SPHARMA(1,1) process projected into $\bigoplus_{n=1}^{8}\mathcal{H}_{n},$
at times   $t=30, 130, 230, 330,$  \linebreak $430, 530, 630, 730, 830,930.$

\clearpage

\begin{figure}[!h]
\begin{center}
\includegraphics[width=3cm,height=2.3cm]{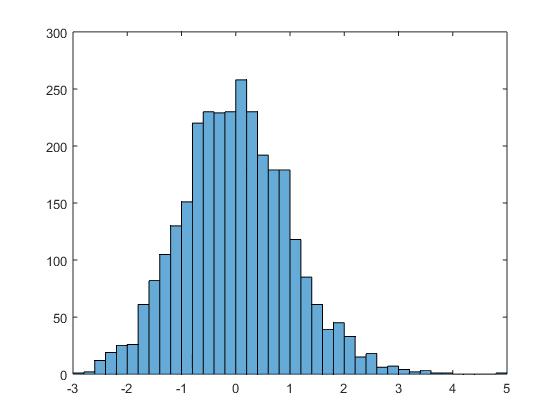}
\includegraphics[width=3cm,height=2.3cm]{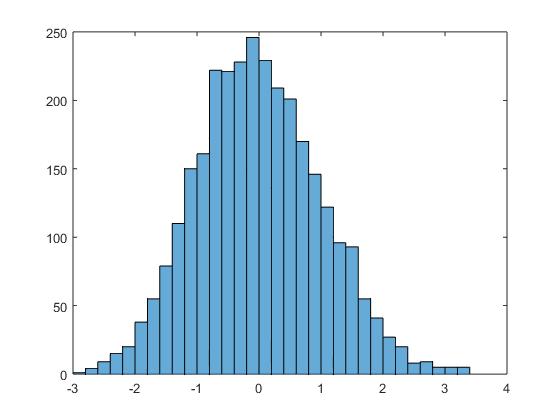}
\includegraphics[width=3cm,height=2.3cm]{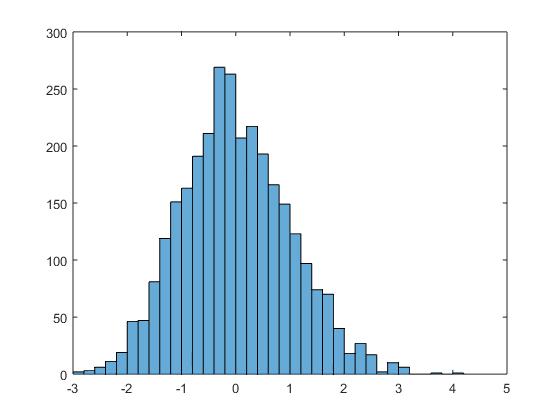}
\includegraphics[width=3cm,height=2.3cm]{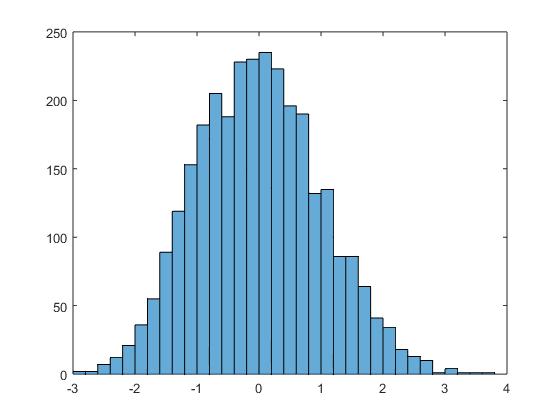}
\includegraphics[width=3cm,height=2.3cm]{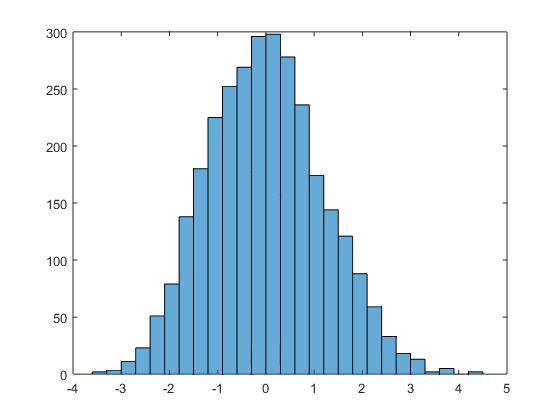}
\includegraphics[width=3cm,height=2.3cm]{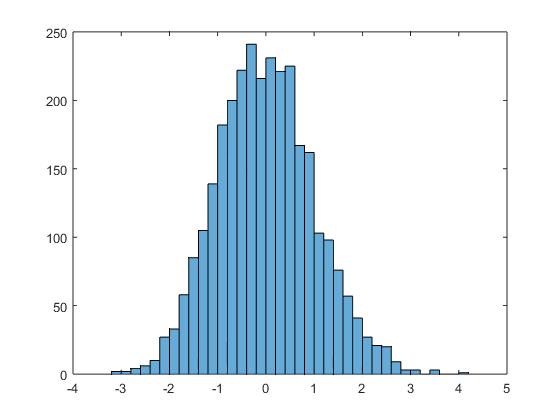}
\includegraphics[width=3cm,height=2.3cm]{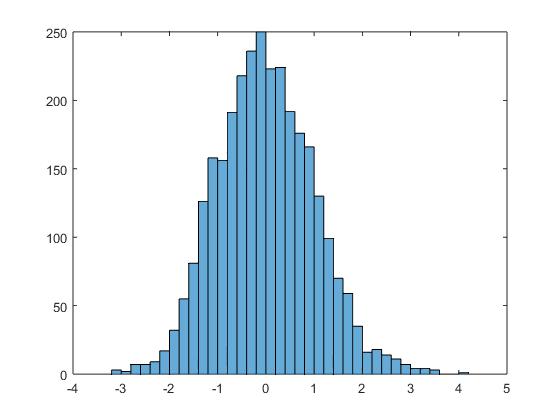}
\includegraphics[width=3cm,height=2.3cm]{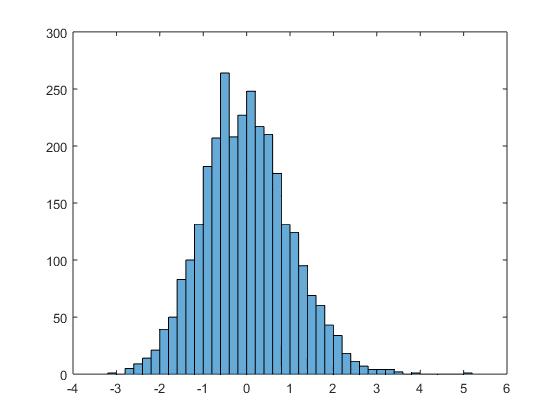}
 \end{center}
\caption{\scriptsize{Empirical projections of the probability measure of standardized $\mathcal{S}_{B_{T}},$ $B_{T}=T^{-1/4},$ under $H_{0},$   into the eigenspaces $\mathcal{H}_{n}\otimes \mathcal{H}_{n},$ for  $n=1,2,3,4,5,6,7,8,$ respectively displayed from the left to the right, and from the top to the bottom,   for functional samples size $T=1000$ and  $R=3000$ repetitions}}
\label{emd1}
\end{figure}

\begin{figure}[!h]
\begin{center}
\includegraphics[width=3cm,height=2.3cm]{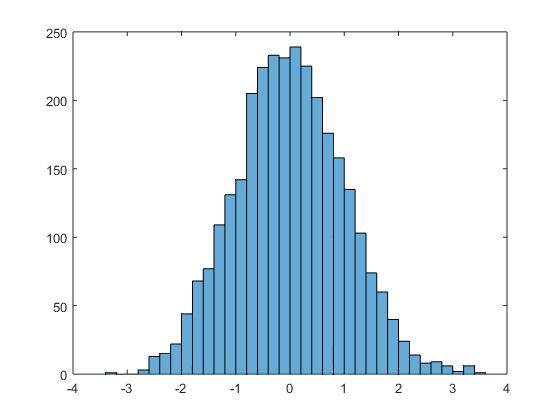}
\includegraphics[width=3cm,height=2.3cm]{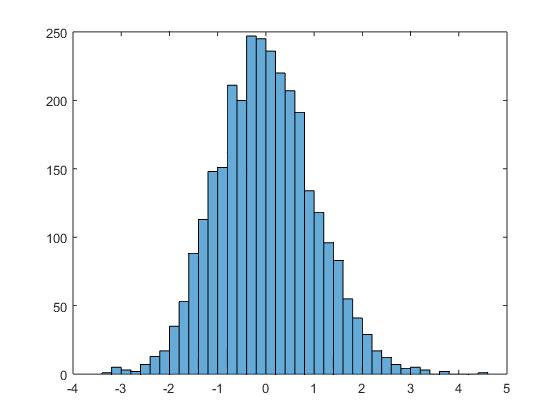}
\includegraphics[width=3cm,height=2.3cm]{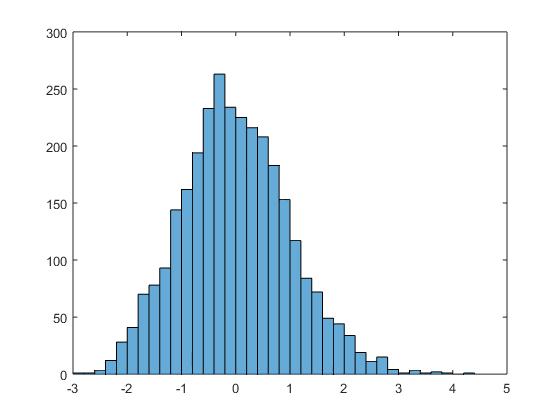}
\includegraphics[width=3cm,height=2.3cm]{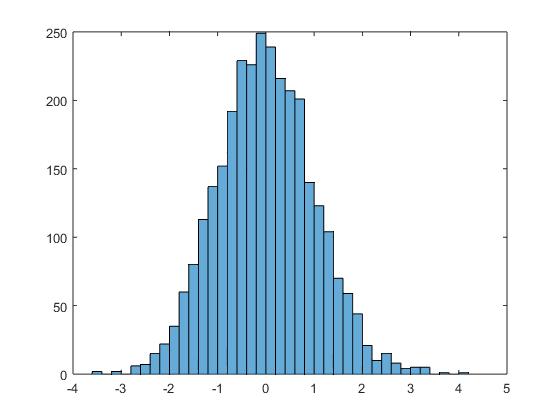}
\includegraphics[width=3cm,height=2.3cm]{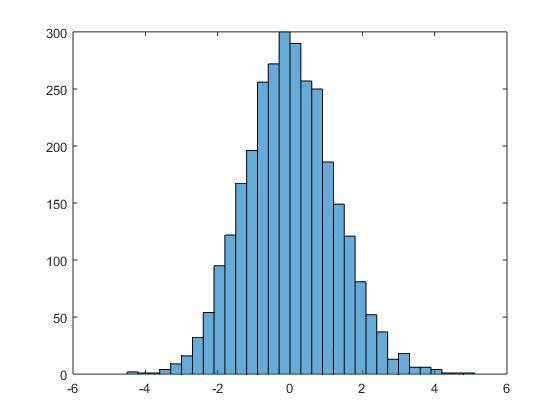}
\includegraphics[width=3cm,height=2.3cm]{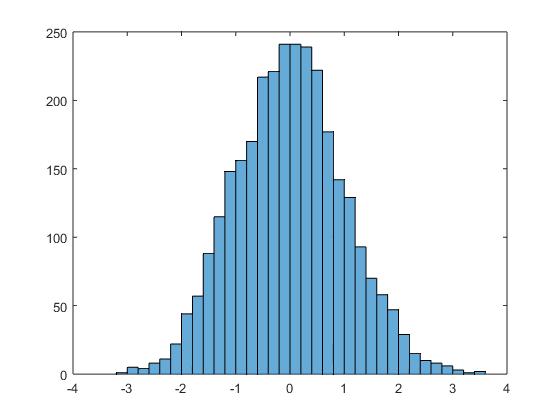}
\includegraphics[width=3cm,height=2.3cm]{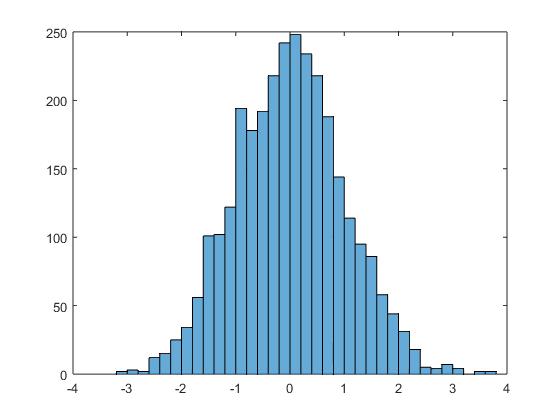}
\includegraphics[width=3cm,height=2.3cm]{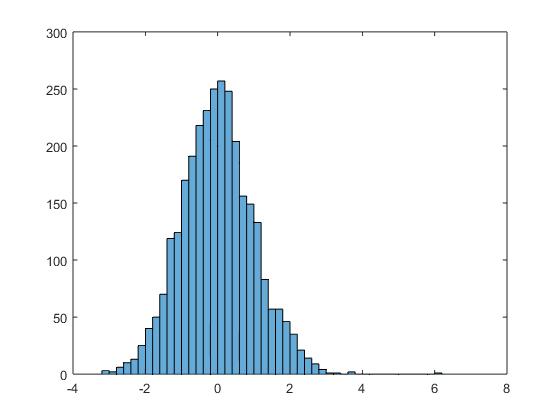}
 \end{center}
\caption{\scriptsize{Empirical projections of the probability measure of standardized $\mathcal{S}_{B_{T}},$ $B_{T}=T^{-1/4},$ under $H_{0},$   into the eigenspaces $\mathcal{H}_{n}\otimes \mathcal{H}_{n},$ for  $n=1,2,3,4,5,6,7,8,$ respectively displayed from the left to the right, and from the top to the bottom,   for functional samples size $T=2000$ and  $R=3000$ repetitions}}
\label{emd2}
\end{figure}

\begin{figure}[!h]
\begin{center}
\includegraphics[width=3cm,height=2.3cm]{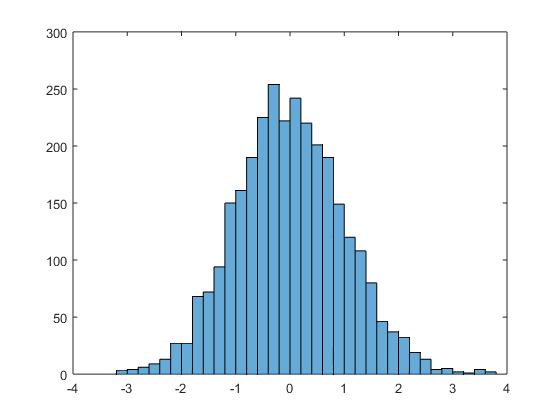}
\includegraphics[width=3cm,height=2.3cm]{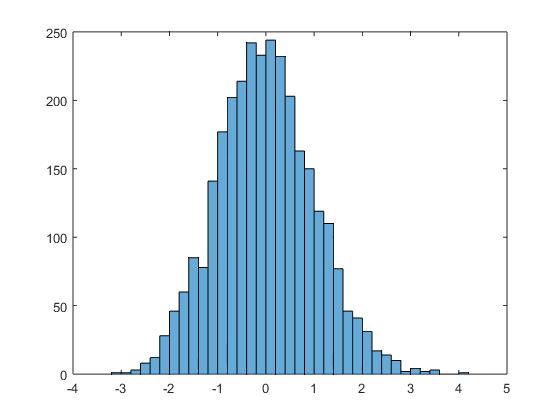}
\includegraphics[width=3cm,height=2.3cm]{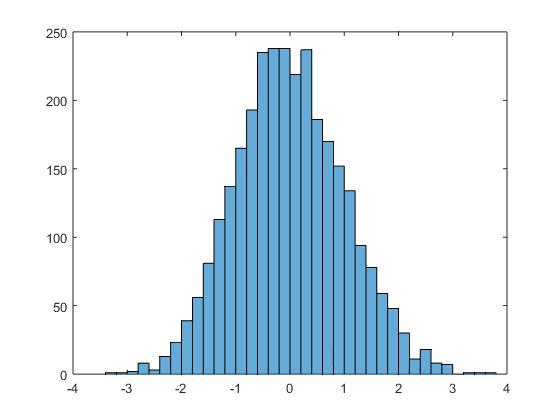}
\includegraphics[width=3cm,height=2.3cm]{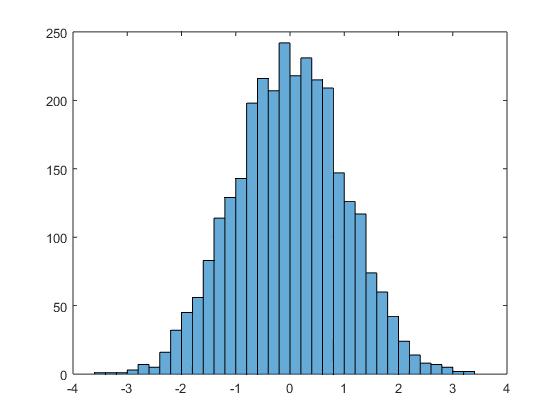}
\includegraphics[width=3cm,height=2.3cm]{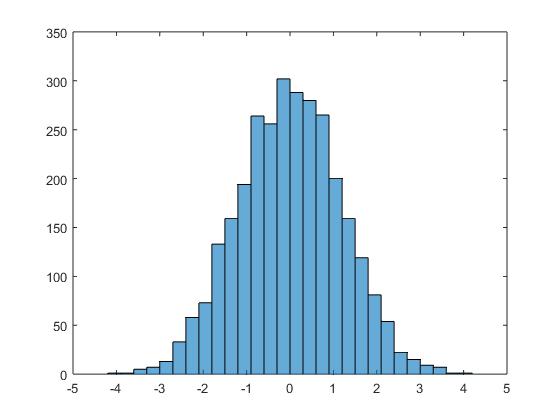}
\includegraphics[width=3cm,height=2.3cm]{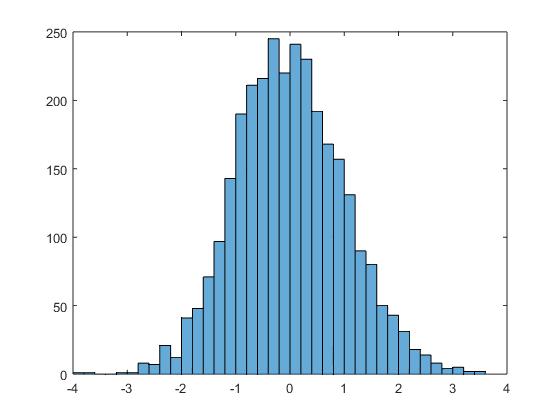}
\includegraphics[width=3cm,height=2.3cm]{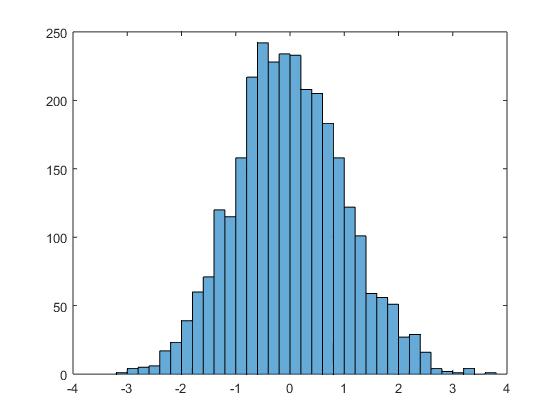}
\includegraphics[width=3cm,height=2.3cm]{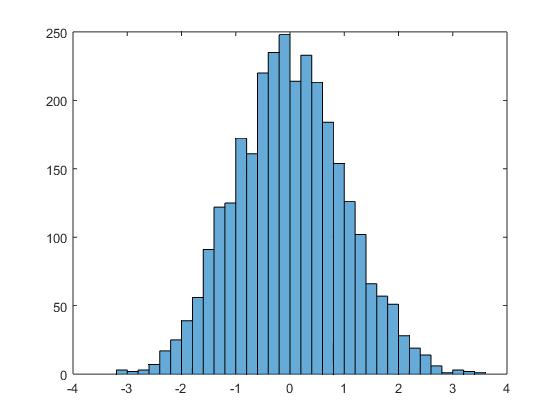}
 \end{center}
\caption{\scriptsize{Empirical projections of the probability measure of standardized $\mathcal{S}_{B_{T}},$ $B_{T}=T^{-1/4},$ under $H_{0},$   into the eigenspaces $\mathcal{H}_{n}\otimes \mathcal{H}_{n},$ for  $n=1,2,3,4,5,6,7,8,$ respectively displayed from the left to the right, and from the top to the bottom,   for functional samples size $T=3000$ and  $R=3000$ repetitions}}
\label{emd3}
\end{figure}

For each $n=1,\dots,8,$ the empirical distribution of  the centered and standardized projection into  $\mathcal{H}_{n}\otimes \mathcal{H}_{n}$ of $\mathcal{S}_{B_{T}}$
is displayed in Figure \ref{emd1} for functional sample size $T=1000$ and $R=3000$ repetitions,  in Figure \ref{emd2}  for functional sample size $T=2000$ and $R=3000$ repetitions,   and  in Figure \ref{emd3} for functional sample size $T=3000$ and $R=3000$ repetitions. These empirical distributions approximate the   support and shape of   a standard Gaussian probability density. The empirical  standardization displays a decreasing pattern over the spherical scale  $n,$ meaning that the respective  limit one--dimensional  Gaussian probability measures of these projections   have decreasing support. According to  Theorem 1.2.1 in \cite{Prato02}, this property is satisfied by the    infinite product Gaussian measure on   $(\mathbb{R}^{\infty},\mathcal{B}(\mathbb{R}^{\infty})),$  whose restriction to  $L^{2}(\mathbb{S}_{d},d\nu, \mathbb{C})$ is identified in the $\ell^{2}$--sense with the  probability measure of  the limit  Gaussian random element  in Theorem \ref{lemvfh0}.


\subsection{Consistency  of the test}
\label{consec}
Under $H_{1},$ for each $\omega \in [-\pi,\pi],$ the   eigenvalues $\left\{ f_{n}(\omega ),\ n\in \mathbb{N}_{0}\right\}$ satisfy (see equations (\ref{eqscm1}), (\ref{ffsvint}) and (\ref{reatlrd}))
\begin{eqnarray}
f_{n}(\omega ) &=&\frac{\lambda_{n}(\mathcal{R}^{\eta}_{0})}{2\pi}\left|\frac{\Psi_{q,n}(\exp(-i\omega ))}{\Phi_{p,n}(\exp(-i\omega ))}\right|^{2}\left|1-\exp\left(-i\omega \right)\right|^{-\alpha (n)},\ n\in \mathbb{N}_{0}.
\label{eqsim2bb}
\end{eqnarray}
 \noindent  Again,  we consider the projection of $X$  into $\bigoplus_{n=1}^{8}\mathcal{H}_{n}.$ Three multifractional integration operators, applied to SPHARMA(1,1) process  generated in Section \ref{sgd},
 are considered in  Sections \ref{ex1}--\ref{ex3}.  In   Example 1,   $\alpha (n)$ is decreasing over $n,$    in Example 2 $\alpha (n)$ is increasing over $n,$ and non--monotone  in  Example 3.  Note that,   under the generated Gaussian scenario,   condition (\ref{shconvrm}) implies that  Assumption I holds. Furthermore, condition  (\ref{ca}) also holds since $\varphi_{1}$ lies in  the unit ball of the space $\mathcal{L}(L^{2}(\mathbb{S}_{2}, d\nu,\mathbb{C})),$ and $\psi_{1}$ belongs to the trace class.

\subsubsection{Example 1}
\label{ex1}
\begin{figure}[!h]
\begin{center}
\includegraphics[width=10cm,height=11cm]{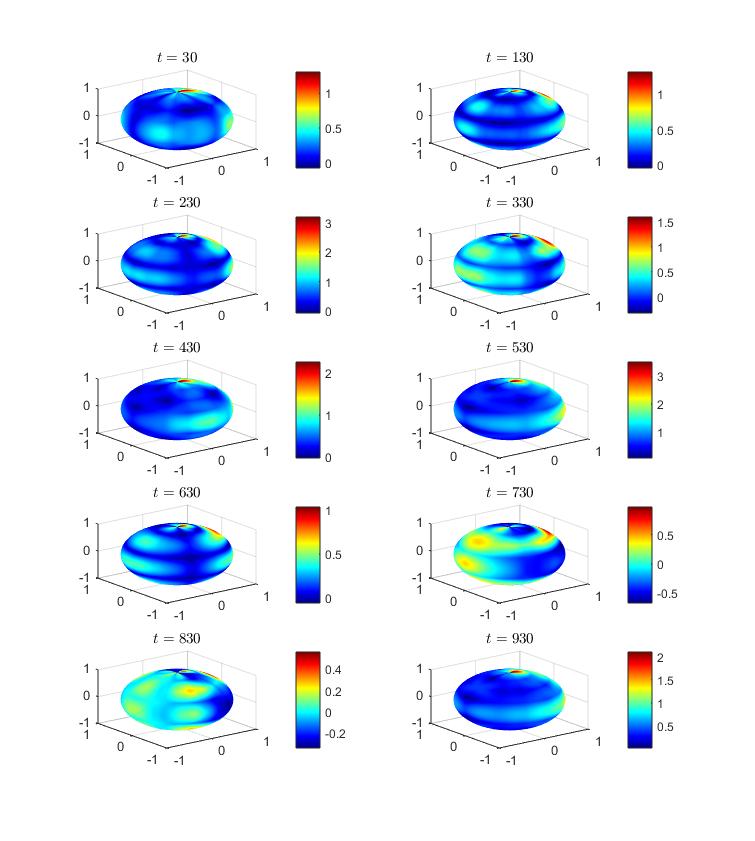}
\end{center}
\caption{\scriptsize{\emph{Example 1}. One sample realization at times $t=30, 130, 230, 330,  430, 530, 630, 730, 830,930$ of   multifractionally integrated  SPHARMA(1,1) process projected into $\bigoplus_{n=1}^{8}\mathcal{H}_{n}$ }}
\label{f6app2}
\end{figure}

 Theorem \ref{th2} is now illustrated in the case where  the largest dependence range is displayed by the process projected   into the eigenspace  $\mathcal{H}_{1}.$ Figure \ref{f6app2}  displays a   sample realization of the corresponding multifractionally integrated SPHARMA (1,1) process $X$  projected into $\bigoplus_{i=1}^{8}\mathcal{H}_{i}.$

\begin{figure}[!h]
\begin{center}
\includegraphics[width=2.5cm,height=2.5cm]{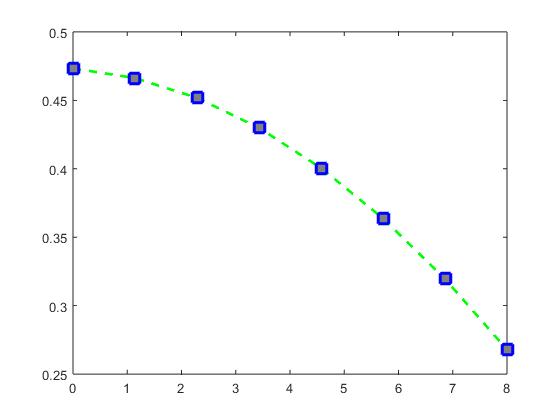}
\includegraphics[width=2.5cm,height=2.5cm]{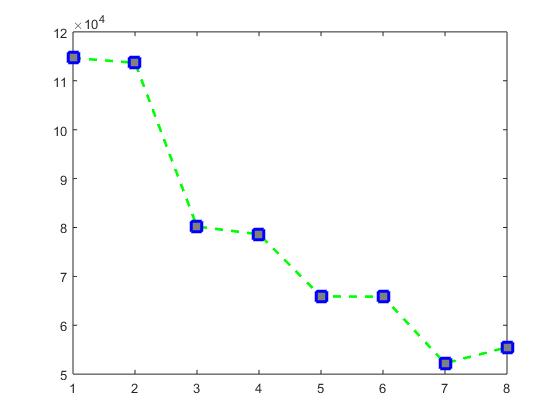}
\includegraphics[width=2.5cm,height=2.5cm]{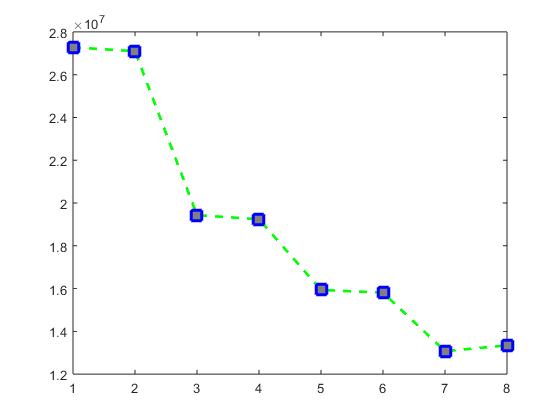}
\includegraphics[width=2.5cm,height=2.5cm]{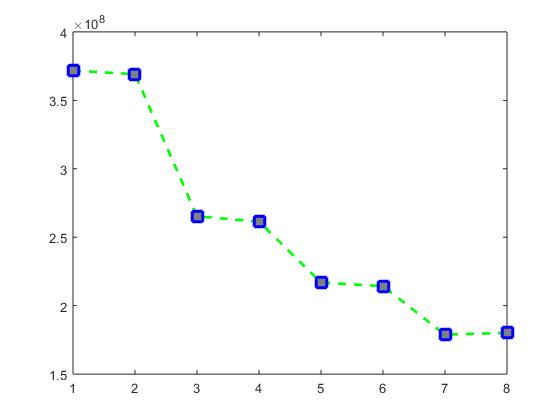}
\end{center}
\caption{\scriptsize{\emph{Example 1}.  Eigenvalues $\alpha (n),$ $n=1,2,3,4,5,6,7,8,$  of  LRD operator $\mathcal{A},$ $L_{\alpha }= 0.4733,$ and $l_{\alpha }= 0.2678$  (plot at the left--hand side).
Sample values of the test operator statistic $\mathcal{S}_{B_{T}},$ $B_{T}= T^{-1/4},$ projected
into  $\mathcal{H}_{n}\otimes \mathcal{H}_{n},$ $n=1,\dots,8$  for functional sample sizes $T=1000,10000,30000$ (three plots at the right--hand side)}}
\label{SBTDIVEX1}
\end{figure}

  In this example,
  $L_{\alpha }= 0.4733,$  $l_{\alpha }= 0.2678,$ and   $\alpha (n)=l_{\alpha }= 0.2678,$ $n\geq 9$
(see plot at the left--hand side  of Figure  \ref{SBTDIVEX1}).

The a.s. divergence of $\mathcal{S}_{B_{T}},$ for $B_{T}= T^{-1/4},$ in the  Hilbert--Schmidt operator norm  (see  Table \ref{T1}) is also reflected in the observed  increasing sample values of each one of its projections   into   $\mathcal{H}_{n}\otimes \mathcal{H}_{n},$ $n=1,\dots,8,$ for the  increasing  functional samples sizes $T=1000, 10000, 30000$   (see the three plots at the right--hand side of Figure \ref{SBTDIVEX1}).

\subsubsection{Example 2}
\label{ex2}

\begin{figure}[!h]
\begin{center}
\includegraphics[width=10cm,height=11cm]{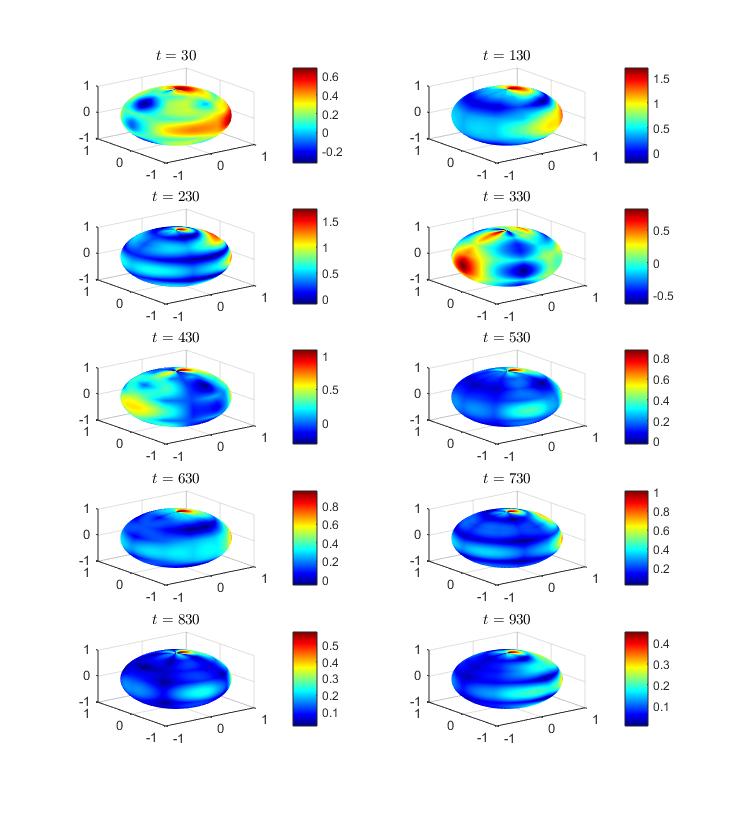}
\end{center}
\caption{\scriptsize{\emph{Example 2}. One sample realization at times $t=30, 130, 230, 330,  430, 530, 630, 730, 830,930$ of   multifractionally integrated  SPHARMA(1,1) process  projected into $\bigoplus_{n=1}^{8}\mathcal{H}_{n}$  }}
\label{f6app2ex2}
\end{figure}

The dominant subspace in this example, where the projected process displays the largest dependence range,  is  eigenspace $\mathcal{H}_{8}.$ One sample realization of the generated   multifractionally integrated  SPHARMA (1,1) process, projected into $\bigoplus_{n=1}^{8}\mathcal{H}_{n},$   is displayed in Figure \ref{f6app2ex2}.

\begin{figure}[!h]
\begin{center}
\includegraphics[width=2.5cm,height=2.5cm]{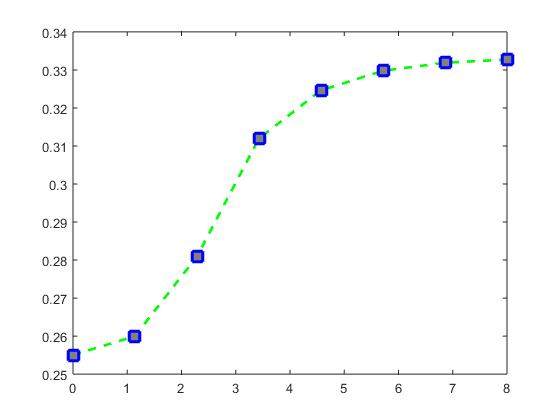}
\includegraphics[width=2.5cm,height=2.5cm]{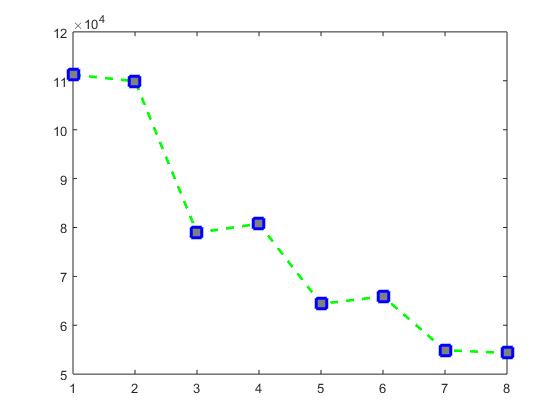}
\includegraphics[width=2.5cm,height=2.5cm]{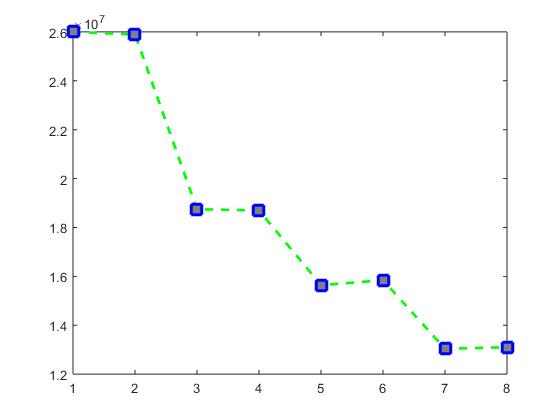}
\includegraphics[width=2.5cm,height=2.5cm]{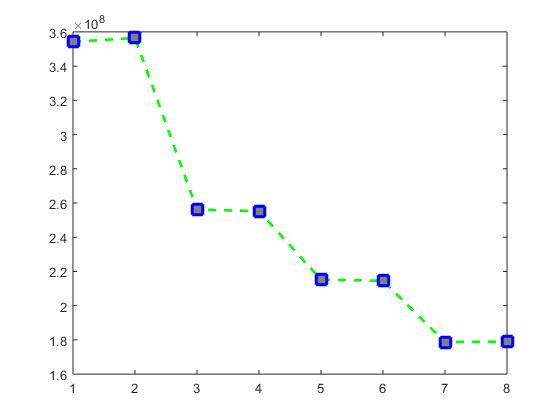}
\end{center}
\caption{\scriptsize{\emph{Example 2}.  Eigenvalues $\alpha (n),$ $n=1,2,3,4,5,6,7,8,$ of  LRD operator $\mathcal{A},$ $L_{\alpha }=   0.3327,$ and $l_{\alpha }=0.2550$ (plot at the left--hand side). Sample values of the kernel of the  test operator statistic $\mathcal{S}_{B_{T}},$ $B_{T}= T^{-1/4},$ projected
into  $\mathcal{H}_{n}\otimes \mathcal{H}_{n},$ $n=1,\dots,8,$    for functional sample sizes $T=1000,10000,30000$ (three plots at the right--hand side)}}
\label{f6app2ex2b}
\end{figure}

The LRD operator eigenvalues $\alpha (n),$ $n=1,2,3,4,5,6,7,8,$ are given in the plot at the left--hand side of  Figure \ref{f6app2ex2b},
 where  $L_{\alpha }=   0.3327,$  $l_{\alpha }=0.2550,$ and $\alpha (n)=l_{\alpha }= 0.2550,$ $n\geq 9.$
The a.s. divergence of our test statistic operator in the Hilbert--Schmidt operator norm (see  also Table \ref{T1}) is   illustrated  in the three plots at the right--hand side of such Figure \ref{f6app2ex2b}, in terms of the sample values of each projection of  $\mathcal{S}_{B_{T}},$ $B_{T}= T^{-1/4},$ into $\mathcal{H}_{n}\otimes \mathcal{H}_{n},$ $n=1,\dots, 8,$  for increasing
functional samples sizes $T=1000, 10000, 30000.$

\subsubsection{Example 3}
\label{ex3}
\begin{figure}[!h]
\begin{center}
\includegraphics[width=10cm,height=11cm]{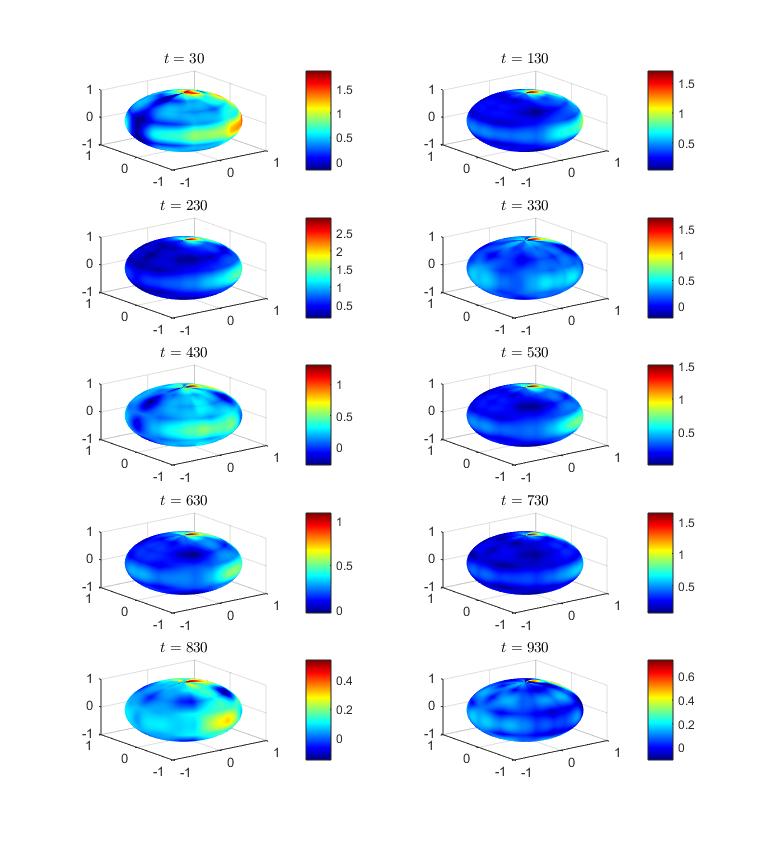}
\end{center}
\caption{\scriptsize{\emph{Example 3}. One sample realization at times $t=30, 130, 230, 330,  430, 530, 630, 730, 830,930$ of  multifractionally integrated  SPHARMA(1,1) process projected into $\bigoplus_{n=1}^{8}\mathcal{H}_{n}$ }}
\label{f6appex3BSR}
\end{figure}

In this third  example,    the dominant   subspace   is  eigenspace  $\mathcal{H}_{5}$ of the Laplace--Beltrami operator $\Delta_{2}.$
One sample realization of the generated multifractionally integrated SPHARMA (1,1) process,  projected into $\bigoplus_{n=1}^{8}\mathcal{H}_{n},$  is displayed in Figure \ref{f6appex3BSR}.

\begin{figure}[!h]
\begin{center}
\includegraphics[width=2.5cm,height=2.5cm]{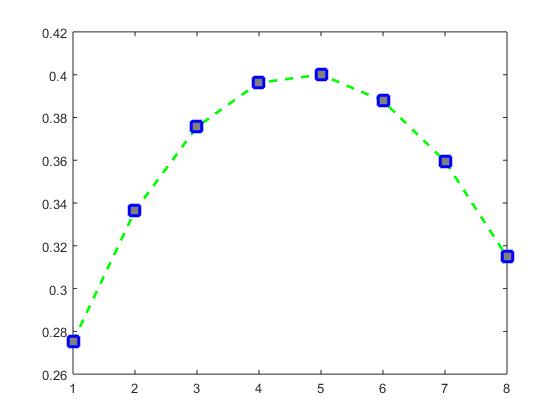}
\includegraphics[width=2.5cm,height=2.5cm]{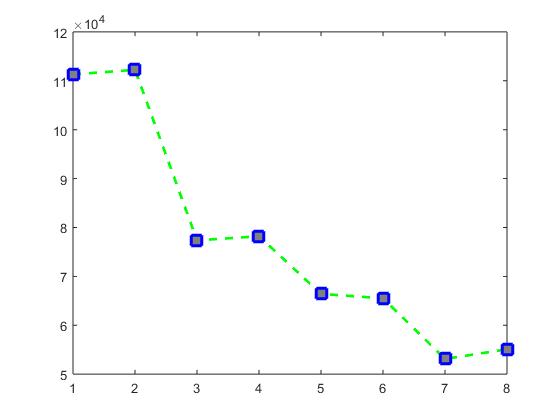}
\includegraphics[width=2.5cm,height=2.5cm]{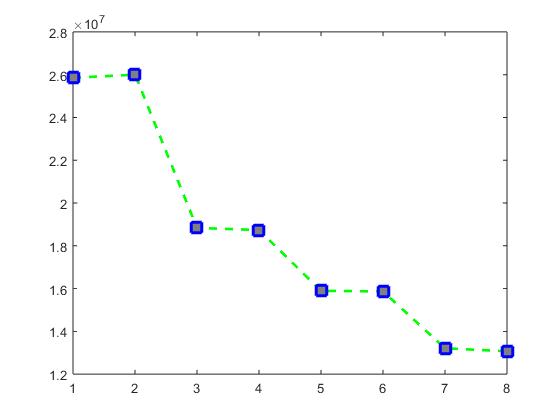}
\includegraphics[width=2.5cm,height=2.5cm]{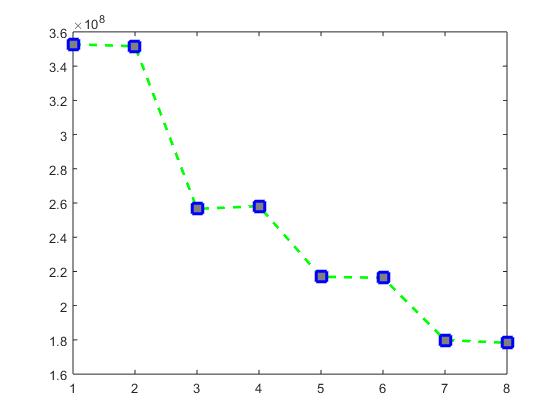}
\end{center}
\caption{\scriptsize{\emph{Example 3}.  Eigenvalues $\alpha (n),$ $n=1,2,3,4,5,6,7,8,$ of the LRD operator $\mathcal{A},$ $L_{\alpha }= 0.4000,$ and $l_{\alpha }= 0.2753$ (plot at the left--hand side). Sample values of the kernel of the test operator statistic $\mathcal{S}_{B_{T}},$ $B_{T}= T^{-1/4},$ projected
into  $\mathcal{H}_{n}\otimes \mathcal{H}_{n},$ $n=1,\dots,8,$  for functional sample sizes $T=1000,10000,30000$ (three plots at the right--hand side)}}
\label{f6app2ex3B}
\end{figure}

The eigenvalues $\alpha (n),$ $n=1,2,3,4,5,6,7,8,$ of
LRD operator $\mathcal{A}$ are showed in the plot at the left--hand side of  Figure \ref{f6app2ex3B} with
  $L_{\alpha }= 0.4000,$ and  $l_{\alpha }= 0.2753=\alpha (n),$ $n\geq 9.$   The sample values of the projections of  $\mathcal{S}_{B_{T}},$ $B_{T}= T^{-1/4},$ into $\mathcal{H}_{n}\otimes \mathcal{H}_{n},$  $n=1,\dots,8,$ for functional samples sizes $T=1000, 10000, 30000,$  can be found in the three plots at the right--hand side of   Figure \ref{f6app2ex3B}   (see also  Table \ref{T1}).

\subsubsection{Almost surely divergence  of $\mathcal{S}_{B_{T}}$  in $\mathcal{S}(L^{2}(\mathbb{M}_{d},d\nu, \mathbb{C}))$ norm under $H_{1}$}

 The  observed values of the Hilbert-Schmidt operator norm
of $\mathcal{S}_{B_{T}},$ for $B_{T}= T^{-1/4},$  projected into   $\bigoplus_{n=1}^{8}\mathcal{H}_{n}\otimes \mathcal{H}_{n},$  is displayed in
Table \ref{T1}, for the three numerical examples generated, and  for the  functional sample sizes   $T=1000, 5000,$ \linebreak  $10000,30000, 50000,100000.$
 \begin{table*}
\caption{\textbf{\emph{Hilbert-Schmidt operator norm of projected $\mathcal{S}_{B_{T}},$ $\beta =1/4$} }}
\label{T1}
\begin{tabular}{@{}lrrrrc@{}}
\hline
 {\bf T}& \multicolumn{1}{c}{{\bf Example 1}} &\multicolumn{1}{c}{{\bf Example 2}}    &\multicolumn{1}{c}{{\bf Example 3}} \\
  \hline
1000 &    2.3036e+05 &     2.2595e+05     &   1.9934e+05\\
5000 &   1.0612e+07  &       1.0223e+07&    9.0697e+06
\\
10000 &    5.5172e+07 &    5.3770e+07 &   4.7303e+07\\
30000&    7.5695e+08 &     7.3377e+08 & 6.4742e+08\\
50000& 2.5516e+09 &  2.4764e+09 & 2.1844e+09\\
100000 & 1.3256e+10 & 1.2892e+10 &    1.2906e+10\\
  \hline
\end{tabular}
\end{table*}
One can observe in  Table \ref{T1}  the increasing sample values of the  Hilbert--Schmidt operator norm of the projected $\mathcal{S}_{B_{T}}$ as $T$ increases, in all the examples under  $l_{\alpha} >1/4,$    with $B_{T}= T^{-\beta},$
 $\beta =1/4,$   satisfying  $TB_{T}\to \infty,$ $T\to \infty.$
 Spherical sample patterns and scales induced by the multifractional integration operator (see Figures \ref{f6app2}, \ref{f6app2ex2} and \ref{f6appex3BSR}) have   no significant effect (see Table \ref{T1}),  when the condition $l_{\alpha} >1/4$ is satisfied under  the bandwidth parameter modelling $B_{T}= T^{-\beta },$ $\beta \in (0,1).$  This fact is also  reflected in Figures \ref{SBTDIVEX1}, \ref{f6app2ex2b} and \ref{f6app2ex3B}, respectively,  where decreasing patterns, and almost the same divergence rates   are displayed  by the sample values of $\mathcal{S}_{B_{T}}$ projected into  $\mathcal{H}_{n}\otimes \mathcal{H}_{n},$  for  $n=1,2,3,4,5,6,7,8,$   in all the  examples. However, the scenario under which   $\alpha (n)$ crosses   the threshold $1/2$  at some spherical scale  $n$   requires a separated analysis, as briefly discussed in Example 4 in the next  section (see Figure \ref{f6app2ex3BS62}).

   \subsection{Example 4}
   \label{Sec5.3}
\begin{figure}[!h]
\begin{center}
\includegraphics[width=10cm,height=11cm]{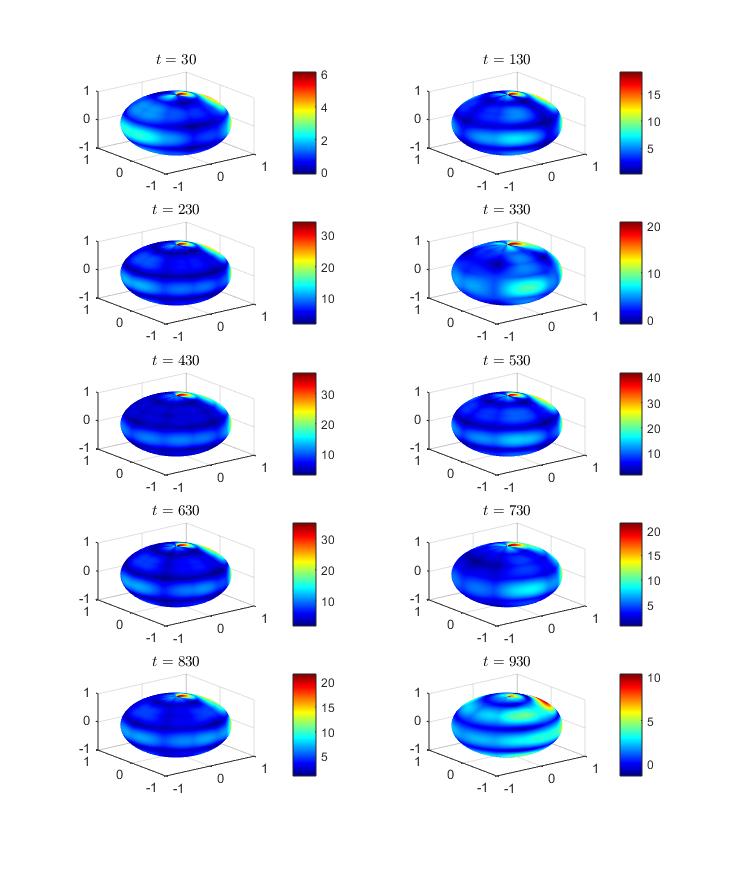}
\end{center}
\caption{\scriptsize{\emph{Example 4}. One sample realization at times $t=30, 130, 230, 330,  430, 530, 630, 730, 830,930$ of   multifractionally integrated  SPHARMA(1,1) process  projected into $\bigoplus_{n=1}^{8}\mathcal{H}_{n}$}}
\label{ex4r1}
\end{figure}

   Our numerical analysis is   extended here  beyond the restriction $L_{\alpha} <1/2.$
   Specifically, this section shows some preliminary numerical results regarding  the effect of higher levels of singularity at zero frequency when  $L_{\alpha }>1/2,$
   i.e., $\|\mathcal{A}\|_{\mathcal{L}(L^{2}(\mathbb{S}_{2},d\nu,\mathbb{C}))}>1/2,$  corresponding to a stronger persistency in time of the projected process into the dominant subspace (see Figure \ref{ex4r1}).

\begin{figure}[!h]
\begin{center}
\includegraphics[width=2.5cm,height=2.5cm]{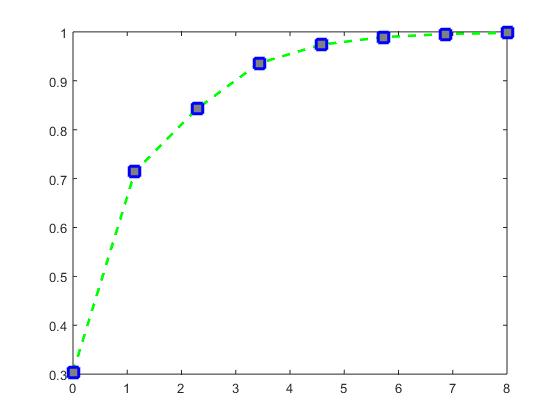}
\includegraphics[width=2.5cm,height=2.5cm]{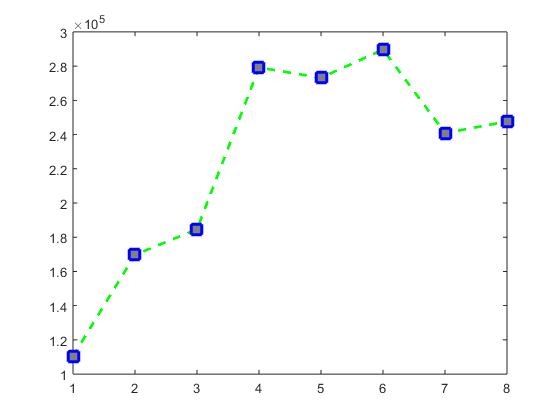}
\includegraphics[width=2.5cm,height=2.5cm]{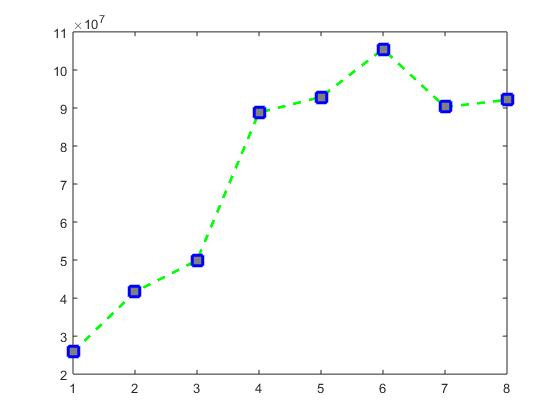}
\includegraphics[width=2.5cm,height=2.5cm]{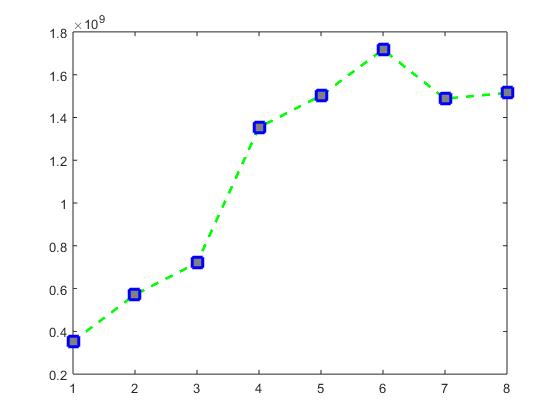}
\end{center}
\caption{\scriptsize{\emph{Example 4}.  Eigenvalues $\alpha (n),$ $n=1,2,3,4,5,6,7,8,$ of  LRD operator $\mathcal{A},$ $L_{\alpha }= 0.9982$ and   $l_{\alpha }=0.3041$ (plot at the left--hand side).
Sample values of the  kernel of the test operator statistic $\mathcal{S}_{B_{T}},$ for  $B_{T}= T^{-1/4},$ projected
into  $\mathcal{H}_{n}\otimes \mathcal{H}_{n},$ $n=1,\dots,8,$  for  the  functional sample sizes $T=1000,10000,30000$ (three plots at the right--hand side) }}
\label{f6app2ex3BS62}
\end{figure}

Under this scenario, conditions (\ref{shconvrm}), and (\ref{ca}) in Theorem \ref{consisimse}, are not satisfied. Indeed, we are out of the scenario where the summability in time of the square of the  Hilbert--Schmidt operator norms  of the elements of the covariance operator family holds. Then, new technical tools are required to address  the asymptotic analysis in the spectral domain of this family of manifold supported functional time series displaying stronger levels of persistency in time.
Let us again consider $\mathcal{S}_{B_{T}},$  for $B_{T}= T^{-1/4},$    projected into $\bigoplus_{n=1}^{8}\mathcal{H}_{n}\otimes \mathcal{H}_{n}.$ In this example,  the multifractional integration of  SPHARMA(1,1) process generated in Section \ref{sgd} has been achieved in terms of LRD operator $\mathcal{A}$ having   eigenvalues displayed at the left--hand side of  Figure  \ref{f6app2ex3BS62}, with $L_{\alpha }= 0.9982$ and   $l_{\alpha }=0.3041,$ and $\mathcal{H}_{8}$  being the dominant subspace.  The same functional sample sizes as in Examples 1--3 have been considered. One can observe, in the  three plots displayed at the right--hand side of  Figure  \ref{f6app2ex3BS62} ,  that the decreasing patterns over $n=1,\dots, 8,$ displayed in Figures \ref{SBTDIVEX1}, \ref{f6app2ex2b} and \ref{f6app2ex3B} do not hold in this example.  Table \ref{T3} also  illustrates a faster increasing than in Examples 1--3 of $\left\|\mathcal{S}_{B_{T}}\right\|_{\mathcal{S}(L^{2}(\mathbb{M}_{d},d\nu, \mathbb{C}))},$ for functional sample sizes ranging from $1000$ to $100000,$  under  $l_{\alpha }>1/4,$ and $B_{T}= T^{-1/4}.$

\begin{table*}
\caption{\textbf{\emph{Example 4. $\left\|\mathcal{S}_{B_{T}}\right\|_{\mathcal{S}(L^{2}(\mathbb{M}_{d},d\nu, \mathbb{C}))},$
} of $\mathcal{S}_{B_{T}}$ projected into $\bigoplus_{n=1}^{8}\mathcal{H}_{n}\otimes \mathcal{H}_{n}$\\  ($\beta =1/4,$ $L_{\alpha }= 0.9982$ and   $l_{\alpha }=0.3041$)   } }
\label{T3}
\begin{tabular}{@{}lrrrrc@{}}
\hline
{\bf Sample}\\
{\bf Size}\\
\hline
 \multicolumn{1}{c}{1000}   &   \multicolumn{1}{c}{5000} &  \multicolumn{1}{c}{10000}& \multicolumn{1}{c}{30000} & \multicolumn{1}{c}{50000 } & \multicolumn{1}{c}{100000} \\
  \hline
   6.5651e+05 &  3.8623e+07 & 2.2172e+08 & 3.5383e+09 & 1.2688e+10  &   7.2258e+10   \\
     \hline
\end{tabular}
\end{table*}
\subsection{Empirical size and power analysis}
\label{sesptp}
 The  empirical size and power properties of the testing approach presented are now illustrated. We have applied the  random projection methodology. Tables \ref{T3bb} and \ref{T4} display the numerical results for $8$ random functional  directions (see
 equation (\ref{ptsopa})).
  Model SPHARMA(1,1) generated in Section \ref{sgd} has been considered in the computation of the empirical size of the test. Multifractionally integrated SPHARMA(1,1) model, generated in Section \ref{ex1}, defines the scenario under the alternative to compute the empirical power. For each one of the eight random directions tested, we have analyzed  the functional samples sizes $T=50,100, 500, 1000,$  and, for each functional sample size, we have considered $R=500,1000, 3000$ repetitions.

The empirical size properties of the proposed  testing  procedure  are quite robust, as one can observe  in the  numerical results displayed  in Table \ref{T3bb}. Specifically, since we are working with finite sample sizes ($T=50,100, 500, 1000$), despite the statistical distance to the normal distribution which holds  asymptotically,   small deviations are observed  from the theoretical $\alpha $ value for all number of repetitions $R$ considered.   One can also observe, in Table \ref{T4}, the increasing patterns displayed by the empirical power with respect to the functional sample sizes tested in all  random directions. Note that these empirical power values    are  in the  interval $[0.776,1].$ In particular,
  since the threshold $T=1000,$ the empirical power is almost $1$  for any of the three  values  of   $R$  studied.

 \begin{table}[!h]
\caption{\textbf{\emph{Empirical size  ($\beta =1/4,$  $\mathbf{k}_{n,j,h,l},$ $n=h=1,2,3,$ $\alpha =0.05$)} }}
\begin{center}
\begin{tabular}{|c|cccccccc|}
\hline
  {\bf R} &  {\bf T} $=$50 & & & && & &  \\ \hline
 500 &0.0280 &    0.0560 &    0.0360&    0.0480 &   0.0600  &  0.0600  &  0.0360  &   0.0520\\
  1000&0.0480 &    0.0420&    0.0320&     0.0380 &   0.0420   & 0.0440   &  0.0320  &   0.0300\\
   3000 &0.0420   & 0.0447 &   0.0453   & 0.0353   & 0.0413   & 0.0400 &   0.0440  &  0.0507\\
  \hline
   {\bf R} &{\bf T} $=$100& & & && & &  \\ \hline
   500&0.0360    & 0.0680  &  0.0440 &   0.0720  &  0.0520  &  0.0240 &   0.0360 &   0.0400\\
  1000&0.0280 &   0.0380  &  0.0380 &   0.0500  &  0.0380  &  0.0740   & 0.0340   & 0.0360\\
   3000 &0.0373  &  0.0507  &  0.0407  &  0.0460  &  0.0440   & 0.0360 &   0.0600  &  0.0480\\
  \hline
  {\bf R} & {\bf T}$=$500 & & & && & & \\ \hline
   500&0.0440 &   0.0520  &  0.0320  &  0.0640 &   0.0320  &  0.0480  &  0.0320  &  0.0480\\
  1000& 0.0420 &   0.0400  &  0.0500   & 0.0460   & 0.0420 &   0.0380  &  0.0540  &  0.0240\\
   3000 & 0.0453 &   0.0393  &  0.0447  &  0.0407  &  0.0427  &  0.0507 &   0.0453  &  0.0553\\
   \hline
  {\bf R} &  {\bf T}$=$1000 & & & && & &  \\ \hline
   500&0.0520  &  0.0360  &  0.0400 &   0.0560  &  0.0600  &  0.0640 &   0.0480   & 0.0520\\
  1000&0.0440  &  0.0380  &  0.0400  &  0.0580  &  0.0500 &   0.0360 &   0.0520 &   0.0400\\
 3000 & 0.0573    & 0.0480   & 0.0507  &  0.0467  &  0.0440  &  0.0453  &  0.0487 &   0.0447\\
   \hline
 \end{tabular}\label{T3bb}
\end{center}
\end{table}

 \begin{table}[!h]
\caption{\textbf{\emph{Empirical power ($\beta =1/4,$  $\mathbf{k}_{n,j,h,l},$ $n=h=1,2,3,$ $\alpha =0.05$) } }}
\begin{center}
\begin{tabular}{|c|cccccccc|}
\hline
   {\bf R} &  {\bf T}$=$50 & & & && & &  \\ \hline
 500&0.9200  &  0.9240 &   0.8640  &  0.8880   & 0.8480 &   0.8160 &   0.7760  &  0.7800\\
  1000& 0.9000  &  0.9000  &  0.8860 &   0.8980 &   0.8000 &   0.8320  &  0.7840  &  0.7840\\
   3000 &0.9247  &   0.9253  &   0.8713  &   0.8760   &  0.8247  &   0.8273   &  0.7980   &  0.8013\\
  \hline
  {\bf R} &{\bf T}$=$100& & & && & &  \\ \hline
   500&0.9920   &  0.9920   &  0.9920   &  0.9880  &   0.9840   &  0.9800  &   0.9880    & 0.9720\\
  1000&0.9880   & 0.9880  &  0.9860  &  0.9840  &  0.9800  &  0.9720  &  0.9740  &  0.9840\\
   3000& 0.9893   & 0.9893  &  0.9920  &  0.9827 &   0.9820  &  0.9773  &  0.9767  &  0.9747\\
  \hline
  {\bf R} & {\bf T}$=$500 & & & && & & \\ \hline
   500&1.0000  &  1.0000   & 1.0000 &   0.9960  &  1.0000  &  0.9960 &   1.0000   & 1.0000\\
  1000&  1.0000 &   1.0000  &  1.0000  &  1.0000  &  0.9980  &  1.0000  &  0.9980   & 1.0000\\
   3000 &1.0000    &  0.9987   &  1.0000  &   1.0000    & 0.9993   &  1.0000    & 0.9993    & 0.9993\\
   \hline
   {\bf R} & {\bf T}$=$1000 & & & && & &  \\ \hline
   500&1 &    1   &  1  &   1    & 1 &    1    & 1     & 1\\
  1000&1 &    1   &  1  &   1    & 1 &    1    & 1     & 1\\
 3000 & 1 &    1   &  1  &   1    & 1 &    1    & 1     & 1\\
   \hline
 \end{tabular}\label{T4}
\end{center}
\end{table}

\section{Final  comments. Reliable inference from our approach}
\label{s6}

The simulation study  illustrates   six key aspects of our approach, briefly summarized in points (i)--(vi) below:
\begin{itemize} \item[(i)] The   tight property  under $H_{0}$ of the random projection sequence \begin{eqnarray}&&\left\langle \sqrt{B_{T}T}(\widehat{\mathcal{F}}_{\omega }^{(T)}-E[\widehat{\mathcal{F}}_{\omega }^{(T)}]),S_{n,j}^{d}\otimes \overline{S_{h,l}^{d}}\right\rangle_{\mathcal{S}(L^{2}(\mathbb{M}_{d},d\nu, \mathbb{C}))}\nonumber\\
&&j=1,\dots, \Gamma(n,d),\ l=1,\dots,\Gamma (h,d),\ n,h\in \mathbb{N}_{0},
\nonumber\end{eqnarray}
\noindent allows   the application of  Prokhorov Theorem to prove the convergence, as $T\to \infty,$  of  $\mathcal{S}_{{B}_{T}}-E[\mathcal{S}_{{B}_{T}}]$ to $\widehat{\mathcal{F}}_{0},$
in the space $\mathcal{L}^{2}_{\mathcal{S}(L^{2}(\mathbb{M}_{d},d\nu, \mathbb{C}))}(\Omega ,\mathcal{A},P).$  In particular,
the asymptotic Gaussian distribution of $\mathcal{S}_{B_{T}}$ under $H_{0}$ follows from this result. This result is  illustrated in Section \ref{sgd}  from Theorem 1.2.1 in \cite{Prato02}.
\item[(ii)]
The  crucial role played by the design of the  test statistic operator  $\mathcal{S}_{B_{T}}$  in the derivation of the  conditions assumed to obtain consistency of the test (see Proposition \ref{cor1} and  Theorem \ref{consisimse}). The simulation study  also  reveals that the additional conditions assumed  in Theorem \ref{th2} lead to    universal a.s. divergence rates.  These rates are  not affected by the localization of the dominant eigenspace, or the value of the parameter $\beta \in (0,1)$ chosen, under the bandwidth parameter scenario $B_{T}=T^{-\beta }$ when $l_{\alpha }>1/4.$
Table \ref{T2} visualizes this fact for parameter values $\beta =0.2, 0.55,  0.9.$ Specifically, in the three examples  analyzed, the  sample values of  $\left\|\frac{\mathcal{S}_{B_{T}}}{(TB_{T})^{1/2}}\right\|_{\mathcal{S}(L^{2}(\mathbb{M}_{d},d\nu, \mathbb{C}))},$
   projected into $\bigoplus_{n=1}^{8}\mathcal{H}_{n}\otimes \mathcal{H}_{n},$  are displayed under these three bandwidth parameter scenarios,  for  functional sample sizes $T=1000,50000,100000.$   Note that, although, as expected,  the sample values of
$\left\|\frac{\mathcal{S}_{B_{T}}}{(TB_{T})^{1/2}}\right\|_{\mathcal{S}(L^{2}(\mathbb{M}_{d},d\nu, \mathbb{C}))}$   slightly increase   when  $\beta $   increases (see Examples 1--3),
   no significant differences are observed in the sample divergence rate of  $\left\|\frac{\mathcal{S}_{B_{T}}}{(TB_{T})^{1/2}}\right\|_{\mathcal{S}(L^{2}(\mathbb{M}_{d},d\nu, \mathbb{C}))},$ between  the three values of parameter  $\beta $ analyzed.  Furthermore,  under condition  $l_{\alpha} >1/4,$ when $B_{T}=T^{-\beta },$ $\beta \in (0,1),$ one can  observe the invariance of the   sample divergence rate against  the location of the dominant eigenspace. Keeping in mind that, as illustrated in  \cite{Ovalle23},  the regular spectral factor $\mathcal{M}_{\omega },$ which here   is represented by the  SPHARMA(p,q) functional spectrum, has not impact in the asymptotic analysis, one can conclude the findings in Table \ref{T2} are representative.
\begin{table}
\caption{\textbf{\emph{$\left\|\frac{\mathcal{S}_{B_{T}}}{(TB_{T})^{1/2}}\right\|_{\mathcal{S}(L^{2}(\mathbb{M}_{d},d\nu, \mathbb{C}))},$
}  $\mathcal{S}_{B_{T}}$ projected into $\bigoplus_{n=1}^{8}\mathcal{H}_{n}\otimes \mathcal{H}_{n}$}}\label{T2}
\begin{tabular}{@{}lcrcrrr@{}}
\hline
 $\boldsymbol{B_{T}=T^{-\beta }}$ &{\bf T}  &
\multicolumn{1}{c}{{\bf Example 1}} &
{\bf Example 2}&
\multicolumn{1}{c}{\bf Example 3} &
 \\
\hline
& $1000$ &   1.6885(1.0e+04) &      1.6384(1.0e+04)     &    1.6438(1.0e+04)\\
$\beta = 0.2$& $50000$ & 4.3549(1.0e+07) & 4.2251(1.0e+07)  &4.2290(1.0e+07)\\
& $100000$ & 1.7493(1.0e+08)   &   1.6944(1.0e+08)    &     1.6989(1.0e+08) \\[6pt]
\hline
& $1000$ &  1.8067 (1.0e+04) &     1.7733(1.0e+04)      &     1.7693(1.0e+04)\\
$\beta =0.55$ & $50000$ & 4.4789(1.0e+07)  & 4.3747(1.0e+07)    & 4.3510(1.0e+07) \\
& $100000$ &    1.7984(1.0e+08)&   1.7470(1.0e+08) &      1.7476(1.0e+08) \\[6pt]
   \hline
& $1000$ &   2.0296(1.0e+04) &      2.0109(1.0e+04)      &   1.9993(1.0e+04) \\
$\beta = 0.9$ & $50000$ & 4.5133(1.0e+07)&   4.4271(1.0e+07) &  4.4138(1.0e+07)\\
& $100000$ & 1.8040 (1.0e+08) & 1.7518(1.0e+08)  &    1.7624(1.0e+08)\\
\hline
\end{tabular}
\end{table}

\item[(iii)] As expected, when higher orders of singularity are displayed at zero frequency, beyond the restriction $L_{\alpha }<1/2,$ a faster divergence of the  sample values of the Hilbert--Schmidt operator norm of $\mathcal{S}_{{B}_{T}},$ and of  its diagonal projections,  is observed.
\item[(iv)] The testing approach adopted shows good empirical size and power properties for finite functional samples, as reported in Tables \ref{T3bb} and \ref{T4}. Namely,
 Table \ref{T3bb}, in  the  simulation study undertaken, shows empirical test sizes  very close to the theoretical value $\alpha =0.05$ for the minimum sample size $T=50$ considered, and   for the number of repetitions $R=500, 1000, 3000.$ Table \ref{T4}  displays, for  $T=50,$ empirical powers in the interval $(0.7760, 0.9253)$  at the eight random directions tested, and for the number of repetitions $R=500, 1000, 3000.$
  Summarizing, as illustrated in Section 5.4, for relatively small functional sample sizes, reliable inference based on our functional spectral nonparametric approach  is  possible. On the other hand, asymptotic properties like consistency of the test are verified to hold, as given in  Table 1, displaying  the increasing order of magnitude of the   Hilbert--Schmidt operator norm of our test statistics, which  is around  $10^{5},$ for the smallest  functional sample size $T=1000$  in all examples.   By the same reasons explained in (ii), i.e., the invariance of the results displayed in Tables \ref{T3bb} and \ref{T4} against the location of the dominant eigenspace ($\mathcal{H}_{1}, \mathcal{H}_{8}$ and $\mathcal{H}_{5}$ respectively in Examples 1,2,3), and  the absence of asymptotic  impact of the choice of the SPHARMA(p,q) functional spectrum, one can conclude  the representativeness of the numerical results  reflected in Tables   \ref{T3bb} and \ref{T4}.

\item[(v)]
    We remark that all computations involved in the simulation study undertaken, in particular, in the implementation of the proposed inference tools and testing approach, have been achieved in terms of a unique orthonormal basis, given by the eigenfunctions of the Laplace Beltrami operator. Thus, we have worked under  the scenario where the eigenfunctions of the elements of the  covariance and spectral density operator families are known. This fact constitutes an important advantage of the analyzed setting, avoiding the use of empirical eigenfunction bases. We have also worked under the  context of fully observed functional data. The case of sparse discretely observed and contaminated functional data can be addressed   from the nonparametric series least--squares estimation of our functional data set, and  plug--in implementation of our test statistics  (see, e.g., \cite{Zhou24} and \cite{RMC2025}). Specifically, under suitable  restrictions on the local  H\"older regularity of our functional data set, and on the  supremum norm of the sieve basis elements, as well as on the pure point spectral properties of the autocovariance matrix of the random  sieve basis (involved in the nonparametric  series  least--squares  estimation of the discretely observed  functional data set), one can derive similar asymptotic results. To this aim, a suitable  manifold random uniform sampling design must be  considered. The time--varying manifold sampling frequency  must display a faster divergence rate than the time--varying sieve basis dimension, but slower than the functional sample size. Under these conditions,  similar  results on   asymptotic $L^{2}$  bias analysis  can be  obtained under $H_{1},$ depending on the almost  surely uniform  convergence rates to the  theoretical values of our functional data, involved in  their nonparametric series least--squares estimation.

\item[(vi)] As commented  in the Introduction,  the presented approach allows   the statistical inference from spatiotemporal  data sets  embedded into the sphere.  Thus, the non--euclidean  spatial statistical analysis of such data sets can be performed enhancing geometrical  interpretability,  avoiding the   usual transformations of longitudes and latitudes of the  data  to work in a cartesian reference  coordinate system. Under the invariance properties assumed in our setting, the case of discretely observed functional data can be addressed considering  the sieve basis constructed from the eigenfunctions of the Laplace Beltrami operator. This sieve basis allows an easier geometrical implementation  (see, e.g., \cite{KurisuMatsuda25}, and \cite{Zhou24} for  alternative sieve bases in nonparametric series least--squares ridge regression  in an euclidean spatial setting).  An   important dimension reduction is obtained in terms of this sieve basis, which is crucial in the reconstruction of high--dimensional data sets. This fact constitutes another remarkable feature of our approach  that reduces   computational burden, allowing the implementation  of resampling techniques  in real data applications. Note that temporal information can be incorporated under our functional time series framework, extending recent developments in the  purely spatial statistical euclidean  context  (see, e.g., \cite{KurisuMatsuda25}).  Just to mention a motivating data example for implementation of our approach,     the authors in \cite{KurisuMatsuda25}  analyze a  georeferenced spatiotemporal discretely observed  population data set to predict  nighttime population in Tokyo. They implement series  spatial ridge regression  estimation after monthly averaging the data, ignoring time information that   is crucial in this prediction problem. Our  functional spectral nonparametric approach allows time information to be processed in an efficient way in an non--euclidean setting as commented in  (iv) and  (v).

\end{itemize}

\noindent \textbf{Acknowledgements}.
This work has been supported in part by projects PID2022--142900NB-I00 and PID2020-116587GB-I00, financed by
MCIU/AEI/10.13039/\linebreak 501100011033 and by FEDER UE,  and CEX2020-001105-M  MCIN/AEI/10.13039/\linebreak 501100011033), as well as supported by grant ED431C 2021/24 (Grupos competitivos) financed by Xunta de Galicia through European Regional Development Funds (ERDF).

\section*{Appendix. Proof of the results}
\begin{appendix}
\section{Proofs of the results in  Section~2}
\subsection{Proof of Theorem \ref{lemvfh0}}
\label{apendicea1}
\begin{proof}
From Lemma 2.1,  \begin{equation}\sqrt{B_{T}T}(\widehat{f}_{\omega_{j} }^{(T)}-E[\widehat{f}_{\omega_{j}}^{(T)}])\to_{D} \widehat{f}_{\omega_{j}},\ j=1,\dots, J,\label{cd}\end{equation}
\noindent where $\to_{D}$ denotes the convergence in distribution. Here,    $\widehat{f}_{\omega_{j} },$ $j=1,\dots,J,$ are jointly zero--mean complex Gaussian elements in $\mathcal{S}(L^{2}(\mathbb{M}_{d},d\nu, \mathbb{C}))$  $\equiv L^{2}(\mathbb{M}_{d}^{2},d\nu\otimes d\nu , \mathbb{C}),$ with covariance kernel (\ref{cobk}).

Let us consider $\{S_{n,j}^{d}, \ j=1,\dots, \Gamma (n,d),\ n\in \mathbb{N}_{0}\},$
 the  orthonormal basis of   eigenfunctions of the Laplace--Beltrami operator $\Delta_{d}$ on $L^{2}\left(\mathbb{M}_{d},d\nu , \mathbb{C}\right).$  From   equation (\ref{cobk}), applying invariance property leading to
 \begin{eqnarray}& &
  f_{\omega }(\tau,\sigma)
\underset{\mathcal{S}(L^{2}(\mathbb{M}_{d},d\nu; \mathbb{C}))}{=}\sum_{n\in \mathbb{N}_{0}}f_{n}(\omega)\sum_{j=1}^{\Gamma (n,d)}S_{n,j}^{d}\otimes \overline{S_{n,j}^{d}}(\tau,\sigma),\quad  (\tau,\sigma)\in \mathbb{M}_{d}^{2},\nonumber\\
  \label{sdo2b}
\end{eqnarray}
\noindent  we obtain, for $\tau_{1},\sigma_{1},\tau_{2},\sigma_{2}\in \mathbb{M}_{d},$
\begin{eqnarray}
&&\mbox{cov}(\widehat{f}_{\omega}(\tau_{1},\sigma_{1}),\widehat{f}_{\omega}(\tau_{2},\sigma_{2}))=2\pi\|W\|_{L^{2}(\mathbb{R})}^{2}\nonumber\\
&&\times \left[\sum_{n,h\in \mathbb{N}_{0}}
\sum_{j=1}^{\Gamma (n,d)}\sum_{l=1}^{\Gamma (h,d)}f_{n}(\omega)f_{h}(\omega)S_{n,j}^{d}(\tau_{1})\overline{S_{h,l}^{d}(\sigma_{1})S_{n,j}^{d}(\tau_{2})}S_{h,l}^{d}(\sigma_{2})\right.\nonumber\\
&&\left.+\eta (2\omega )\sum_{n,h\in \mathbb{N}_{0}}
\sum_{j=1}^{\Gamma (n,d)}\sum_{l=1}^{\Gamma (h,d)}f_{n}(\omega)f_{h}(\omega)  S_{n,j}^{d}(\tau_{1})\overline{S_{h,l}^{d}(\sigma_{1})}S_{h,l}^{d}(\tau_{2})\overline{S_{n,j}^{d}(\sigma_{2})}\right].\nonumber\\
 \label{sdo33}
\end{eqnarray}

Under $H_{0},$ from Theorem D2 in the Supplementary Material of  \cite{Panaretos13},  keeping in mind (\ref{sdo2b}),
\begin{eqnarray}
&&\mbox{cov}\left( \widehat{f}_{\omega}^{(T)}(\tau_{1} ,\sigma_{1}),\widehat{f}_{\omega}^{(T)}(\tau_{2},\sigma_{2})\right)\nonumber\\
&&=\frac{2\pi}{T}\int_{-\pi}^{\pi} W^{(T)}(\omega-\alpha )W^{(T)}(\omega-\alpha )f_{\alpha}(\tau_{1},\tau_{2})\overline{f_{\alpha }(\sigma_{1},\sigma_{2})}d\alpha\nonumber\\
&&+\frac{2\pi}{T}\int_{-\pi}^{\pi}W^{(T)}(\omega-\alpha )W^{(T)}(\omega+\alpha ) f_{\alpha }(\tau_{1},\sigma_{2})\overline{f_{\alpha }(\sigma_{1},\tau_{2})}d\alpha
\nonumber\\
&&+\mathcal{O}(B_{T}^{-2}T^{-2})+\mathcal{O}(T^{-1})\nonumber\\
&&= \frac{2\pi}{T}\sum_{n,h\in \mathbb{N}_{0}}
\sum_{j=1}^{\Gamma (n,d)}\sum_{l=1}^{\Gamma (h,d)}\int_{-\pi}^{\pi} W^{(T)}(\omega-\alpha )W^{(T)}(\omega-\alpha )
f_{n}(\alpha )f_{h}(\alpha)
d\alpha
\nonumber\\
&&\hspace*{4cm}\times S_{n,j}^{d}(\tau_{1})\overline{S_{h,l}^{d}(\sigma_{1})S_{n,j}^{d}(\tau_{2})}S_{h,l}^{d}(\sigma_{2})\nonumber\\
&&
\hspace*{2cm} +\frac{2\pi}{T}\int_{-\pi}^{\pi}W^{(T)}(\omega-\alpha )W^{(T)}(\omega+\alpha )f_{n}(\alpha )f_{h}(\alpha)
d\alpha\nonumber\\
&&\times   S_{n,j}^{d}(\tau_{1})\overline{S_{h,l}^{d}(\sigma_{1})}S_{h,l}^{d}(\tau_{2})
\overline{S_{n,j}^{d}(\sigma_{2})}+\mathcal{O}(B_{T}^{-2}T^{-2})+\mathcal{O}(T^{-1}).
 \label{covkernbb}
\end{eqnarray}

From  (\ref{covkernbb}),    applying Cauchy--Schwartz inequality, and the orthonormality of the basis of eigenfunctions of the Laplace Beltrami operator, the following inequalities hold:  For every $j=1,\dots,\Gamma (n,d),$  $n\in \mathbb{N}_{0},$  and
$l=1,\dots ,\Gamma (h,d),$ $h\in \mathbb{N}_{0},$
\begin{eqnarray}&&
E\left[ \left|\left\langle \sqrt{B_{T}T}(\widehat{\mathcal{F}}_{\omega }^{(T)}-E[\widehat{\mathcal{F}}_{\omega}^{(T)}]),S_{n,j}^{d}\otimes \overline{S_{h,l}^{d}}\right\rangle_{\mathcal{S}(L^{2}(\mathbb{M}_{d},d\nu, \mathbb{C}))}\right|^{2}\right]
\nonumber\\
&& = 2\pi\int_{-\pi}^{\pi}W\left(\frac{\omega -\alpha }{B_{T}}\right)\left[W\left(\frac{\omega-\alpha }{B_{T}}\right)+W\left(\frac{\omega+\alpha }{B_{T}}\right)\left\langle S_{n,j}^{d}, \overline{S_{h,l}^{d}}\right\rangle_{L^{2}(\mathbb{M}_{d}, d\nu,\mathbb{C})}\right.\nonumber\\
&&\left.\times\left\langle \overline{S_{n,j}^{d}}, S_{h,l}^{d} \right\rangle_{L^{2}(\mathbb{M}_{d}, d\nu,\mathbb{C})}\right] f_{n}(\alpha )f_{h}(\alpha )\frac{d\alpha }{B_{T}}+\mathcal{O}(B_{T}^{-2}T^{-2})+\mathcal{O}(T^{-1})\nonumber\\
&&\leq 2\pi\int_{-\pi}^{\pi}W\left(\frac{\omega -\alpha }{B_{T}}\right)\left[W\left(\frac{\omega-\alpha }{B_{T}}\right)+W\left(\frac{\omega+\alpha }{B_{T}}\right)\right]f_{n}(\alpha )f_{h}(\alpha )\frac{d\alpha }{B_{T}}\nonumber\\
&&+\mathcal{O}(B_{T}^{-2}T^{-2})+\mathcal{O}(T^{-1})
\leq \mathcal{N}_{1}\left[\sum_{\tau \in \mathbb{Z}}\|\mathcal{R}_{\tau}\|_{L^{1}(L^{2}(\mathbb{M}_{d}, d\nu, \mathbb{C}))}\right]^{2}+\varepsilon (T)<\infty,
\label{eqtight}
\end{eqnarray}
\noindent  under  SRD,  and for certain   $\mathcal{N}_{1} >0,$ and   $\varepsilon (T)>0,$ with  $\varepsilon (T)\to 0,$
as $T\to \infty.$
As before,  $\widehat{\mathcal{F}}_{\omega }^{(T)}$ denotes the weighted periodogram operator with kernel
 $\widehat{f}^{(T)}_{\omega}.$  In equation (\ref{eqtight}),  we have considered $T$ sufficiently large to apply the identity $W(x)=1/B_{T}W(x/B_{T}),$ for $B_{T}<1,$ and $x\in [-\pi,\pi]$ (see Lemma F11 in the Supplementary Material of \cite{Panaretos13}).  Thus, under $H_{0},$  assuming the  conditions  in Lemma 2.1,  the sequence  $\sqrt{B_{T}T}(\widehat{\mathcal{F}}_{\omega}^{(T)}-E[\widehat{\mathcal{F}}_{\omega}^{(T)}])$  is tight.
  Hence, the convergence, as $T\to \infty,$  of $\sqrt{B_{T}T}(\widehat{\mathcal{F}}_{\omega}^{(T)}-E[\widehat{\mathcal{F}}_{\omega}^{(T)}])$ to the Gaussian random operator  $\widehat{\mathcal{F}}_{\omega}$ with kernel $\widehat{f}_{\omega }$ (see equation (\ref{cd})), in the norm of the space
$\mathcal{L}^{2}_{\mathcal{S}(L^{2}(\mathbb{M}_{d},d\nu, \mathbb{C}))}(\Omega ,\mathcal{A},P),$  follows from Prokhorov Theorem.
Here, \linebreak $\mathcal{L}^{2}_{\mathcal{S}(L^{2}(\mathbb{M}_{d},d\nu; \mathbb{C}))}(\Omega ,\mathcal{A},\mathcal{P})$ denotes the space of zero--mean second--order  \linebreak  $\mathcal{S}(L^{2}(\mathbb{M}_{d},d\nu; \mathbb{C}))$--valued random variables with the norm $\sqrt{E\|\cdot\|^{2}_{\mathcal{S}(L^{2}(\mathbb{M}_{d},d\nu; \mathbb{C}))}}.$

Let us consider \begin{eqnarray}&&C_{n,j,h,l}(\omega ,\xi)\nonumber\\
&&=\mbox{Cov}\left(\sqrt{TB_{T}}\widehat{\mathcal{F}}_{\omega }^{(T)}-\widehat{\mathcal{F}}_{0},\sqrt{TB_{T}}\widehat{\mathcal{F}}_{\xi }^{(T)}-\widehat{\mathcal{F}}_{0}\right)\left(S_{n,j}^{d}\otimes \overline{S_{h,l}^{d}}\right)\left(S_{n,j}^{d}\otimes \overline{S_{h,l}^{d}}\right),\nonumber\end{eqnarray}
\noindent for  $j=1,\dots, \Gamma(n,d),$ and $l=1,\dots,\Gamma (h,d),$  $n,h\in \mathbb{N}_{0},$ where for a bounded linear operator $\mathcal{A}$  on a separable Hilbert space $H,$   $\mathcal{A}(\varphi)(\phi)=\left\langle \mathcal{A}(\varphi),\phi\right\rangle_{H},$ for every $\varphi,\phi\in \mbox{Dom}(\mathcal{A}).$ Here, as before,  $\widehat{\mathcal{F}}_{0}$ is the Gaussian random element with random kernel $\widehat{f}_{0}$ introduced in equation (\ref{cd}) for $\omega_{j}=0.$ Then, applying  Cauchy--Schwartz inequality in $\mathcal{L}^{2}(\Omega ,\mathcal{A},P),$ the space of  complex--valued zero--mean second--order random variables on $(\Omega ,\mathcal{A},P),$ and Jensen's inequality,  we obtain
\begin{eqnarray}
&&E\left\|\left(\mathcal{S}_{B_{T}}-E[\mathcal{S}_{B_{T}}]\right)-\widehat{\mathcal{F}}_{0}\right\|_{\mathcal{S}(L^{2}(\mathbb{M}_{d},d\nu, \mathbb{C}))}^{2}\nonumber\\ &&=\sum_{n,h\in \mathbb{N}_{0}}
\sum_{j=1}^{\Gamma (n,d)}\sum_{l=1}^{\Gamma (h,d)}\int_{[-\pi,\pi]^{2}}\delta_{T}(0-\omega )\delta_{T}(0-\xi)C_{n,j,h,l}(\omega ,\xi)d\omega d\xi\nonumber\\
&&\leq \sum_{n,h\in \mathbb{N}_{0}}
\sum_{j=1}^{\Gamma (n,d)}\sum_{l=1}^{\Gamma (h,d)}\int_{[-\pi,\pi]^{2}}\delta_{T}(0-\omega )\delta_{T}(0-\xi)\sqrt{C_{n,j,h,l}(\omega , \omega)C_{n,j,h,l}(\xi,\xi)}d\omega d\xi\nonumber\\
&&\leq \sum_{n,h\in \mathbb{N}_{0}}
\sum_{j=1}^{\Gamma (n,d)}\sum_{l=1}^{\Gamma (h,d)}
\left[\int_{[-\pi,\pi]^{2}}\delta_{T}(0-\omega )\delta_{T}(0-\xi)C_{n,j,h,l}(\omega , \omega)C_{n,j,h,l}(\xi,\xi)
d\omega d\xi\right]^{1/2}\nonumber\\
&&=\sum_{n,h\in \mathbb{N}_{0}}
\sum_{j=1}^{\Gamma (n,d)}\sum_{l=1}^{\Gamma (h,d)}
\int_{[-\pi,\pi]}\delta_{T}(0-\omega )C_{n,j,h,l}(\omega , \omega)d\omega
\nonumber\\
&&=\int_{[-\pi,\pi]}\delta_{T}(0-\omega )E\left\|\sqrt{TB_{T}}\left(\widehat{\mathcal{F}}_{\omega }^{(T)}- E\left[\widehat{\mathcal{F}}_{\omega }^{(T)}\right]\right)-\widehat{\mathcal{F}}_{0}\right\|^{2}_{\mathcal{S}(L^{2}(\mathbb{M}_{d},d\nu, \mathbb{C}))}d\omega \nonumber\\
&&\leq \mathcal{N}_{2}\left[\left(\sum_{\tau\in \mathbb{Z}}\left\|\mathcal{R}_{\tau }\right\|_{L^{1}(L^{2}(\mathbb{M}_{d},d\nu, \mathbb{C}))}
\right)^{2}+\left(\sum_{\tau\in \mathbb{Z}}\left\|\mathcal{R}_{\tau }\right\|_{L^{1}(L^{2}(\mathbb{M}_{d},d\nu, \mathbb{C}))}\right)^{4}\right]<\infty,
\nonumber\\\label{compshnorm2bb}
\end{eqnarray}
\noindent  under $H_{0},$ for certain positive constant $\mathcal{N}_{2},$ where the last inequality follows from equations (\ref{sdo33}) and (\ref{covkernbb}), for $T$ sufficiently large.

Applying Dominated Convergence Theorem, we then obtain
\begin{eqnarray}&&
\lim_{T\to \infty}E\left\|\left(\mathcal{S}_{B_{T}}-E[\mathcal{S}_{B_{T}}]\right)-\widehat{\mathcal{F}}_{0}\right\|_{\mathcal{S}(L^{2}(\mathbb{M}_{d},d\nu, \mathbb{C}))}^{2}\nonumber\\
&&=\lim_{T\to \infty}\sum_{n,h\in \mathbb{N}_{0}}
\sum_{j=1}^{\Gamma (n,d)}\sum_{l=1}^{\Gamma (h,d)}\int_{[-\pi,\pi]^{2}}\delta_{T}(0-\omega)\delta_{T}(0-\xi)\nonumber\\ &&\times \mbox{Cov}\left(\sqrt{TB_{T}}\widehat{\mathcal{F}}_{\omega }^{(T)}-\widehat{\mathcal{F}}_{0},\sqrt{TB_{T}}\widehat{\mathcal{F}}_{\xi }^{(T)}-\widehat{\mathcal{F}}_{0}\right)\left(S_{n,j}^{d}\otimes \overline{S_{h,l}^{d}}\right)\left(S_{n,j}^{d}\otimes \overline{S_{h,l}^{d}}\right)d\omega d\xi\nonumber\\
&&=\sum_{n,h\in \mathbb{N}_{0}}
\sum_{j=1}^{\Gamma (n,d)}\sum_{l=1}^{\Gamma (h,d)}\lim_{T\to \infty}\int_{[-\pi,\pi]^{2}}\delta_{T}(0-\omega)\delta_{T}(0-\xi)\nonumber\\ &&\hspace*{-1cm}\times \mbox{Cov}\left(\sqrt{TB_{T}}\widehat{\mathcal{F}}_{\omega }^{(T)}-\widehat{\mathcal{F}}_{0},\sqrt{TB_{T}}\widehat{\mathcal{F}}_{\xi }^{(T)}-\widehat{\mathcal{F}}_{0}\right)\left(S_{n,j}^{d}\otimes \overline{S_{h,l}^{d}}\right)\left(S_{n,j}^{d}\otimes \overline{S_{h,l}^{d}}\right)d\omega d\xi=0, \nonumber\\
\label{idopk2}
 \end{eqnarray}
\noindent in view of the  convergence in $\mathcal{L}^{2}_{\mathcal{S}(L^{2}(\mathbb{M}_{d},d\nu; \mathbb{C}))}(\Omega,\mathcal{A},P)$ of $\sqrt{TB_{T}}\left(\widehat{\mathcal{F}}_{0}^{(T)}-E\left[\widehat{\mathcal{F}}_{0}^{(T)}\right]\right)$ to
$\widehat{\mathcal{F}}_{0}.$ Thus, the convergence in distribution of $\mathcal{S}_{B_{T}}-E[\mathcal{S}_{B_{T}}]$ to $\widehat{\mathcal{F}}_{0}$
holds.

 \end{proof}

\section{Proofs of the results in Section~3}

\subsection{Proof of Lemma \ref{lem1}}
\label{ap2}
\begin{proof}
 Let us consider
  \begin{eqnarray}&&
 \hspace*{-1cm}\left\|\int_{-\pi}^{\pi}\left[\mathcal{F}_{\omega }-\mathcal{F}^{(T)}_{\omega }\right]d\omega\right\|^{2}_{\mathcal{S}(L^{2}(\mathbb{M}_{d},d\nu, \mathbb{C}))}
 \nonumber\\ &&=\sum_{n\in \mathbb{N}_{0}}\sum_{j=1}^{\Gamma (n,d)}\int_{[-\pi,\pi]^{2}}\left[f_{n}(\omega )-f_{n}^{(T)}(\omega )\right]
 \overline{f_{n}(\xi )}d\xi d\omega \nonumber\\ &&
 +\int_{[-\pi,\pi]^{2}}\left[f_{n}^{(T)}(\omega )-f_{n}(\omega )\right] \overline{f_{n}^{(T)}(\xi )}d\xi d\omega,
  \label{orderofconvmpo2}
   \end{eqnarray}
\noindent where the sequence of functions
$$f_{n}^{(T)}(\omega )=\int_{-\pi}^{\pi} F_{T}(\omega - \xi)
f_{n}(\xi) d\xi,\quad \forall \omega \in [-\pi,\pi],\ n\in \mathbb{N}_{0},$$
\noindent defines the  frequency--varying  pure point spectra of the operator family \linebreak $\left\{\mathcal{F}_{\omega}^{(T)}=E_{H_{1}}\left[\mathcal{P}_{ \omega }^{(T)}\right],\ \omega \in [-\pi,\pi]\right\},$
for every $T\geq 2,$ with $\mathcal{P}_{ \omega }^{(T)}$ denoting the periodogram operator (see equation (\ref{periodogro})), and  $E_{H_{1}}$ denoting  the  expectation under the  alternative $H_{1}.$

Under $H_{1},$ for every $n\in \mathbb{N}_{0},$   $f_{n}(\cdot)\in L^{1}([-\pi,\pi]).$ Applying well--known properties of F\'ejer kernel, we then obtain,
as $T\to \infty,$
 \begin{eqnarray}&&
 f_{n}^{(T)}(\omega )\to  f_{n}(\omega ),\quad \forall \omega \in [-\pi,\pi]\backslash \Lambda_{0},\ \mbox{with} \ \int_{\Lambda_{0}}d\omega =0.
 \label{fkp}
 \end{eqnarray}
 \noindent  From Young convolution inequality in $L^{2}([-\pi,\pi]),$ for each  $n\in \mathbb{N}_{0},$ and  $T\geq 2,$
   \begin{equation}
 \int_{[-\pi,\pi]}\left|f_{n}^{(T)}(\omega )\right|^{2}d\omega \leq  \int_{[-\pi,\pi]}\left|f_{n}(\omega )\right|^{2}d\omega <\infty,
 \label{ycitc}
 \end{equation}
 \noindent  under $H_{1},$   since $l_{\alpha },L_{\alpha }\in (0,1/2).$  Thus, $f_{n}^{(T)} \in L^{2}([-\pi,\pi]),$ for every $n\in \mathbb{N}_{0},$ and $T.$

   To apply Dominated Convergence Theorem in (\ref{orderofconvmpo2}), the following additional  inequalities are considered, obtained from triangle inequality,   Young convolution inequality for functions in $L^{1}([-\pi,\pi]),$ and  Jensen's   inequality, keeping in mind that $f_{n}(\omega )\geq 0,$ a.s. in $\omega \in [-\pi,\pi],$ for every $n\in \mathbb{N}_{0},$
   \begin{eqnarray}
   &&\left|\int_{[-\pi,\pi]^{2}}\overline{f_{n}^{(T)}(\xi )}f_{n}^{(T)}(\omega )d\xi d\omega \right|
      \leq
   \int_{[-\pi,\pi]^{2}}\left|\overline{f_{n}^{(T)}(\xi )}f_{n}^{(T)}(\omega )\right|d\xi d\omega \nonumber\\
   &&\leq
    \int_{[-\pi,\pi]^{2}}\left|f_{n}(\xi )f_{n}(\omega )\right|d\xi d\omega =\left[\int_{[-\pi,\pi]}f_{n}(\omega )d\omega \right]^{2}\leq \int_{[-\pi,\pi]}\left|f_{n}(\omega )\right|^{2}d\omega .\nonumber\\
        \label{orderofconvmpo3}
   \end{eqnarray}

   Also, in a similar way,
   \begin{eqnarray}
   &&\left|\int_{[-\pi,\pi]^{2}}\overline{f_{n}^{(T)}(\xi )}f_{n}(\omega )d\xi d\omega \right|\leq \int_{[-\pi,\pi]}\left|f_{n}(\omega )\right|^{2}d\omega
   \nonumber\\
   &&\left|\int_{[-\pi,\pi]^{2}}\overline{f_{n}(\xi )}f_{n}^{(T)}(\omega )d\xi d\omega \right|\leq \int_{[-\pi,\pi]}\left|f_{n}(\omega )\right|^{2}d\omega
    \nonumber\\
   &&\left|\int_{[-\pi,\pi]^{2}}f_{n}(\xi )f_{n}(\omega )d\xi d\omega \right|\leq \int_{[-\pi,\pi]}\left|f_{n}(\omega )\right|^{2}d\omega .
      \label{efproofl2}
   \end{eqnarray}

   Under $H_{1},$  from  equation (\ref{shconvrm}),
   \begin{equation}\sum_{n\in \mathbb{N}_{0}}\Gamma (n,d)\int_{[-\pi,\pi]}\left|f_{n}(\omega )\right|^{2}d\omega =\int_{[-\pi,\pi]}\left\|\mathcal{F}_{\omega }\right\|^{2}_{\mathcal{S}(L^{2}(\mathbb{M}_{d},d\nu, \mathbb{C}))}d\omega <\infty.\label{efproofl2bb}\end{equation}
   \noindent From equations (\ref{fkp})--(\ref{efproofl2bb}), one can apply Dominated Convergence Theorem in equation (\ref{orderofconvmpo2}),
     obtaining
   \begin{eqnarray}&&
   \hspace*{-1cm}\lim_{T\to \infty}
   \left\|\int_{-\pi}^{\pi}\left[\mathcal{F}_{\omega }-\mathcal{F}^{(T)}_{\omega }\right]d\omega\right\|^{2}_{\mathcal{S}(L^{2}(\mathbb{M}_{d},d\nu, \mathbb{C}))}\nonumber\\
    &&=\sum_{n\in \mathbb{N}_{0}}\sum_{j=1}^{\Gamma(n,d)}\int_{[-\pi,\pi]^{2}}\lim_{T\to \infty}\overline{f_{n}^{(T)}(\xi )}\left[f_{n}^{(T)}(\omega )-f_{n}(\omega )\right] d\xi d\omega \nonumber\\
 &&+\int_{[-\pi,\pi]^{2}} \lim_{T\to \infty}\overline{f_{n}(\xi )}\left[f_{n}(\omega )-f_{n}^{(T)}(\omega )\right]d\xi d\omega =0.
 \label{efproofl3}
   \end{eqnarray}

The rate of convergence to zero of the bias is now obtained   in the time domain.   Let $B_{n}$  be defined as \begin{equation}B_{n}(t)=\int_{-\pi}^{\pi}\exp(it\omega )f_{n}(\omega )d\omega ,\quad t\in \mathbb{Z},\ n\in \mathbb{N}_{0}.
   \label{eqift}
   \end{equation}
   \noindent
   The function sequence  $$\left\{ \mathbb{I}_{[-(T-1),T-1]}(t)\frac{T-|t|}{T}B_{n}(t),\ t\in \mathbb{Z},\   n\in \mathbb{N}_{0}\right\}_{T\geq 2}$$ \noindent  pointwise  converges, as $T\to \infty ,$ to $B_{n}(t)$ with rate of convergence $T^{-1},$ and  satisfies,
  for every  $T\geq 2,$
   \begin{eqnarray}&&
   \left|  \mathbb{I}_{[-(T-1),T-1]}(t)\frac{T-|t|}{T}B_{n}(t)\right|^{2}\leq  \left|B_{n}(t)\right|^{2}.\label{dfsint}\end{eqnarray}

  From  (\ref{efproofl2bb}) and   Parseval identity (see equation (\ref{eqift})),
    \begin{eqnarray}
    && \sum_{t\in \mathbb{Z}}\sum_{n\in \mathbb{N}_{0}}\Gamma (n,d)\left|B_{n}(t)\right|^{2}=\sum_{t\in \mathbb{Z}}\|\mathcal{R}_{t}\|^{2}_{\mathcal{S}(L^{2}(\mathbb{M}_{d},d\nu, \mathbb{C}))}\nonumber\\
    &&\hspace*{3.7cm}
    =\int_{-\pi}^{\pi}\left\|\mathcal{F}_{\omega }\right\|_{\mathcal{S}(L^{2}(\mathbb{M}_{d},d\nu, \mathbb{C}))}^{2}d\omega
        <\infty.\nonumber\\
   \label{dtchhbb}
   \end{eqnarray}
   From equations (\ref{dfsint}) and  (\ref{dtchhbb}), Dominated Convergence Theorem then  leads to
   \begin{eqnarray}
   &&\lim_{T\to \infty}\sum_{t\in \mathbb{Z}}\left\|\mathcal{R}_{t}-\mathbb{I}_{[-(T-1),T-1]}(t)\frac{T-|t|}{T}\mathcal{R}_{t}\right\|^{2}_{\mathcal{S}(L^{2}(\mathbb{M}_{d},d\nu, \mathbb{C}))}\nonumber\\
   &&=\sum_{t\in \mathbb{Z}}
  \sum_{n\in \mathbb{N}_{0}}\Gamma (n,d)\lim_{T\to \infty} \left| B_{n}(t)- \mathbb{I}_{[-(T-1),T-1]}(t)\frac{T-|t|}{T}B_{n}(t)\right|^{2}=  0,\nonumber\\ \label{dtc}
   \end{eqnarray}
        \noindent and  $ \sum_{t\in \mathbb{Z}}
 \left\|\mathcal{R}_{t}-\mathbb{I}_{[-(T-1),T-1]}(t)\frac{T-|t|}{T}\mathcal{R}_{t}\right\|^{2}_{\mathcal{S}(L^{2}(\mathbb{M}_{d},d\nu, \mathbb{C}))}=\mathcal{O}(T^{-2}).$ Hence, the desired result follows  from Parseval identity.

  \end{proof}
\subsection{Proof of Corollary \ref{lem3}}
\label{aplem3}

\begin{proof}

Applying  Lemma 3.1, and Lemmas  F10 and F12 of Appendix F in the Supplementary Material of    \cite{Panaretos13},
\begin{eqnarray}&&
\int_{-\pi }^{\pi}E_{H_{1}}[\widehat{\mathcal{F}}^{(T)}_{\omega }]d\omega \nonumber\\
&&
\underset{\mathcal{S}(L^{2}(\mathbb{M}_{d},d\nu, \mathbb{C}))}{=}\int_{\mathbb{R}}W(\xi)\int_{-\pi }^{\pi}\left[
\mathcal{F}_{\omega -\xi B_{T}}+\mathcal{O}(T^{-1})\right]d\omega d\xi+\mathcal{O}(B_{T}^{-1}T^{-1})\nonumber\\
&&\underset{\mathcal{S}(L^{2}(\mathbb{M}_{d},d\nu, \mathbb{C}))}{=}\int_{-\pi }^{\pi}\int_{\mathbb{R}}W(\xi)
\mathcal{F}_{\omega -\xi B_{T}}d\xi d\omega+\mathcal{O}(T^{-1})+\mathcal{O}(B_{T}^{-1}T^{-1}),\label{ec23}
\end{eqnarray}

\noindent as  we wanted to prove.

\end{proof}

\subsection{Proof of Lemma \ref{lem4cs}}
\label{a3b}
\begin{proof}

Under \emph{Assumption I},   there exists an orthonormal basis $\left\{ \phi_{n},\ n\in \mathbb{N}\right\}$ of  \linebreak $L^{2}\left(\mathbb{M}_{d}^{2},\otimes_{i=1}^{2}\nu(dx_{i}), \mathbb{R}\right)$
  such that  (see \cite{Helgason59})
 \begin{eqnarray}
 &&\int_{\mathbb{M}_{d}}\mbox{cum}\left(X_{u_{1}},  X_{u_{2}},
 X_{u_{3}},  X_{0}\right)(\tau_{1}, \tau_{2}, \tau_{3}, \tau_{4})\phi_{n} (\tau_{3},\tau_{4})\nu(d\tau_{3})\nu(d\tau_{4})\nonumber\\
 &&=B_{n}(u_{1},u_{2},u_{3})\phi_{n} (\tau_{1},\tau_{2}),\quad \forall (\tau_{1},\tau_{2})\in \mathbb{M}_{d}\times\mathbb{M}_{d},\ u_{1},u_{2},u_{3}\in \mathbb{Z},\quad n\geq 1.\nonumber\\
 \label{eql6p1}
 \end{eqnarray}

Furthermore,
 \begin{eqnarray}
 &&
\int_{[-\pi ,\pi ]^{3}}
\mbox{cum}\left(\widetilde{X}_{\omega_{1}}^{(T)}(\tau_{1}), \widetilde{X}_{\omega_{2}}^{(T)}(\tau_{2}),
\widetilde{X}_{\omega_{3}}^{(T)}(\tau_{3}), \widetilde{X}_{\omega_{4}}^{(T)}(\tau_{4})\right)
d\omega_{1}d\omega_{2}d\omega_{3}\nonumber\\
&&\underset{\mathcal{S}\left(L^{2}\left(\mathbb{M}_{d}^{2},\otimes_{i=1}^{2}\nu(dx_{i}), \mathbb{C}\right)\right)}{=}
\frac{1}{(2\pi T)^{2}}\int_{[-\pi ,\pi ]^{3}}\sum_{t_{1},t_{2},t_{3},t_{4}=0}^{T-1}\exp\left( -i\sum_{j=1}^{3}(t_{j}-t_{4})\omega_{j}\right)
\nonumber\\
&&
\times \exp\left(-it_{4}\sum_{j=1}^{4}\omega_{j}\right) \mbox{cum}\left(X_{t_{1}-t_{4}}(\tau_{1}),  X_{t_{2}-t_{4}}(\tau_{2}),
 X_{t_{3}-t_{4}}(\tau_{3}),  X_{0}(\tau_{4})\right)\prod_{j=1}^{3}d\omega_{j}\nonumber\\
&&\underset{\mathcal{S}\left(L^{2}\left(\mathbb{M}_{d}^{2},\otimes_{i=1}^{2}\nu(dx_{i}), \mathbb{C}\right)\right)}{=}\int_{[-\pi ,\pi ]^{3}}\frac{1}{(2\pi T)^{2}}\sum_{u_{1},u_{2},u_{3}=-(T-1)}^{T-1}\exp\left( -i\sum_{j=1}^{3}u_{j}\omega_{j}\right)
 \nonumber\\
 &&
 \hspace*{2cm} \times \mbox{cum}\left(X_{u_{1}}(\tau_{1}),  X_{u_{2}}(\tau_{2}),
 X_{u_{3}}(\tau_{3}),  X_{0}(\tau_{4})\right)\nonumber\\
 &&\hspace*{-0.5cm}\times \sum_{t\in \mathbb{Z}}h^{(T)}(u_{1}+t)h^{(T)}(u_{2}+t)h^{(T)}(u_{3}+t)h^{(T)}(t)\exp\left(-it\left(\sum_{j=1}^{4}\omega_{j}\right)\right)
 \prod_{j=1}^{3}d\omega_{j},\nonumber\\
 \label{eql6p2}
 \end{eqnarray}
 \noindent with $h(t)=1,$ $0\leq t\leq T,$ and $h(t)=0,$ otherwise. In (\ref{eql6p2}),  we have considered the change of variable $u_{j}=t_{j}-t_{4},$  $j=1,2,3,$
and $t=t_{4}.$
 Denote, for every $n\geq 1,$  and $(\omega_{1},\omega_{2},\omega_{3})\in [-\pi,\pi]^{3},$ $$f_{n}(\omega_{1},\omega_{2},\omega_{3})=\frac{1}{(2\pi)^{3}}\sum_{u_{1},u_{2},u_{3}\in \mathbb{Z}}\exp\left( -i\sum_{j=1}^{3}\omega_{j}u_{j}\right)B_{n}(u_{1},u_{2},u_{3}),$$

 \noindent where, for  $u_{1},u_{2},u_{3}\in \mathbb{Z},$  $\{B_{n}(u_{1},u_{2},u_{3}),\ n\geq 1\}$ satisfies (\ref{eql6p1}). From equations  (\ref{eql6p1})--(\ref{eql6p2}),
applying Fourier transform inversion formula, for each
  $n\geq 1,$
 \begin{eqnarray}
 &&
T \int_{[-\pi ,\pi ]^{3}\times \mathbb{M}^{4}_{d}}
\mbox{cum}\left(\widetilde{X}_{\omega_{1}}^{(T)}(\tau_{1}), \widetilde{X}_{\omega_{2}}^{(T)}(\tau_{2}),
\widetilde{X}_{\omega_{3}}^{(T)}(\tau_{3}), \widetilde{X}_{\omega_{4}}^{(T)}(\tau_{4})\right)\nonumber\\
&&\hspace*{4cm}\times \phi_{n}(\tau_{1},\tau_{2})\phi_{n}(\tau_{3},\tau_{4})
\prod_{j=1}^{4}d\tau_{j} \prod_{i=1}^{3}d\omega_{i}\nonumber\\
&&=\frac{(2\pi)^{3}}{(2\pi )^{2} T}\int_{[-\pi,\pi]^{6}}\sum_{u_{1},u_{2},u_{3}=-(T-1)}^{T-1}\exp\left(-i\sum_{j=1}^{3}u_{j}(\omega_{j}-\xi_{j})\right)\nonumber\\
&&\times
\sum_{t\in \mathbb{Z}}h^{(T)}(u_{1}+t)h^{(T)}(u_{2}+t)h^{(T)}(u_{3}+t)h^{(T)}(t)\exp\left(-it\left(\sum_{j=1}^{4}\omega_{j}\right)\right)\nonumber\\
&&\hspace*{3cm}\times f_{n}(\xi_{1},\xi_{2},\xi_{3})\prod_{j=1}^{3}d\xi_{j}\prod_{i=1}^{3}d\omega_{i}\nonumber\\
&&=\frac{2\pi}{ T}\int_{[-\pi,\pi]^{6}}\left[\sum_{u_{1},u_{2},u_{3}=-(T-1)}^{T-1}\exp\left(-i\sum_{j=1}^{3}u_{j}(\omega_{j}-\xi_{j})\right)\right.\nonumber\\
&&\left.\times \sum_{t\in \mathbb{Z}}\exp\left(-it\left(\sum_{j=1}^{4}\omega_{j}\right)\right)
h\left(t+\max_{j=1,2,3}|u_{j}|\right)\right]f_{n}\left(\xi_{1},\xi_{2},\xi_{3}\right)\prod_{j=1}^{3}d\xi_{j}\prod_{i=1}^{3}d\omega_{i}.\nonumber\\
\label{eql6p3b}
 \end{eqnarray}

As  $T\to \infty,$ uniformly in $\omega_{4}\in [-\pi,\pi],$\begin{eqnarray}
&&\frac{1}{T}\left[\sum_{u_{1},u_{2},u_{3}=-(T-1)}^{T-1}\exp\left(-i\sum_{j=1}^{3}u_{j}(\omega_{j}-\xi_{j})\right)\right.\nonumber\\
&&\hspace*{2cm}\left.
\times \sum_{t\in \mathbb{Z}}\exp\left(-it\left(\sum_{j=1}^{4}\omega_{j}\right)\right)
h\left(t+\max_{j=1,2,3}|u_{j}|\right)\right]\to \delta (\boldsymbol{\omega}-\boldsymbol{\xi}),\nonumber\\
\label{diracd}
\end{eqnarray}
\noindent where $\delta (\boldsymbol{\omega}-\boldsymbol{\xi})=\prod_{j=1}^{3}\delta (\omega_{j}-\xi_{j})$ denotes the Dirac Delta distribution, defining the  kernel of the identity operator on $L^{2}([-\pi,\pi]^{3}).$ Using the notation
 \begin{eqnarray}
 &&\delta_{T} (\boldsymbol{\omega}-\boldsymbol{\xi}):=
\frac{1}{T}\left[\sum_{u_{1},u_{2},u_{3}=-(T-1)}^{T-1}\exp\left(-i\sum_{j=1}^{3}u_{j}(\omega_{j}-\xi_{j})\right)\right.\nonumber\\
&&\hspace*{3cm}\left.
\times \sum_{t\in \mathbb{Z}}\exp\left(-it\left(\sum_{j=1}^{4}\omega_{j}\right)\right)
h\left(t+\max_{j=1,2,3}|u_{j}|\right)\right],
\label{diracd}
\end{eqnarray}
\noindent equation (\ref{eql6p3b})  can be rewritten as
\begin{eqnarray}
 &&
T \int_{[-\pi ,\pi ]^{3}\times \mathbb{M}^{4}_{d}}
\mbox{cum}\left(\widetilde{X}_{\omega_{1}}^{(T)}(\tau_{1}), \widetilde{X}_{\omega_{2}}^{(T)}(\tau_{2}),
\widetilde{X}_{\omega_{3}}^{(T)}(\tau_{3}), \widetilde{X}_{\omega_{4}}^{(T)}(\tau_{4})\right)\nonumber\\
&&\hspace*{4cm}\times \phi_{n}(\tau_{1},\tau_{2})\phi_{n}(\tau_{3},\tau_{4})
\prod_{j=1}^{4}d\tau_{j} \prod_{i=1}^{3}d\omega_{i}\nonumber\\
&&=2\pi \int_{[-\pi,\pi]^{6}}
\delta_{T} (\boldsymbol{\omega}-\boldsymbol{\xi})f_{n}\left(\boldsymbol{\xi}\right)d\boldsymbol{\xi}d\boldsymbol{\omega },\quad
n\geq 1.
\label{eql6p3}
 \end{eqnarray}
 Note that, for  $T\geq T_{0},$  with $T_{0}$ sufficiently large,
\begin{eqnarray}
 &&
\left|\delta_{T} (\boldsymbol{\omega}-\boldsymbol{\xi})f_{n}(\boldsymbol{\xi})\right|\leq \left|f_{n}(\boldsymbol{\xi})\right|,
\quad   \boldsymbol{\omega}\neq \boldsymbol{\xi},
\label{neq}
\end{eqnarray}
\noindent since $\delta_{T} (\boldsymbol{\omega}-\boldsymbol{\xi})\to 0,$ $T\to \infty,$ for every
$(\boldsymbol{\omega}, \boldsymbol{\xi})\in [-\pi,\pi]^{6}\backslash\Lambda,$ with $\Lambda =\{ (\boldsymbol{\omega}, \boldsymbol{\xi})\in [-\pi,\pi]^{6};\ \boldsymbol{\omega}= \boldsymbol{\xi}
\}\subset [-\pi,\pi]^{6}.$  Under \emph{Assumption I}, applying Parseval identity,
\begin{eqnarray}
&&\sum_{n\geq 1}\int_{[-\pi,\pi]^{6}}|f_{n}(\omega_{1},\omega_{2},\omega_{3})|\prod_{j=1}^{6}d\omega_{j}\nonumber\\
&&\leq (2\pi)^{3}\sum_{n\geq 1}\int_{[-\pi,\pi]^{3}}|f_{n}(\omega_{1},\omega_{2},\omega_{3})|^{2}\prod_{j=1}^{3}d\omega_{j}
\nonumber\\
&&
= (2\pi)^{3}\int_{[-\pi,\pi]^{3}}\|\mathcal{F}_{\omega_{1},\omega_{2},\omega_{3}}\|^{2}_{\mathcal{S}\left(L^{2}\left(\mathbb{M}_{d}^{2},\otimes_{i=1}^{2}\nu(dx_{i}), \mathbb{C}\right)\right)}\prod_{j=1}^{3}d\omega_{j}\nonumber\\
&&=(2\pi)^{3}\sum_{t_{1},t_{2},t_{3}\in
\mathbb{Z}}\left\|\mbox{cum}\left(X_{t_{1}}, X_{t_{2}}, X_{t_{3}},  X_{0}\right)\right\|_{L^{2}(\mathbb{M}_{d}^{4},\otimes_{i=1}^{4} d\nu(x_{i}), \mathbb{R})}^{2}<\infty.
\label{piproof}
\end{eqnarray}
\noindent  Hence, from (\ref{neq})--(\ref{piproof}), applying  Dominated Convergence Theorem,
\begin{eqnarray}
&&\lim_{T\to \infty}\int_{[-\pi,\pi]^{6}}\delta_{T} (\boldsymbol{\omega}-\boldsymbol{\xi})f_{n}\left(\boldsymbol{\xi}\right)d\boldsymbol{\xi}d\boldsymbol{\omega }\nonumber\\
&&=\int_{[-\pi,\pi]^{6}}\lim_{T\to \infty}\delta_{T} (\boldsymbol{\omega}-\boldsymbol{\xi})f_{n}\left(\boldsymbol{\xi}\right)d\boldsymbol{\xi}d\boldsymbol{\omega }\nonumber\\
&&=\int_{[-\pi,\pi]^{3}}f_{n}\left(\boldsymbol{\omega}\right)d\boldsymbol{\omega },\ n\geq 1.
\label{iddd3}
\end{eqnarray}
\noindent and, as $T\to \infty,$  \begin{eqnarray}
&&\left|\int_{[-\pi,\pi]^{6}}\delta_{T} (\boldsymbol{\omega}-\boldsymbol{\xi})f_{n}\left(\boldsymbol{\xi}\right)d\boldsymbol{\xi}d\boldsymbol{\omega }-\int_{[-\pi,\pi]^{3}}f_{n}\left(\boldsymbol{\omega}\right)d\boldsymbol{\omega }\right|=\mathcal{O}(T^{-1}).\label{eqbb}
\end{eqnarray}

Therefore,  from (\ref{eql6p3}), (\ref{piproof}) and (\ref{eqbb}), uniformly in $\omega_{4}\in [-\pi,\pi],$
\begin{eqnarray}
&& \lim_{T\to \infty}\sum_{n}
\left|T \int_{[-\pi ,\pi ]^{3}\times \mathbb{M}^{4}_{d}}
\mbox{cum}\left(\widetilde{X}_{\omega_{1}}^{(T)}(\tau_{1}), \widetilde{X}_{\omega_{2}}^{(T)}(\tau_{2}),
\widetilde{X}_{\omega_{3}}^{(T)}(\tau_{3}), \widetilde{X}_{\omega_{4}}^{(T)}(\tau_{4})\right)\right.\nonumber\\
&&\hspace*{2cm}\left.\times \phi_{n}(\tau_{1},\tau_{2})\phi_{n}(\tau_{3},\tau_{4})
\prod_{j=1}^{4}d\tau_{j}  -2\pi f_{n}\left(\omega_{1},\omega_{2},\omega_{3}\right)\prod_{i=1}^{3}d\omega_{i}\right|
=0.\nonumber
\end{eqnarray}

\noindent It then follows that, as $T\to \infty,$ the  norm

\begin{eqnarray}&&
\hspace*{-1.25cm}\left\| \int_{[-\pi ,\pi ]^{3}}
\left[T\mbox{cum}\left(\widetilde{X}_{\omega_{1}}^{(T)}, \widetilde{X}_{\omega_{2}}^{(T)},
\widetilde{X}_{\omega_{3}}^{(T)}, \widetilde{X}_{\omega_{4}}^{(T)}\right)
-2\pi\mathcal{F}_{\omega_{1},\omega_{2},\omega_{3}}\right]\prod_{j=1}^{3}d\omega_{j}
\right\|_{\mathcal{S}\left(L^{2}\left(\mathbb{M}_{d}^{2},\otimes_{i=1}^{2}\nu(dx_{i}), \mathbb{C}\right)\right)}\nonumber\\
\end{eqnarray}
\noindent goes to zero,  with \begin{eqnarray}&&\hspace*{-1.5cm}T\int_{[-\pi,\pi]^{3}}\mbox{cum}\left(\widetilde{X}_{\omega_{1}}^{(T)}, \widetilde{X}_{\omega_{2}}^{(T)},
\widetilde{X}_{\omega_{3}}^{(T)}, \widetilde{X}_{\omega_{4}}^{(T)}\right)\prod_{i=1}^{3}d\omega_{i}
\nonumber\\
&&\hspace*{1.5cm}=2\pi\int_{[-\pi,\pi]^{3}}\mathcal{F}_{\omega_{1},\omega_{2},\omega_{3}}
\prod_{j=1}^{3}d\omega_{j}+
\mathcal{O}(T^{-1}),\nonumber\end{eqnarray}
\noindent    in the norm of the space $\mathcal{S}\left(L^{2}\left(\mathbb{M}_{d}^{2},\otimes_{i=1}^{2}\nu(dx_{i}), \mathbb{C}\right)\right),$ where  $\mathcal{F}_{\omega_{1},\omega_{2},\omega_{3}}$ denotes the cumulant spectral  density operator of order $4$ of $X$ under $H_{1},$ introduced in equation (\ref{csdoo4}).
\end{proof}

\section{Proofs of the results in  Section~4}

\subsection{Proof of Proposition  \ref{cor1}}\label{apcb}
\begin{proof}  Under $H_{1},$ we have  $0<l_{\alpha }\leq \alpha (n)\leq L_{\alpha }<1/2,$ for every $n\in \mathbb{N}_{0}.$  From Lemma 3.1, considering  $T$ sufficiently large,
\begin{eqnarray}&&
\left|\int_{[-\sqrt{B_{T}}/2, \sqrt{B_{T}}/2]}E_{H_{1}}[\widehat{f}^{(T)}_{n}(\omega )]\frac{d\omega}{\sqrt{B_{T}}} \right|
\nonumber\\
&&
=\int_{-\pi}^{\pi}\frac{1}{B_{T}}\left[\int_{[-\sqrt{B_{T}}/2, \sqrt{B_{T}}/2]}W\left(\frac{\omega -\alpha }{B_{T}}\right)\frac{d\omega}{\sqrt{B_{T}}}
\right]f^{(T)}_{n}(\alpha )d\alpha \nonumber\\
&&+
\mathcal{O}(B_{T}^{-1}T^{-1})\nonumber\\
&&\simeq \frac{1}{\sqrt{B_{T}}}\int_{-\pi}^{\pi}\frac{1}{\sqrt{B_{T}}}W\left(\frac{ -\alpha }{B_{T}}\right)f_{n}(\alpha )d\alpha +
\mathcal{O}(B_{T}^{-1}T^{-1})+\mathcal{O}(T^{-1})\nonumber\\
&& \geq g(T)=\mathcal{O}(B_{T}^{-1/2-l_{\alpha }}),\quad \forall n\in \mathbb{N}_{0},
\label{eqfcpr1bh}
\end{eqnarray}
\noindent where  $\left\{ \widehat{f}^{(T)}_{n}(\omega ), \ n\in \mathbb{N}_{0}\right\}$  and $\left\{ f^{(T)}_{n}(\omega ) , \ n\in \mathbb{N}_{0}\right\}$  respectively denote  the frequency varying eigenvalues of the weighted periodogram operator $\widehat{\mathcal{F}}_{\omega }^{(T)}$ and the mean operator $\mathcal{F}_{\omega }^{(T)}=E\left[\mathcal{P}_{\omega }^{(T)}\right].$   Here,   $a_{T}\simeq b_{T}$ means that the two sequences $\{ a_{T},\ T>0\}$ and $\{b_{T},\ T>0\}$  have the same limit as $T\to \infty.$
From (\ref{eqfcpr1bh}),
\begin{eqnarray}
&&\left\| \int_{[-\sqrt{B_{T}}/2, \sqrt{B_{T}}/2]}E_{H_{1}}[\widehat{\mathcal{F}}_{\omega  }^{(T)}]\frac{d\omega}{\sqrt{B_{T}}}\right\|_{\mathcal{S}(L^{2}(\mathbb{M}_{d},d\nu, \mathbb{C}))}
\nonumber\\
&&
 \geq
\left\| \int_{[-\sqrt{B_{T}}/2, \sqrt{B_{T}}/2]}E_{H_{1}}[\widehat{\mathcal{F}}_{\omega  }^{(T)}]\frac{d\omega}{\sqrt{B_{T}}}\right\|_{\mathcal{L}(L^{2}(\mathbb{M}_{d},d\nu, \mathbb{C}))}
\nonumber\\
&&
= \sup_{n\in \mathbb{N}_{0}}\left|\int_{-\pi}^{\pi}\frac{1}{B_{T}}\left[\int_{[-\sqrt{B_{T}}/2, \sqrt{B_{T}}/2]}W\left(\frac{\omega -\alpha }{B_{T}}\right)\frac{d\omega}{\sqrt{B_{T}}}
\right]f^{(T)}_{n}(\alpha )d\alpha +
\mathcal{O}(B_{T}^{-1}T^{-1})\right|\nonumber\\
&&
\geq g(T)=\mathcal{O}(B_{T}^{-1/2-l_{\alpha }}),\quad T\to \infty,
\label{eqaltdiv}
\end{eqnarray}\noindent   where $\mathcal{L}(L^{2}(\mathbb{M}_{d},d\nu, \mathbb{C}))$ denotes the space of bounded linear operators on $L^{2}(\mathbb{M}_{d},d\nu, \mathbb{C}).$
\end{proof}
\subsection{Proof of Theorem \ref{consisimse}}
\label{appdd}
\begin{proof}
  From   Lemmas 3.1  and  3.3,   applying trace formula, for $T$ sufficiently large,
  \begin{eqnarray}
  &&\hspace*{-1.5cm}\int_{-\pi}^{\pi}E_{H_{1}}\left\|\widehat{\mathcal{F}}_{\omega }^{(T)}-E_{H_{1}}[\widehat{\mathcal{F}}_{\omega }^{(T)}]\right\|_{\mathcal{S}(L^{2}(\mathbb{M}_{d},d\nu, \mathbb{C}))}^{2}d\omega \nonumber\\
  &&\hspace*{-1.5cm}\leq \frac{2\pi}{TB_{T}}\int_{-\pi}^{\pi}  \sum_{n,h\in \mathbb{N}_{0}}\sum_{j=1}^{\Gamma (n,d)}\sum_{l=1}^{\Gamma (h,d)}\left| \int_{-\pi}^{\pi} \frac{1}{\sqrt{B_{T}}}W\left(\frac{\omega -\alpha }{B_{T}}\right)\frac{1}{\sqrt{B_{T}}} W\left(\frac{\omega -\alpha }{B_{T}}\right)
f_{n}(\alpha )f_{h}(\alpha)
d\alpha \right|d\omega \nonumber\\
&&\hspace*{-1.5cm}+\frac{2\pi}{TB_{T}}\int_{-\pi}^{\pi}\sum_{n,h\in \mathbb{N}_{0}}\sum_{j=1}^{\Gamma (n,d)}\sum_{l=1}^{\Gamma (h,d)}\left|\int_{-\pi}^{\pi} \frac{1}{\sqrt{B_{T}}}W\left(\frac{\omega -\alpha }{B_{T}}\right) \frac{1}{\sqrt{B_{T}}}W\left(\frac{\omega +\alpha }{B_{T}}\right)f_{n}(\alpha )f_{h}(\alpha)
d\alpha\right|d\omega\nonumber\\
&&\hspace*{-1.5cm}+\mathcal{O}(B_{T}^{-2}T^{-2})+\mathcal{O}(T^{-1})\nonumber\\
&&\simeq \frac{2\pi}{TB_{T}}\int_{-\pi}^{\pi}  \sum_{n,h\in \mathbb{N}_{0}}\sum_{j=1}^{\Gamma (n,d)}\sum_{l=1}^{\Gamma (h,d)}\left| \int_{-\pi}^{\pi} \delta_{T} (\omega-\alpha )
f_{n}(\alpha )f_{h}(\alpha)
d\alpha \right|d\omega \nonumber\\
&&\hspace*{-1.5cm}+\frac{2\pi}{TB_{T}}\int_{-\pi}^{\pi}\sum_{n,h\in \mathbb{N}_{0}}\sum_{j=1}^{\Gamma (n,d)}\sum_{l=1}^{\Gamma (h,d)}\left|\int_{-\pi}^{\pi}\delta_{T} (\omega-\alpha )f_{n}(\alpha )f_{h}(\alpha)
d\alpha\right|d\omega\nonumber\end{eqnarray}
 \begin{eqnarray}
&&+\mathcal{O}(B_{T}^{-2}T^{-2})+\mathcal{O}(T^{-1})\nonumber\\
&&\simeq  \frac{2\pi}{TB_{T}} \sum_{n,h\in \mathbb{N}_{0}}\sum_{j=1}^{\Gamma (n,d)}\sum_{l=1}^{\Gamma (h,d)}2\int_{[-\pi,\pi]}
f_{n}(\omega )f_{h}(\omega)
d\omega +\mathcal{O}(B_{T}^{-2}T^{-2})+\mathcal{O}(T^{-1})\nonumber\\
&&=  h(T)=\mathcal{O}(B_{T}^{-1}T^{-1}),\ T\to \infty,\nonumber\\
\label{eqpruebath2}
\end{eqnarray}
\noindent where, as before, $a_{T}\simeq b_{T}$ means that the two sequences $\{ a_{T},\ T>0\}$ and $\{b_{T},\ T>0\}$  have the same limit as $T\to \infty.$

\noindent Note  that,  under $H_{1},$ equation (\ref{eqpruebath2}) follows   from    condition $$\int_{[-\pi,\pi]}\|\mathcal{M}_{\omega }\|^{2}_{L^{1}(L^{2}(\mathbb{M}_{d},d\nu, \mathbb{C}))}|\omega |^{-2L_{\alpha }}d\omega <\infty,$$ \noindent  since
 \begin{eqnarray}
&&
\sum_{n,h\in \mathbb{N}_{0}}\sum_{j=1}^{\Gamma (n,d)}\sum_{l=1}^{\Gamma (h,d)}\int_{[-\pi,\pi]}
f_{n}(\omega )f_{h}(\omega)
d\omega \nonumber\\
&&
\leq \sum_{n,h\in \mathbb{N}_{0}}\sum_{j=1}^{\Gamma (n,d)}\sum_{l=1}^{\Gamma (h,d)}\int_{[-\pi,\pi]}M_{n}(\omega ) M_{h}(\omega )|\omega|^{-2L_{\alpha } }d\omega \nonumber\\
&&=
\int_{[-\pi,\pi]}\|\mathcal{M}_{\omega }\|_{L^{1}(L^{2}(\mathbb{M}_{d},d\nu,\mathbb{C}))}^{2}|\omega|^{-2L_{\alpha } }d\omega=\mathcal{O}(1).
\nonumber
 \end{eqnarray}

\end{proof}

\subsection{Proof of Corollary \ref{cor22}}
\label{apseccon}
\begin{proof}
Applying triangle and Jensen inequalities, we obtain from Corollary 3.2  and Theorem 4.2,
\begin{eqnarray}
&&\left\|\int_{-\pi}^{\pi}E_{H_{1}}\left[\widehat{\mathcal{F}}_{\omega }^{(T)}-\int_{-\pi}^{\pi}W(\xi)\mathcal{F}_{\omega -B_{T}\xi }d\xi \right]d\omega \right\|_{\mathcal{S}(L^{2}(\mathbb{M}_{d},d\nu, \mathbb{C}))}\nonumber\\
&&\leq \left\|\int_{-\pi}^{\pi}E_{H_{1}}\left[\widehat{\mathcal{F}}_{\omega }^{(T)}-E_{H_{1}}\left[\widehat{\mathcal{F}}_{\omega }^{(T)}\right]
\right]d\omega\right\|_{\mathcal{S}(L^{2}(\mathbb{M}_{d},d\nu, \mathbb{C}))}\nonumber\\
&&\hspace*{1.5cm}+\left\|\int_{-\pi}^{\pi}\left[E_{H_{1}}\left[\widehat{\mathcal{F}}_{\omega }^{(T)}\right]
-\int_{-\pi}^{\pi}W(\xi)\mathcal{F}_{\omega -B_{T}\xi }d\xi \right]d\omega \right\|_{\mathcal{S}(L^{2}(\mathbb{M}_{d},d\nu, \mathbb{C}))}\nonumber\\
&&\leq \left[\int_{-\pi}^{\pi}E_{H_{1}}\left\|\widehat{\mathcal{F}}_{\omega }^{(T)}-E_{H_{1}}\left[\widehat{\mathcal{F}}_{\omega }^{(T)}\right]\right\|_{\mathcal{S}(L^{2}(\mathbb{M}_{d},d\nu, \mathbb{C}))}^{2}
d\omega \right]^{1/2}\nonumber\\
&&\hspace*{1.5cm}+\left\|\int_{-\pi}^{\pi}\left[E_{H_{1}}\left[\widehat{\mathcal{F}}_{\omega }^{(T)}\right]
-\int_{-\pi}^{\pi}W(\xi)\mathcal{F}_{\omega -B_{T}\xi }d\xi \right] d\omega \right\|_{\mathcal{S}(L^{2}(\mathbb{M}_{d},d\nu, \mathbb{C}))}\nonumber\\
&&=\mathcal{O}(T^{-1/2}B_{T}^{-1/2}),\quad T\to \infty.
\label{imsecompcc}
    \end{eqnarray}

\end{proof}

\subsection{Proof of Theorem \ref{th2}}
\label{apd}
\begin{proof}
The proof of this result shares some ideas with the proof of  Theorem 2  of  \cite{Harrisetal08}, formulated in the time domain  for    real--valued time series.
Specifically, the test statistic operator $\mathcal{S}_{B_{T}}$ is  rewritten as
\begin{eqnarray}
&&\hspace*{-1cm}\mathcal{S}_{B_{T}}=\sqrt{B_{T}T}\int_{[-\sqrt{B_{T}}/2, \sqrt{B_{T}}/2]}E_{H_{1}}\left[\widehat{\mathcal{F}}_{\omega }^{(T)} \right]\frac{d\omega}{\sqrt{B_{T}}}
\nonumber\\
&& \circ
\left[\mathbb{I}_{L^{2}(\mathbb{M}_{d},d\nu, \mathbb{C})}+\left[\int_{[-\sqrt{B_{T}}/2, \sqrt{B_{T}}/2]}\left(\widehat{\mathcal{F}}_{\omega }^{(T)}-
E_{H_{1}}\left[\widehat{\mathcal{F}}_{\omega}^{(T)}\right]\right)\frac{d\omega}{\sqrt{B_{T}}}  \right]\right.
\nonumber\\
&&
\left.\circ \left[\int_{[-\sqrt{B_{T}}/2, \sqrt{B_{T}}/2]}E_{H_{1}}\left[\widehat{\mathcal{F}}_{\omega }^{(T)} \right]\frac{d\omega}{\sqrt{B_{T}}} \right]^{-1}\right],
\label{st}\end{eqnarray}
\noindent where $\circ$ means the composition of operators,
$\mathbb{I}_{L^{2}(\mathbb{M}_{d},d\nu, \mathbb{C})}$ denotes  the identity operator on the space
 $L^{2}(\mathbb{M}_{d},d\nu, \mathbb{C}),$  and $\left[\int_{[-\sqrt{B_{T}}/2, \sqrt{B_{T}}/2]}E_{H_{1}}\left[\widehat{\mathcal{F}}_{\omega}^{(T)}\right]\frac{d\omega}{\sqrt{B_{T}}} \right]^{-1}$ is  the inverse of operator
  $\int_{[-\sqrt{B_{T}}/2, \sqrt{B_{T}}/2]}E_{H_{1}}\left[\widehat{\mathcal{F}}_{\omega}^{(T)}\right]\frac{d\omega}{\sqrt{B_{T}}}.$

Our strategy  in the proof  of this result consists of   first proving, under $H_{1},$ the divergence, in the norm of the space $\mathcal{S}(L^{2}(\mathbb{M}_{d},d\nu, \mathbb{C})),$  of operator $$\sqrt{B_{T}T}\int_{[-\sqrt{B_{T}}/2, \sqrt{B_{T}}/2]}E_{H_{1}}\left[\widehat{\mathcal{F}}_{\omega }^{(T)} \right]\frac{d\omega}{\sqrt{B_{T}}}.$$
Then, under the conditions of Theorem 4.2,  we derive the convergence to zero, as $T\to \infty,$ of  random operator \begin{eqnarray}&&\hspace*{-1cm}\left[\int_{[-\sqrt{B_{T}}/2, \sqrt{B_{T}}/2]}\left[\widehat{\mathcal{F}}_{\omega}^{(T)}-
E_{H_{1}}\left[\widehat{\mathcal{F}}_{\omega}^{(T)}\right]\right]\frac{d\omega}{\sqrt{B_{T}}}\right] \circ \left[\int_{[-\sqrt{B_{T}}/2, \sqrt{B_{T}}/2]}E_{H_{1}}\left[\widehat{\mathcal{F}}_{\omega}^{(T)}\right]\frac{d\omega}{\sqrt{B_{T}}}  \right]^{-1},\nonumber\\
\label{eqrothm}\end{eqnarray}
\noindent   in the space $\mathcal{L}^{2}_{\mathcal{S}(L^{2}(\mathbb{M}_{d},d\nu, \mathbb{C}))}(\Omega ,\mathcal{A},P),$ which holds with a  suitable  rate under
 $l_{\alpha }>1/4$ and  $B_{T}=T^{-\beta },$ $\beta \in (0,1),$ allowing the application of Borell Cantelli Lemma to ensure almost surely convergence.
  Specifically,  from Proposition  4.1, as $T\to \infty,$
\begin{eqnarray}&&\left\|\sqrt{B_{T}T}\int_{[-\sqrt{B_{T}}/2, \sqrt{B_{T}}/2]}E_{H_{1}}\left[\widehat{\mathcal{F}}_{\omega }^{(T)}\right]\frac{d\omega}{\sqrt{B_{T}}} \right\|_{\mathcal{S}(L^{2}(\mathbb{M}_{d},d\nu, \mathbb{C}))}
\nonumber\\
&& \hspace*{4cm}\geq g(T)=\mathcal{O}\left(T^{1/2}B_{T}^{-l_{\alpha }}\right).\nonumber\\
\label{eqcor1}
\end{eqnarray}

 For the random operator in  (\ref{eqrothm}),   the following inequality holds:
\begin{eqnarray}
&&E_{H_{1}}\left\|\int_{[-\sqrt{B_{T}}/2, \sqrt{B_{T}}/2]}\left[\widehat{\mathcal{F}}_{\omega }^{(T)}-
E_{H_{1}}\left[\widehat{\mathcal{F}}_{\omega }^{(T)}\right]\right]\frac{d\omega}{\sqrt{B_{T}}}
\right.
\nonumber\\
&&
\hspace*{2cm}\circ \left.\left[\int_{[-\sqrt{B_{T}}/2, \sqrt{B_{T}}/2]}E_{H_{1}}\left[\widehat{\mathcal{F}}_{\omega}^{(T)}\right]\frac{d\omega}{\sqrt{B_{T}}}  \right]^{-1}
\right\|^{2}_{\mathcal{S}(L^{2}(\mathbb{M}_{d},d\nu, \mathbb{C}))}\nonumber\\
&&\leq \left\|
\left[\int_{[-\sqrt{B_{T}}/2, \sqrt{B_{T}}/2]}E_{H_{1}}\left[\widehat{\mathcal{F}}_{\omega}^{(T)}\right]\frac{d\omega}{\sqrt{B_{T}}} \right]^{-1} \right\|^{2}_{\mathcal{L}(L^{2}(\mathbb{M}_{d},d\nu, \mathbb{C}))}\nonumber\\
&&\hspace*{2cm} \times
E_{H_{1}}\left\|\int_{[-\sqrt{B_{T}}/2, \sqrt{B_{T}}/2]}\left[\widehat{\mathcal{F}}_{\omega}^{(T)}-
E_{H_{1}}\left[\widehat{\mathcal{F}}_{\omega}^{(T)}\right]\right]\frac{d\omega}{\sqrt{B_{T}}} \right\|^{2}_{\mathcal{S}(L^{2}(\mathbb{M}_{d},d\nu, \mathbb{C}))}.
\label{eqseconsumnum}
\end{eqnarray}
From equation  (\ref{eqaltdiv}),  as $T\to \infty,$\begin{eqnarray}&&\left\|\left[\int_{[-\sqrt{B_{T}}/2, \sqrt{B_{T}}/2]}E_{H_{1}}\left[\widehat{\mathcal{F}}_{\omega}^{(T)}\right]\frac{d\omega}{\sqrt{B_{T}}}\right]^{-1}
\right\|^{2}_{\mathcal{L}(L^{2}(\mathbb{M}_{d},d\nu, \mathbb{C}))}\leq b(T)=\mathcal{O}\left(B_{T}^{2l_{\alpha }+1}\right),\label{rest1}\nonumber\\
\end{eqnarray}
 \noindent  and, from  Theorem 4.2,
 \begin{eqnarray}&&\hspace*{-1cm}\int_{[-\sqrt{B_{T}}/2, \sqrt{B_{T}}/2]}\mbox{Var}_{H_{1}}\left(\widehat{\mathcal{F}}_{\omega }^{(T)}\right)\frac{d\omega}{\sqrt{B_{T}}}
\nonumber\\
&&\hspace*{-0.5cm}
\leq\int_{[-\pi,\pi]}\mbox{Var}_{H_{1}}\left(\widehat{\mathcal{F}}_{\omega }^{(T)}\right)\frac{d\omega}{\sqrt{B_{T}}} \leq u(T)= \mathcal{O}(T^{-1}B_{T}^{-1-1/2}),\quad   T\to \infty.
\label{rest2}
\end{eqnarray}

For each $T\geq 2,$ applying  Jensen inequality, in terms of  the uniform probability measure on the interval $[-\sqrt{B_{T}}/2, \sqrt{B_{T}}/2],$
\begin{eqnarray}&&
E_{H_{1}}\left\|\int_{[-\sqrt{B_{T}}/2, \sqrt{B_{T}}/2]}\left[\widehat{\mathcal{F}}_{\omega}^{(T)}-
E_{H_{1}}\left[\widehat{\mathcal{F}}_{\omega}^{(T)}\right]\right]\frac{d\omega}{\sqrt{B_{T}}} \right\|^{2}_{\mathcal{S}(L^{2}(\mathbb{M}_{d},d\nu, \mathbb{C}))}\nonumber\\
&&=\left\|\int_{[-\sqrt{B_{T}}/2, \sqrt{B_{T}}/2]}\left[\widehat{\mathcal{F}}_{\omega}^{(T)}-
E_{H_{1}}\left[\widehat{\mathcal{F}}_{\omega}^{(T)}\right]\right]\frac{d\omega}{\sqrt{B_{T}}} \right\|_{
\mathcal{L}^{2}_{\mathcal{S}(L^{2}(\mathbb{M}_{d},d\nu, \mathbb{C}))}(\Omega ,\mathcal{A},P)}^{2}\nonumber\\
&&=\varphi_{H_{1}}\left( E_{\mathcal{U}([-\sqrt{B_{T}}/2, \sqrt{B_{T}}/2])}\left[\widehat{\mathcal{F}}_{\omega}^{(T)}-
E_{H_{1}}\left[\widehat{\mathcal{F}}_{\omega}^{(T)}\right]\right]\right)\nonumber\\
&&\leq E_{\mathcal{U}([-\sqrt{B_{T}}/2, \sqrt{B_{T}}/2])}\left[\varphi_{H_{1}}\left(\left[\widehat{\mathcal{F}}_{\omega}^{(T)}-
E_{H_{1}}\left[\widehat{\mathcal{F}}_{\omega}^{(T)}\right]\right]\right)\right]\nonumber\\
&&= \int_{[-\sqrt{B_{T}}/2, \sqrt{B_{T}}/2]}E_{H_{1}}\left\|\widehat{\mathcal{F}}_{\omega}^{(T)}-
E_{H_{1}}\left[\widehat{\mathcal{F}}_{\omega}^{(T)}\right]\right\|^{2}_{\mathcal{S}(L^{2}(\mathbb{M}_{d},d\nu, \mathbb{C}))}\frac{d\omega}{\sqrt{B_{T}}},\label{rest2bb}
\end{eqnarray}
\noindent where  $E_{\mathcal{U}([-\sqrt{B_{T}}/2, \sqrt{B_{T}}/2])}$ denotes expectation under the uniform probability measure on the interval $[-\sqrt{B_{T}}/2, \sqrt{B_{T}}/2],$  and
$\varphi_{H_{1}}(\cdot ) =\|\cdot\|_{\mathcal{L}^{2}_{\mathcal{S}(L^{2}(\mathbb{M}_{d},d\nu, \mathbb{C}))}(\Omega ,\mathcal{A},P)}^{2}= E_{H_{1}}\|\cdot \|^{2}_{\mathcal{S}(L^{2}(\mathbb{M}_{d},d\nu, \mathbb{C}))}$ is a convex function. Thus,  from equations (\ref{eqseconsumnum})--(\ref{rest2bb}),    \begin{eqnarray}&&\hspace*{-2cm} E_{H_{1}}\left\|\int_{[-\sqrt{B_{T}}/2, \sqrt{B_{T}}/2]}\left[\widehat{\mathcal{F}}_{\omega}^{(T)}-
E_{H_{1}}\left[\widehat{\mathcal{F}}_{\omega}^{(T)}\right]\right]\frac{d\omega}{\sqrt{B_{T}}} \right.\nonumber\\
&&\left.\hspace*{0.5cm}\circ \left[\int_{[-\sqrt{B_{T}}/2, \sqrt{B_{T}}/2]}E_{H_{1}}\left[\widehat{\mathcal{F}}_{\omega}^{(T)}\right]\frac{d\omega}{\sqrt{B_{T}}}  \right]^{-1}
\right\|^{2}_{\mathcal{S}(L^{2}(\mathbb{M}_{d},d\nu, \mathbb{C}))}
\nonumber\\
&& \hspace*{2cm} \leq h(T)=\mathcal{O}(T^{-1}B_{T}^{2l_{\alpha }-1/2}),\quad  T\to \infty.
\label{sicv}
\end{eqnarray}
 From  equation (\ref{sicv}), applying Chebyshev's inequality,
 \begin{eqnarray}
 &&
 P\left[\left\|\int_{[-\sqrt{B_{T}}/2, \sqrt{B_{T}}/2]}\left[\widehat{\mathcal{F}}_{\omega}^{(T)}-
E_{H_{1}}\left[\widehat{\mathcal{F}}_{\omega}^{(T)}\right]\right]\frac{d\omega}{\sqrt{B_{T}}} \right.\right.\nonumber\\
&&\left.\left.\hspace*{0.5cm}\circ \left[\int_{[-\sqrt{B_{T}}/2, \sqrt{B_{T}}/2]}E_{H_{1}}\left[\widehat{\mathcal{F}}_{\omega}^{(T)}\right]\frac{d\omega}{\sqrt{B_{T}}}  \right]^{-1}\right\|_{\mathcal{S}(L^{2}(\mathbb{M}_{d},d\nu, \mathbb{C}))}>\varepsilon\right]\nonumber\\
&&\leq E_{H_{1}}\left\|\int_{[-\sqrt{B_{T}}/2, \sqrt{B_{T}}/2]}\left[\widehat{\mathcal{F}}_{\omega}^{(T)}-
E_{H_{1}}\left[\widehat{\mathcal{F}}_{\omega}^{(T)}\right]\right]\frac{d\omega}{\sqrt{B_{T}}} \right.\nonumber\\
&&\left.\hspace*{0.5cm}\circ \left[\int_{[-\sqrt{B_{T}}/2, \sqrt{B_{T}}/2]}E_{H_{1}}\left[\widehat{\mathcal{F}}_{\omega}^{(T)}\right]\frac{d\omega}{\sqrt{B_{T}}}  \right]^{-1}
\right\|^{2}_{\mathcal{S}(L^{2}(\mathbb{M}_{d},d\nu, \mathbb{C}))}/\varepsilon^{2}\nonumber\\
&&\leq h(T)/\varepsilon^{2}=\mathcal{O}(T^{-1}B_{T}^{2l_{\alpha }-1/2}).
\label{feq}
\end{eqnarray}
Since   $l_{\alpha }>1/4,$  hence, $2l_{\alpha }-1/2=\rho>0,$   and,  for $B_{T}=T^{-\beta },$
$T^{-1}B_{T}^{2l_{\alpha }-1/2}=T^{-1-\beta \rho},$ with $\beta \in (0,1),$ and $\rho \in (0,1/2).$ From equation  (\ref{feq}),
 Borel--Cantelli  lemma then  leads to
\begin{eqnarray}
&&\left\|\int_{[-\sqrt{B_{T}}/2, \sqrt{B_{T}}/2]}\left[\widehat{\mathcal{F}}_{\omega}^{(T)}-
E_{H_{1}}\left[\widehat{\mathcal{F}}_{\omega}^{(T)}\right]\right]\frac{d\omega}{\sqrt{B_{T}}}\right. \nonumber\\
&&\hspace*{2cm}\left. \circ \left[\int_{[-\sqrt{B_{T}}/2, \sqrt{B_{T}}/2]}E_{H_{1}}\left[\widehat{\mathcal{F}}_{\omega}^{(T)}\right]\frac{d\omega}{\sqrt{B_{T}}}  \right]^{-1}\right\|_{\mathcal{S}(L^{2}(\mathbb{M}_{d},d\nu, \mathbb{C}))}\to_{a.s.} 0.\nonumber\\
\label{feqfv}
\end{eqnarray}
\noindent  as  $T\to \infty.$
  The a.s. divergence of $\left\|\mathcal{S}_{B_{T}}\right\|_{\mathcal{S}(L^{2}(\mathbb{M}_{d},d\nu, \mathbb{C}))},$ as  $T\to \infty ,$ follows   from   equations  (\ref{st}), (\ref{eqcor1}) and (\ref{feqfv}).

\end{proof}

\subsection{Proof of Lemma \ref{CLTSBT}}
\label{apendicea0}
\begin{proof}
 For $\omega \in (-\pi,\pi )\backslash \{0\},$  consider a Gaussian random  element  $\widehat{\mathcal{F}}_{\omega }$  in the space $\mathcal{S}(L^{2}(\mathbb{M}_{d}, d\nu ,\mathbb{C})),$ with kernel $\widehat{f}_{\omega }$  satisfying  (\ref{cobk}). Hence, from (\ref{sdo33}),
 \begin{eqnarray}
&&\frac{1}{2\pi\|W\|^{2}_{L^{2}(\mathbb{R})}}E\left[\widehat{f}_{\omega }\otimes \widehat{f}_{\omega }\right](\tau_{1} ,\sigma_{1},\tau_{2},\sigma_{2})\nonumber\\
&&=\sum_{n,h\in \mathbb{N}_{0}}\sum_{j=1}^{\Gamma (n,d)}\sum_{l=1}^{\Gamma (h,d)}
f_{n}(\omega )f_{h}(\omega ) S_{n,j}^{d}(\tau_{1})\overline{S_{h,l}^{d}(\sigma_{1})S_{n,j}^{d}(\tau_{2})}S_{h,l}^{d}(\sigma_{2}),
 \label{eqacobb2}\end{eqnarray}
\noindent for every $(\tau_{i},\sigma_{i})\in \mathbb{M}_{d}^{2},$ $i=1,2,$  $\omega \in (-\pi,\pi )\backslash \{0\}.$
 Thus, the diagonal coefficients \linebreak $\left\{\lambda_{n,h}(\omega ),\ n,h\in \mathbb{N}_{0}\right\}=\left\{ f_{n}(\omega )f_{h}(\omega ),\ n,h\in \mathbb{N}_{0}\right\}$ define   the eigenvalues of the autocovariance operator  (\ref{eqacobb2}).  Since $\mathbb{M}^{2}_{d}$ is a compact set, and $\mathcal{R}_{\widehat{\mathcal{F}}_{\omega } ,\widehat{\mathcal{F}}_{\omega }}$ is a trace  positive   semidefinite self--adjoint operator
 under $H_{0},$
the orthogonal expansion \begin{eqnarray}
&&\hspace*{-1cm}\frac{1}{\sqrt{2\pi}\|W\|_{L^{2}(\mathbb{R})}}\widehat{f}_{\omega }(\tau,\sigma)
=\sum_{n,h\in \mathbb{N}_{0}}\sum_{j=1}^{\Gamma (n,d)}
\sum_{l=1}^{\Gamma (h,d)}
\sqrt{f_{n}(\omega)f_{h}(\omega)}
Y_{n,j,h,l}(\omega)S_{n,j}^{d}(\tau)\overline{S_{h,l}^{d}(\sigma)}\nonumber\\
\label{KLSHb}\end{eqnarray}\noindent  holds in the space $\mathcal{L}^{2}_{\mathcal{S}(L^{2}(\mathbb{M}_{d},d\nu; \mathbb{C}))}(\Omega ,\mathcal{A},\mathcal{P}).$    The random Fourier coefficients are given by \begin{eqnarray}&&
Y_{n,j,h,l}(\omega )=\frac{(\sqrt{2\pi}\|W\|_{L^{2}(\mathbb{R})})^{-1}}{\sqrt{f_{n}(\omega  )f_{h}(\omega  )}}\int_{\mathbb{M}_{d}^{2}}\widehat{f}_{\omega }(\tau ,\sigma )\overline{S_{n,j}^{d}}(\tau)S_{h,l}^{d}(\sigma)d\nu (\sigma )d\nu (\tau ),
\nonumber\\
&& j=1,\dots, \Gamma (n,d),\  l=1,\dots, \Gamma(h,d), \ n,h\in \mathbb{N}_{0},\ \omega \in [-\pi,\pi ]\backslash \{0\}.
\label{KLSH2bb}\end{eqnarray}

\end{proof}

\end{appendix}

\end{document}